\documentclass[12pt,english,british]{paper}
\usepackage{mathptmx}
\usepackage[T1]{fontenc}
\usepackage[latin9]{inputenc}
\usepackage[a4paper]{geometry}
\usepackage{float}
\usepackage{amsmath}
\usepackage{graphicx}
\usepackage{amssymb}

\makeatletter

\newcommand{\noun}[1]{\textsc{#1}}

\usepackage{amsmath,amssymb}
 \usepackage{centernot}

\makeatother

\usepackage{babel}

\begin{document}
\[
\mbox{\LARGE\textsc{flows near compact invariant sets}}\]
\[
\mbox{\large Part I}\]
\[
\mbox{\large\textsc{pedro teixeira}}\]
\[
10\mbox{\,\,\ February\,\,}2012\]

\begin{abstract}
In this paper it is proved that near a compact, invariant, proper
subset of a ~$C^{\,0}$ ~flow on a locally compact, connected metric
space, at least one, out of twenty eight relevant dynamical phenomena,
will necessarily occur. ~This result (Theorem 1) shows that assuming
the connectedness of the phase space, implies the existence of a considerably
deeper classification of topological flow behaviour in the vicinity
of compact invariant sets than that described in the classical theorems
of Ura-Kimura and Bhatia. ~The proposed classification brings to
light, in a systematic way, the possibility of occurrence of \emph{orbits
of infinite height} arbitrarily near the compact invariant in question,
and this under relatively simple conditions. ~Singularities of ~$C^{\,\infty}$
~vector fields displaying this strange phenomenon occur in every
dimension ~$n\geq3$ (in this paper, a ~$C^{\,\infty}$ ~flow on
~$\mathbb{S}^{3}$ ~exhibiting such an equilibrium is constructed).
~Near periodic orbits, the same phenomenon is observable already
in dimension 4 ~(and on every manifold of dimension ~$n\geq5)$.~
As a corollary to the main result, an elegant characterization of
the topological Hausdorff structure of the set of all compact minimal
sets of the flow is obtained (Theorem 2). \medskip{}

\emph{keywords: }topological behaviour of ~$C^{\,0}$ ~flows,
compact invariant sets, compact minimal sets, topological Hausdorff
structure, non-hyperbolic singularities and periodic orbits, orbits
of infinite height.

\emph{\hfill{}In memoriam Vladimir I. Arnol'd}

\bigskip{}

\begin{figure}[H]
\noindent \begin{centering}
\includegraphics[scale=0.95]{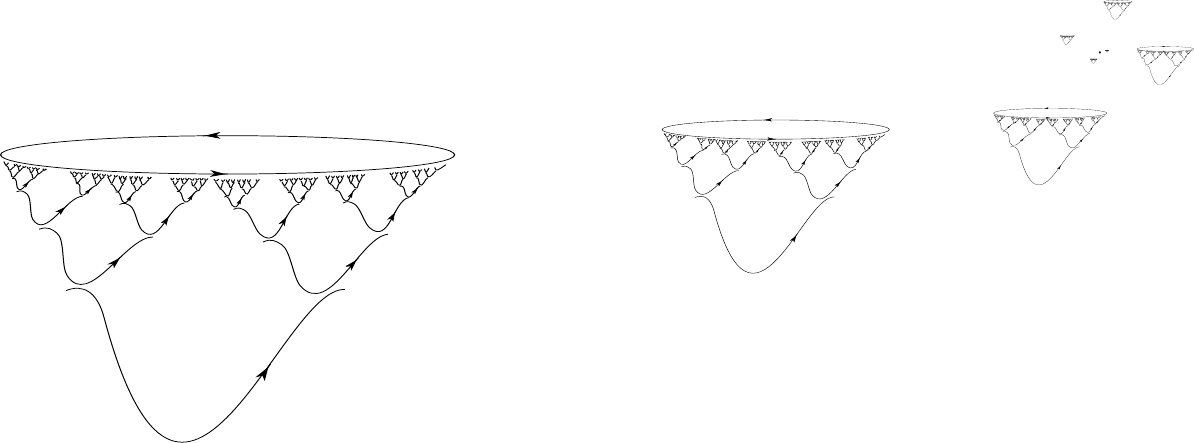}
\par\end{centering}

\caption{{\small Example of the occurrence of orbits of }\emph{\small infinite
height}{\small . This figure illustrates case 8.6 of Theorem 1.}}
\end{figure}

\pagebreak{}

\noun{\normalsize CONTENTS.}{\normalsize }\\
{\normalsize 1.~~Introduction.}\\
{\normalsize 2.~~Definitions and basic results.}\\
{\normalsize 3.~~Presentation of the main result. ~Corollaries.}\\
{\normalsize 4.~~Special orbital structures.}\\
{\normalsize 5.~~Lemmas.}\\
{\normalsize 6.~~The main theorem. ~Topological structure of
~$\mbox{\ensuremath{\mbox{CMin}(M)}.}$}\\
{\normalsize 7.~~Independent realizations. ~Examples.}\\
{\normalsize 8.~~Proof of lemma 7. ~Two topological lemmas.}{\normalsize \par}

\medskip{}

1.\noun{ INTRODUCTION.}

\medskip{}

{\normalsize ~~~The present work establishes a natural classification
of topological behaviour of $C^{\,0}$ ~flows near arbitrary compact
invariant sets ~$K$, ~on locally compact connected metric spaces
~$M$ ~(}\emph{\normalsize e.g.}{\normalsize{} on 2nd countable,
Hausdorff, compact or not, connected manifolds). ~It can be seen
as a considerably deep generalization of a classical topological-dynamical
result of Ura and Kimura \cite{ura} and Bhatia (see }\emph{\normalsize e.g.
}{\normalsize \cite{bhat}, p. 114), when the hypothesis of connectedness
of the phase space is added. ~The }\emph{\normalsize Ura-Kimura-Bhatia
Theorem}{\normalsize{} states that if ~$M$ ~is a locally compact
(but not necessarily connected) metric space and ~$K$ ~is as above,
then at least one of the following four cases takes place:}{\normalsize \par}
\end{abstract}
I.~~~~$K$ ~is an attractor (\emph{i.e.} asymptotically stable)

II.~~~$K$ ~is a repeller (\emph{i.e. }negatively asymptotically
stable)

III. ~there are points ~$x,y\in M\setminus K$ ~such that ~$\emptyset\neq\alpha(x)\subset K$
~and ~$\emptyset\neq\omega(y)\subset K$

IV.~~every neighbourhood of ~$K$ ~contains a compact invariant
set that contains\\
\hspace*{6mm}$K$ ~as a proper subset (\emph{i.e. }~$K$~
is not an isolated invariant set).\medskip{}

~~~While valid for very general flows$\big/$phase spaces and despite
its importance, namely, in persistence theory (see \emph{e.g.} the
preface of \cite{smit}), the above result has, in our opinion, an
obvious serious limitation that hinders the possibility of a natural,
substantial deepening of the classification it proposes: since the
phase space is not assumed to be connected, a (nonvoid) compact invariant
set ~$Q\varsubsetneq M$ ~may be open in ~$M$. ~This makes ~$Q$
~simultaneously an attractor and a repeller, while in fact ~$Q$
~neither attracts nor repels a single point outside itself. ~Actually
as ~$M\setminus Q$ ~is closed, sufficiently near but outside ~$Q$
~the flow is vacuous!

~~Adding the assumption of connectedness of the phase space dramatically
improves the possibility of partially describing the {}``dynamical
landscape'' in the vicinity of a compact invariant set. ~Natural
considerations lead to the identification of twenty five possible
relevant dynamical phenomena that fall under case IV of the \emph{Ura-Kimura-Bhatia
Theorem}. ~Moreover the whole twenty eight cases are distributed
among five groups, two cases belonging to distinct groups being incompatible
\emph{i.e }cannot be simultaneously satisfied. A key role in the classification
is played by compact invariant sets ~$\emptyset\neq K\subsetneq M$
~that are either attractors or repellers or isolated from minimal
sets and stagnant. ~By the later we mean that for some neighbourhood
~$U$ ~of~ $K$, ~$U\setminus K$ ~contains no minimal set of
the flow and in addition condition ~III~ above is satisfied. ~Although
the main result of this paper (Theorem 1) goes much deeper, a {}``flavour''
of some of its most important conclusions is given in the following
corollary:\medskip{}

\noindent \emph{Let} ~$M$ ~\emph{be a locally compact, connected
metric space with a} ~$C^{\,0}$\emph{~flow} ~$\theta$\emph{ ~and}
~$K$~\emph{ a compact, invariant, proper subset of} ~$M$. \emph{~Then:}

\emph{either~~}$\,$\textbf{I})~~~\emph{~K is an attractor}

\emph{or~~~~~~~~}\textbf{II})\emph{~~~~K is a repeller}

\emph{or}~~~~~~~$\,$$\,$\textbf{III})~~~~\emph{K is isolated
from minimal sets and stagnant.}

\emph{or}~~~~~~\textbf{IV})\emph{~~~~there is a nonvoid,
compact, connected invariant set ~}$Q\subset\mbox{bd\,}K$\emph{
~and a sequence ~}$\varLambda_{_{n}}\subset M\setminus K$\emph{
~of compact minimal sets of the flow such that the following three
conditions hold:}

~~$\bullet$~~~$(\varLambda_{_{n}})$ ~~\emph{converges to
~$Q$ ~in the Hausdorff metric}

~~$\bullet$~~~\emph{all ~$\varLambda_{_{n}}'s$ ~belong to
the same one of the following three classes: ~equilibrium orbits,
~periodic orbits, ~aperiodic compact minimals sets}

~~$\bullet$~~~\emph{either all ~$\varLambda_{_{n}}'s$ ~are
attractors or they are repellers or they are all isolated from minimal
sets and stagnant.\medskip{}
}

\emph{or}~~\textbf{V})~~~\emph{for all sufficiently small open
neighbourhood ~$U$ ~of ~$K$, ~the compact minimals sets contained
in ~$U\setminus K$ ~form a nonvoid ~$\mathfrak{c}$-dense in itself
set i.e. any neighbourhood of a compact minimal set ~$\varLambda\subset U\setminus K$
~contains a continuum of compact minimal sets.}%
\footnote{note that ~$M$ ~being locally compact, every sufficiently small
neighbourhood\emph{ }of ~$K$ ~may contain only \emph{compact} minimal
sets, hence in both ~IV~and ~V~ we may actually replace compact
minimal set(s) by minimal set(s) everywhere. %
}

\emph{Finally if none of conditions }~I\emph{ ~to ~}V\emph{~ holds,
then} \emph{orbits of infinite height}%
\footnote{An orbit ~$\gamma_{_{0}}$ ~is of ~\emph{infinite height }if there
is an infinite strict inclusion chain of orbit closures\[
\mbox{cl\,}\gamma_{_{0}}\supsetneq\mbox{cl\,}\gamma_{_{1}}\supsetneq\cdots\cdots\supsetneq\mbox{cl\,}\gamma_{_{n}}\supsetneq\cdots\cdots\]
(some authors call ~$\gamma_{_{0}}$ ~an orbit of \emph{infinite
depth}).%
}\emph{ will necessarily occur arbitrarily near but outside ~$K$,
~more precisely,}

~~~~~\textbf{VI})\emph{~~~~Given any neighbourhood ~$U$
~of ~$K$, ~there is a sequence of orbits ~$\gamma_{_{n}}\subset U\setminus K$~
such that}\[
\mbox{cl}\,\gamma_{_{1}}\supsetneq\mbox{cl}\,\gamma_{_{2}}\supsetneq\cdots\cdots\supsetneq\mbox{cl}\,\gamma_{_{n}}\supsetneq\cdots\cdots\]

\begin{figure}[H]
\noindent \begin{centering}
\includegraphics[scale=0.26]{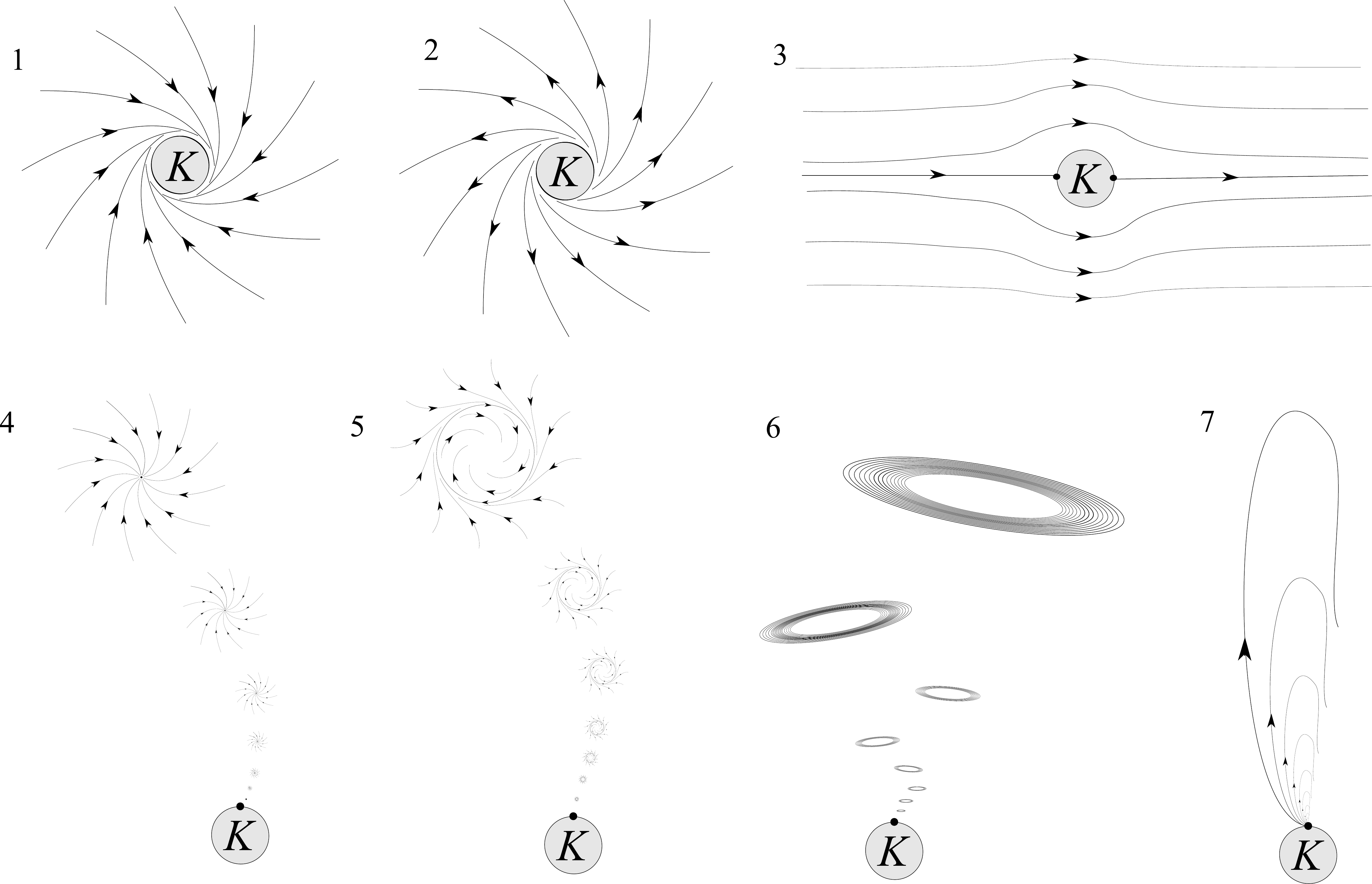}
\par\end{centering}

\caption{{\small Examples: 1 - case I; 2 - case II; 3 - case III; 4\&5 - case
IV; 6 - case V; 7 - case VI.}}
\end{figure}

~~~Therefore, if ~$K$ ~is neither an attractor nor a repeller
and conditions ~III~ and ~IV~ fail, then a {}``super-abundance''
of compact minimal sets (case V) or an outstanding kind of limit behaviour
(case VI) will emerge arbitrarily near (but outside) ~$K$. ~The
possible occurrence of the later disturbing dynamical phenomenon is
not a mere theoretical speculation: in section 7 ~a ~$C^{\,\infty}$
~flow on ~$\mathbb{S}^{3}$ ~exhibiting it is given. ~This is
made possible by the existence of a ~$C^{\,\infty}$ ~flow without
minimal sets on a non-compact surface of infinite genus ~$(C^{\,\infty}\mbox{-})$
embeddable in ~$\mathbb{R}^{3}$ ~(see Beniere and Meigniez \cite{beni}
and the pioneer work of Inaba \cite{inab}). ~In a subsequent paper
\cite{teix} we shall actually show that our classification is both
pertinent and non-redundant: each of the twenty eight cases it describes
admits an independent realization either by a ~$C^{\,\infty}$ ~flow
on a (2nd countable, Hausdorff) compact, connected manifold or by
a $(C^{\,0})$ ~subflow corresponding to a compact, connected invariant
subset of such ~$C^{\,\infty}$~ flow (see section 7). 

~~~Among its many interesting consequences, the classification
theorem%
\footnote{see section 3 for a first illustrated presentation and section 6 for
the full statement.%
} has remarkable, and somewhat unexpected, implications in the topological
structure of the set ~$\mbox{CMin}(M)$ ~of all \emph{compact }minimal
sets of the flow, endowed with the Hausdorff metric ~$d_{_{H}}$
~(see sections 3 and 6). ~As an example, let ~$\mathfrak{A}$ ~be
the set of compact minimal sets that are either attractors or repellers
or isolated from minimal sets and stagnant. ~Then,\pagebreak{}

~~~{}``\emph{If a ~$C^{\,0}$ ~flow on a locally compact, connected
metric space} $M$ \emph{~has only countably many compact minimal
sets and displays no orbits of infinite height then ~$\mathfrak{A}$
~is necessarily open dense in ~}$\mbox{CMin}(M)$\emph{ ~(in the
Hausdorff metric)''.}

If the flow displays uncountably many compact minimal sets, then a
preliminary result (lemma 7) permits to establish a topological decomposition
of ~$\mbox{CMin}(M)$, ~in a certain sense analogue to that of \emph{Cantor-Bendixson
Theorem} for \emph{Polish spaces}:

~~~{}``\emph{If ~}$\mbox{CMin}(M)$ \emph{~is uncountable then
all but a countable number of compact minimal sets of the flow have
a continuum of compact minimal sets on each of their neighbourhoods''. }

This result brings to light the corresponding analogue proposition
~$\varSigma$ ~for the set ~$\mbox{Per}(M)$ ~of all periodic
orbits of the flow:

~~~{}``\emph{If ~}$\mbox{Per}(M)$ \emph{~is uncountable then
all but a countable number of periodic orbits of the flow have a continuum
of periodic orbits on each of their neighbourhoods''. }

There are no counterexamples to ~$\varSigma$ ~within {}``standard''
Dynamical Systems Theory: ~its negation implies the \emph{Negation
of the Continuum Hypothesis ~}(this is actually a purely topological
fact, depending only on the Hausdorff metric separability of ~$\mbox{Per}(M)\,\big)$.\emph{
~}Hence ~$\varSigma$ ~is either demonstrable in ~\emph{Zermelo-Fraenkel+Axiom
of Choice (ZFC) ~}or it is independent of this standard axiomatic
(see section 6).

~~~Despite its topological nature, the greatest interest of Theorem
1 lies in the context of $C^{\, r\geq0}$ ~flows on ~$C^{\, r}$
manifolds.%
\footnote{except if otherwise mentioned, manifolds are always assumed to be
2nd countable, Hausdorff, compact or not, connected and boundaryless.%
}~ Particularly noteworthy is, perhaps, its contribution to the understanding
of what can happen, from the dynamical point of view, in two potential
{}``nightmare'' phenomena of differentiable dynamics: \emph{non-hyperbolic
singularities and periodic orbits}.%
\footnote{in relation to the later, the seven cases ~7.1 ~to ~7.6 ~and ~10.1~
of Theorem 1 (section 6) are \emph{apriori }ruled out (since singularities
cannot accumulate on regular points$\big)$.%
} ~To see how hopeless standard {}``analytic'' methods may be in
the study of the former, even in low dimensions, consider, for example,
the case of complete ~$C^{\,\infty}$ vector fields on ~$\mathbb{R}^{2}$
~having the origin ~$O$ ~as an isolated \emph{flat}%
\footnote{the point ~$O$ ~is a \emph{flat }singularity of ~$X\in\mathfrak{X}^{\infty}(\mathbb{R}^{2})$
~if this vector field vanishes together with its derivatives of all
orders at ~$O$, ~\emph{i.e.} if it has identically null Taylor
expansion at\emph{ }that point.%
}\emph{ }singularity. ~It is not difficult to see that there is a
\emph{continuum }of such vector fields ~$X_{_{i}}$, ~$i\in\mathbb{R}$,
~that are mutually\emph{ }topologically\emph{ non }equivalent%
\footnote{recall that two complete vector fields ~$X,\, Y\in\mathfrak{X}^{\,\infty}(\mathbb{R}^{2})$
\emph{~}are\emph{ topologically equivalent at} ~$O$ ~if there
are open neighbourhoods ~$U$ ~and ~$V$ ~of ~$O$ ~and a homeomorphism
~$\varphi:U\longrightarrow V$, ~fixing ~$O$~ and ~carrying
each maximal segment of $\, X$-orbit contained in ~$U$~ onto a
maximal segment of $\, Y$-orbit contained in ~$V$, ~preserving
time orientation.%
} at\emph{ ~$O$, }~and whose local topological behaviour at ~$O$
~cannot (with the eventual exception of some very general dynamical
properties such as stability) be investigated by standard differential
methods $\big(\, C^{\, r\geq1}$ ~coordinate changes, blow-up desingularizations,
etc$\big)$. ~This shows that already in ~$\mathbb{R}^{2}$, ~there
is a \emph{continuum }of distinct possible topological ~$C^{\,\infty}$
~flow behaviours near an isolated singularity ~$O$, ~that are
practically left in the darkness by {}``analytic'' methods, and
in such cases there seems to be no much alternative to what can be
learned from the topological-dynamical approach.

~~~This paper is organized as follows: after the preliminaries
of section 2 (giving special emphasis to Hausdorff metric concepts),
a provisional version of the main result (Theorem 1) is presented
in section 3. ~Instead of the twenty eight cases of the full statement,
only seventeen cases are distinguished. ~The last of these, condition
H (occurrence of orbits of \emph{infinite height }arbitrarily near
but outside the compact invariant in question)%
\footnote{see condition VI on p. 3.%
} ~encompasses 12 distinct cases of the full statement, requiring
the introduction of certain denumerable collections of orbits displaying
a kind of \emph{{}``}fractal-like\emph{''} structure with respect
to orbital limit relations. ~These collections ($X\mbox{-}trees$,
~$X\mbox{-}\alpha shells$ ~and ~$X\mbox{-\ensuremath{\omega}}shells$),
capturing essential features of the dynamical complexity of this strange
and beautiful phenomenon, are fully presented in section 4.~ Still
on section 3, some important remarks concerning the meaning and scope
of application of Theorem 1 (particularly in the context of flows
on manifolds), as well as some interesting consequences for the topological
structure of the set of all compact minimal sets of the flow, are
given. ~Twelve lemmas, some of them dynamically significant on their
own, are presented in section 5. ~Lemma 2 is a version of the above
mentioned Ura-Kimura-Bhatia Theorem, a proof of it being included
for the sake of completeness. ~The main result of this paper, Theorem
1, is fully stated in section 6. ~As a principal corollary, we obtain
a simple characterization of the topological Hausdorff structure of
the set of all compact minimal sets of the flow ~$\mbox{CMin}(M)$
~(Theorem 2). ~The question of existence of flows displaying, near
a compact invariant set, the types of dynamical behaviour foreseen
in Theorem 1 is addressed in section 7. ~Examples are actually constructed
showing the possibility of independent realization of some of the
cases susceptible of raising greater doubts, namely those involving
the occurrence of orbits of \emph{infinite height}. ~As already mentioned,
a positive answer to the question of existence of independent realizations
for all the twenty eight cases will be given in \cite{teix} under
quite favourable smoothness conditions. ~Finally a proof of Lemma
7 is presented in section 8. ~This result, crucial for the present
work, shows that in a certain sense, the {}``local cardinality''
of an open and dense in itself (with respect to the Hausdorff metric)
nonvoid set of compact minimal sets behaves as that of (nonvoid) perfect
subsets of ~$\mathbb{R}.$

\pagebreak{}

\noun{2.~DEFINITIONS AND BASIC RESULTS.}\medskip{}

~~~Let ~$M$ ~be a metric space with a ~$C^{\,0}$ ~flow ~$\theta:\mathbb{R}\times M\longrightarrow M$
~and ~$K$~ a compact, invariant, proper subset of~ $M$ ~$\big($we
denote {}``the flow ~$\theta$ ~on ~\emph{$M$'' ~}by\emph{
}~$(M,\theta)\big).$ ~A \emph{minimal set }of ~$(M,\theta)$ ~is
a nonvoid, closed, invariant subset of ~$M$~ that contains no proper
subset satisfying these three conditions \emph{i.e.} an orbit closure
that contains no smaller one. ~We reserve the expression ~\emph{periodic
orbit} ~to orbits ~$\mathcal{O}(x)$~ for which ~$\big\{ t\in\mathbb{R}:\,\theta(t,x)=x\big\}=\lambda\mathbf{\mathbb{Z}}$,~
for some ~$\lambda>0.$ ~In this case, the unique ~$\lambda>0$~~satisfying
that identity is the ~\emph{period of}$\,\,$ $\mathcal{O}(x)$~~(of
the \emph{periodic point }~$x).$ ~A minimal set that is neither
an \emph{equilibrium} \emph{orbit}%
\footnote{the orbit of an\emph{ equilibrium point z ~i.e }a singleton ~$\big\{ z\big\}=\big\{\theta(t,z):\, t\in\mathbb{R\big\}}.$%
} nor a \emph{periodic orbit} is called an~ \emph{aperiodic minimal.}

\textbf{Definitions.} ~Let ~$M$,~ $\theta$, ~$K$ ~as above
and~ $x\in M$, ~~$X\subset M$.\medskip{}

\selectlanguage{english}%
~~$\mathcal{N}_{_{X}}:=\mbox{\mbox{the set of neighborhoods of \,\ensuremath{X\mbox{ \,\ in \,}M}.}}$\foreignlanguage{british}{\medskip{}
}

\selectlanguage{british}%
$\begin{array}{ll}
\mathcal{O}(x):=\{\theta(t,x):\, t\in\mathbb{R}\} & \mbox{the orbit of \,\ensuremath{x}}\\
\mathcal{O}^{+}(x):=\{\theta(t,x):\,\, t\geq0\} & \mbox{the positive (half) orbit of \,\ensuremath{x}}\\
\mathcal{O}(X):=\overset{\,}{\underset{\underset{\,}{x\in X}}{\bigcup}\mathcal{O}(x)} & \mbox{the orbital saturation of \,\ensuremath{X}}\\
\mathcal{O}^{+}(X):=\underset{\underset{\,}{x\in X}}{\bigcup}\mathcal{O}^{+}(x) & \mbox{the positive orbital saturation of \,\ensuremath{X}}\\
\mbox{Orb}(X):=\mathcal{\big\{ O}(x):\,\mathcal{O}(x)\subset X\big\}\, & \mbox{the set of orbits contained in \,\ensuremath{X}}\end{array}$\medskip{}

$\mathcal{O}^{-}(x)$ ~and ~$\mathcal{O}^{-}(X)$ ~are the negative
concepts corresponding to ~$\mathcal{O}^{+}(x)$ ~and ~$\mathcal{O}^{+}(X)$.
~When dealing with a unique flow ~$\theta$ ~we write ~$x_{_{t}}$
~for ~$\theta(t,x).$\medskip{}

\hspace{3mm}$\omega(x):=\underset{t>0}{\bigcap}\mathcal{\mbox{cl\,}O}^{+}(x_{_{t}})$
~~the ~$\omega$-limit set of ~$x$.

\hspace{3mm}$\alpha(x):=\underset{t<0}{\bigcap}\mbox{cl\,}\mathcal{O}^{-}(x_{_{t}})$
~~the ~$\alpha$-limit set of ~$x$.\medskip{}

For any orbit ~$\gamma=\mathcal{O}(x)$, ~we define ~$\alpha(\gamma):=\alpha(x)$,
~$\omega(\gamma):=\omega(x)$ ~since points in the same orbit have
the same $\alpha\mbox{-limit}$ ~and ~$\omega\mbox{-limit}$ sets.\medskip{}

\hspace{3mm}$B^{+}(K):=\big\{\, x\in M:\,\emptyset\,\,\neq\,\,\omega(x)\subset K\,\big\}$

\hspace{3mm}$B^{-}(K):=\big\{\, x\in M:\,\emptyset\,\,\neq\,\,\alpha(x)\subset K\,\big\}$

\hspace{3mm}$A^{+}(K):=\big\{\, x\in M:\,\emptyset\,\,\neq\,\,\omega(x)\,\cap\, K\,\,\neq\,\,\omega(x)\,\big\}$

\hspace{3mm}$A^{-}(K):=\big\{\, x\in M:\,\emptyset\,\,\neq\,\,\alpha(x)\,\cap\, K\neq\,\,\alpha(x)\,\big\}$\medskip{}

\emph{i.e}. ~$A^{+}(K)$ ~$\big(\, A^{-}(K)\,\big)$ ~is the set
of points of ~$M$ ~whose ~$\omega$-limit $\big($$\alpha$-limit$\big)$
~set intercepts both ~$K$ ~and ~$M\setminus K$. ~We say ~$K$
~is:\medskip{}

$\bullet$\hspace{3mm}\emph{stable} ~if for any ~$U\in\mathcal{N}_{_{K}}$
~there is a~ $V\in\mathcal{N}_{_{K}}$~ such that ~$\mathcal{O}^{+}(V)\subset U$.

$\bullet$\hspace{3mm}\emph{negatively stable} if for any ~$U\in\mathcal{N}_{_{K}}$
there is a ~$V\in\mathcal{N}_{_{K}}$ such that ~$\mathcal{O}^{-}(V)\subset U$.

$\bullet$\hspace{3mm}\emph{bi-stable} ~if it is both stable and
negatively stable.

$\bullet$\hspace{3mm}\emph{bi-stable in relation to }~$N\subset M$
if for any ~$U\in\mathcal{N}_{_{K}}$~ there is a ~$V\in\mathcal{N}_{_{K}}$~
such that ~$\mathcal{O}(N\,\cap\, V)\subset U$ ~\emph{i.e. }given
~$U\in\mathcal{N}_{_{K}}$, ~any point ~$x\in N$ ~at a sufficiently
small distance from ~$K$ ~has its orbit entirely contained in ~$U$.

$\bullet$\hspace{3mm}an ~\emph{attractor}%
\footnote{If ~$K$ ~is an \emph{attractor }then it is easily seen that ~$B^{+}(K)$
~is an invariant, open subset of ~$M$. ~Idem for ~$B^{-}(K)$
~if ~$K$ ~is a \emph{repeller.}%
} ~if it~is stable and ~$B^{+}(K)\in\mathcal{N}_{_{K}}$.

$\bullet$\hspace{3mm}a ~\emph{repeller} ~if it negatively stable
and ~$B^{-}(K)\in\mathcal{N}_{_{K}}$ ~i.e if it is an \emph{attractor
}in the time reversal flow \foreignlanguage{english}{~$\phi(t,x)=\theta(-t,x)$.}

$\bullet$\hspace{3mm}\emph{stagnant}%
\footnote{we shall often use the expression \emph{non-stagnant }instead of \emph{not
stagnant.}%
} if there are points ~$x,y\in M\setminus K$ ~such that ~$\emptyset\neq\alpha(x)\subset K$
~and ~$\emptyset\neq\omega(y)\subset K$.

$\bullet$\hspace{3mm}\emph{isolated from minimal sets }~if there
is a ~$U\in\mathcal{N}_{_{K}}$ ~such that ~$U\setminus K$ ~contains
no minimal set of ~$(M,\theta).$ ~We also use the abridged terminology
\emph{isolated from minimals}. ~If ~$K$ ~is itself a minimal set
then we say that ~$K$ ~is an \emph{~isolated minimal (set}).

\smallskip{}

\textbf{Definitions.} 

$\begin{array}{lll}
\mbox{C}(X)\,\,\, & := & \mbox{the set of nonvoid, compact subsets of \,}X.\\
\mbox{Ci}(X)\, & := & \mbox{the set of nonvoid, compact, invariant subsets of \,\emph{X}.}\\
\mathfrak{S}(X) & := & \mbox{the set of nonvoid, compact, connected, invariant subsets of \emph{X}.}\\
\mbox{CMin}(X) & := & \mbox{the set of compact minimal sets contained in \,\ensuremath{X}}.\\
\mbox{Eq}(X) & := & \big\{\{x\}:\, x\in M\,\,\,\,\,\mbox{and}\,\,\,\,\,\theta(t,x)=x\mbox{\,\,\ for all \,}t\in\mathbb{R}\big\}=\\
 & \,\,= & \mbox{the set of equilibrium orbits contained in \,}X.\\
\mbox{Per}(X) & := & \mbox{the set of periodic orbits contained in \,\ensuremath{X}}.\\
\mbox{Am}(X) & := & \mbox{the set of compact aperiodic minimal sets contained in \ensuremath{X}.}\end{array}$

\medskip{}

$\mbox{C}(M)$ ~and its subsets are naturally endowed with the \emph{Hausdorff
metric} $d_{_{H}}.$~ To emphasise that this\emph{ }metric is the
one in question, we employ the expressions ~$d_{_{H}}-$\emph{open},
~$d_{_{H}}-$\emph{closed}, ~$d_{_{H}}-$\emph{near,} ~$d_{_{H}}-$\emph{converges}
$\big(\overset{d_{_{H}}}{\longrightarrow}\big)$, ~$d_{_{H}}-$\emph{isolated},
~etc. %
\footnote{Metric concepts in ~$\mbox{\ensuremath{\big(}C}(M),d_{_{H}}\big)$
are distinguished from the corresponding concepts in ~$(M,d)$~
by the subscript ~$_{H}$ ~ $\mbox{e.g.. \,\,}B_{_{H}}(\, X,\,\epsilon\,):=\big\{ Y\in\mbox{C}(M):\, d_{_{H}}(X,Y)<\epsilon\big\};$
~analogously, closure, boundary and interior are denoted by ~$\mbox{cl}_{_{H}},$
$\mbox{bd}_{_{H}},$ $\mbox{int}_{_{H}}$.%
}~ A set ~$\mathfrak{A}\subset\mbox{C}(M)$ ~~$d_{_{H}}-$\emph{accumulates
in ~}$\mathfrak{B}\subset\mbox{C}(M)$ ~if ~$(\mbox{cl}_{_{H}}\mathfrak{A})\,\cap\,\mathfrak{B}\neq\emptyset$.
~A sequence ~$X_{_{n}}\in\mbox{C}(M)$ ~~$d_{_{H}}-$\emph{accumulates
in ~}$\mathfrak{B}\subset\mbox{C}(M)$ ~if it has a subsequence
~$d_{_{H}}-$converging to some ~$X\in\mathfrak{B}.$ ~Working
primarily within the Hausdorff metric, we shall deal essentially with
\emph{equilibrium orbits} rather than with \emph{equilibria. ~}Note,
however, that the set ~$E$ ~of \emph{equilibria} of the flow, endowed
with the metric ~$d$ ~of ~$M$, ~is isometric to the metric space
~$\big[\,\mbox{Eq}(M),\, d_{_{H}}\big]$ ~via the canonical map
~$e\longrightarrow\{e\}$. ~The following classical result, originally
proved in the context of convex body theory, is of central importance
to the present work: 

\medskip{}

\textbf{\emph{Blaschke Theorem}}\emph{: If ~$[\, N,d\,]$ ~is a
compact metric space then so is ~}$\mbox{\ensuremath{\big[}C}(N),d_{_{H}}\big]$\emph{.}\medskip{}

where\emph{ ~}$\mbox{C}(N)$ \emph{~}is the set of nonvoid, compact
subsets of\emph{ ~$N$ ~}(see \emph{e.g.} \cite{bura},\emph{ }p.253)\emph{.}
~If ~$N$~ is a compact metric space and ~$\mathfrak{C}$ ~is
a $d_{_{H}}-$closed (and thus compact) subset of ~$\mbox{C}(N),$
~then the consequent possibility of selecting from a given sequence
~$\varLambda_{_{n}}\in\mathfrak{C}$, ~a subsequence $d_{_{H}}-$converging
to some ~$\varLambda\in\mathfrak{C}$, ~will be referred as ~\emph{Blaschke
Principle.~} To avoid double indices, we will often suppose that
the selected subsequence is ~$(\varLambda_{_{n}})$ ~itself.~ Again
in the metric space ~\emph{M}, ~if ~$N\subset M$ ~is compact
then the continuity of the flow implies that ~$\mbox{Ci}(N)$ ~is
$d_{_{H}}-$closed in ~$\mbox{C}(N)$ ~and thus compact; a simple
argument%
\footnote{Suppose ~$Y_{_{n}}\in\mbox{Cc}(N)$ ~and ~$Y_{_{n}}\longrightarrow X$.
~If ~$X\in\mbox{C}(N)\setminus\mbox{Cc}(N)$ ~then $X=X_{_{_{0}}}\,\cup\, X_{_{_{1}}}$
~where ~$X_{_{_{0}}},\,\,\, X_{_{1}}$ ~are disjoint, nonvoid compacts.
~Letting ~$\lambda:=\mbox{dist}$$\big(X_{_{_{0}}},\, X_{_{_{1}}}\big)/2$,
~because of its connectedness, each ~$Y_{_{n}}$ ~has at least
one point ~$z_{_{n}}$ ~in the compact ~$S(X_{_{_{0}}},\lambda):=\big\{ y\in M:\,\mbox{dist}(y,X_{_{_{0}}})=\lambda\big\}$.
~Taking a convergent subsequence ~$z_{_{n_{j}}}\longrightarrow z\in\mbox{ \ensuremath{S(X_{_{0}},\lambda)}}$
~it follows that ~$z\in\mbox{lim}\, Y_{_{n}}=X=X_{_{_{0}}}\,\cup\, X_{_{_{1}}}$
which is absurd in virtue of the definition of ~$\lambda$ ~$\big(\,$as
~$S(X_{_{0}},\lambda)\,\cap\, X{}_{_{0}}=\emptyset=S(X_{_{0}},\lambda)\,\cap\, X_{_{1}}\,\big)$.%
} shows that ~$\mbox{Cc}(N),$ ~the set of nonvoid, compact, connected
subsets of ~$N$ ~is also compact, hence ~$\mathfrak{S}(N)=\mbox{Ci}(N)\,\cap\,\mbox{Cc}(N)$
~is compact\emph{.} ~Observe that while ~$\mbox{Ci}(N),$ $\mbox{Cc}(N),$
$\mathfrak{S}(N)$ ~and ~$\mbox{Eq}(N)$ ~are ~$d_{_{H}}-$closed
in~ $\mbox{C}(N)$ ~and thus compact, ~$\mbox{CMin}(N),$ ~$\mbox{Per}(N)$
~and ~$\mbox{Am}(N)$ ~in general are not. ~Note that ~$\mbox{CMin}(N)=\mbox{Eq}(N)\,\sqcup\,\mbox{Per}(N)\,\sqcup\,\mbox{Am}(N)\subset\mathfrak{S}(N)\subset\mbox{C}(N)$
~$\big(\,\sqcup$ ~denotes disjoint union$\big)$.

\medskip{}

\emph{Remark}. ~The reader should keep in mind the following basic
facts as they will often be implicitly used without mention.~ Suppose
~$N\subset M$ ~is compact. ~If ~$\mathcal{O}^{+}(x)\subset N$~
then ~$\omega(x)$ ~and ~$\mbox{cl}\,\mathcal{O}^{+}(x)=\mathcal{O}^{+}(x)\,\cup\,\omega(x)$
~both belong to ~$\mathfrak{S}(N)$ ~and in particular are nonvoid.
~Analogue fact holds for ~$\mathcal{O}^{-}(x)$, ~$\alpha(x)$
~and ~$\mbox{cl}\,\mathcal{O}^{-}(x)$ ~when ~$\mathcal{O}^{-}(x)\subset N$.
~Also ~$\gamma\in\mbox{Orb}(N)$ ~implies ~$\mbox{cl\,}\gamma=\gamma\,\cup\,\alpha(\gamma)\,\cup\,\omega(x)\in\mathfrak{S}(N)$.
~If ~$N$ ~is a nonvoid, compact invariant set then it contains
at least one compact minimal set of the flow. ~If ~$X$ ~is a\emph{
}minimal set\emph{ }and ~$K$~ is a\emph{ }closed invariant set,
then either ~$X\subset K$ ~or ~$X\subset M\setminus K$, ~since
the set of closed invariant sets is closed under interceptions. ~If
~$N\subset M$ ~is invariant then ~$\mbox{cl\,}N$, ~$\mbox{bd\,}N$
~and ~$\mbox{int\,}N$ ~(respectively, topological closure, boundary
and interior of ~$N$) are also invariant.

\medskip{}

\textbf{Definitions. ~}A set ~$\mathfrak{C}\subset\mbox{C}(M)$
~is ~$d_{_{H}}-$\emph{dense in itself} ~if every ~$\varLambda\in\mathfrak{C}$
~is not ~$d_{_{H}}\mbox{-}$ \emph{isolated}~in ~$\mathfrak{C}$
~\emph{i.e.} if~$\varLambda\in\mbox{cl}_{_{H}}\big(\mathfrak{C}\setminus\{\varLambda\}\big)$,
~for all ~$\varLambda\in\mathfrak{C}$. ~$\varLambda\in\mbox{C}(M)$
~is a ~$\mathfrak{c}-$\emph{condensation element} \emph{of }~$\mathfrak{C}$~
if for every~ $\epsilon>0$\[
\#\,\big(B_{_{H}}(\varLambda,\epsilon)\,\cap\,\mathfrak{C}\big)\geq\mathfrak{c}\]
where ~$\mathfrak{c}$ ~denotes the cardinal of the \emph{continuum}.
~A set ~$\mathfrak{C}\subset\mbox{C}(M)$ ~is ~$\mathfrak{c}-$\emph{dense
in itself} ~\emph{ }if every~ $\varLambda\in\mathfrak{C}$ ~is
a ~$\mathfrak{c}-$condensation element\emph{ of} ~$\mathfrak{C}.$

(note that in this paper, except if otherwise mentioned, the expression~
\emph{$\mathfrak{c}-$dense in itself }~always respects to the Hausdorff
metric ~$d_{_{H}}$ ~and the same applies to \emph{$\mathfrak{c}-$condensation
element})\emph{.}

\emph{Important remark. }If ~$M$ ~is a locally compact, connected
metric space, then ~$M$ ~is necessarily \emph{separable} (see \emph{e.g.}
\cite{lima}, p.278) and thus has at most a \emph{continuum} of points.%
\footnote{Let ~$N\subset M$ ~be a countable dense subset. ~We can associate
with each ~$z\in M$ ~a distinct sequence ~$x_{_{n}}\in N$ ~converging
to ~$z$. ~Since the cardinal of the set of all sequences of points
of ~$N$ ~is at most ~$\#\big(\mathbb{N}^{\mathbb{N}}\big)=\#\mathbb{R}=\mathfrak{c}$,
~it follows that ~$\#M\leq\mathfrak{c}$ ~\emph{i.e. ~$M$} ~has,
at most, a \emph{continuum }of points. %
}~ Therefore there is, at most, a \emph{continuum }of orbits in the
flow ~$(M,\theta)$ ~and also, at most, a \emph{continuum }of minimal
sets (distinct minimal sets are disjoint), thus if ~$\mathfrak{C}\subset\mbox{CMin}(M)$
~then the inequality above reduces to\[
\#\,\big(B_{_{H}}(\varLambda,\epsilon)\,\cap\,\mathfrak{C}\big)=\mathfrak{c}\]
Hence for such phase space ~$M$, ~a ~$\mathfrak{c}$-\emph{dense
in itself} ~set of compact minimal sets is either empty or has the
cardinal of the \emph{continuum}. ~In chapter 5 (Lemma 6) we shall
actually see that a set of compact minimal sets ~$\mathfrak{C}\subset\mbox{CMin}(M)$
~is ~$\mathfrak{c}$-\emph{dense in itself} ~iff every neighbourhood
~$U_{_{\varLambda}}\subset M$ ~of each ~$\varLambda\in\mathfrak{C}$
~contains a \emph{continuum} of elements of ~$\mathfrak{C}$, ~showing
that in this particularly important case, we may actually think in
terms of the simpler metric ~$d$ ~of ~$M$, ~instead of the Hausdorff
metric ~$d_{_{H}}$ ~of ~$\mbox{C}(M)$.\medskip{}

\textbf{Definition.} ~For each ~$\mathfrak{C}\subset2^{M}\,(=$
the set of subsets of ~$M\,$)\emph{~}~and ~$A\subset M,$

\hspace{3mm}$\mathfrak{C}^{*}:=\bigcup\mathfrak{C}=\underset{\varGamma\in\mathfrak{C}}{\bigcup}\,\,\varGamma$~~~(the
set of points of \emph{~M~} belonging to elements of \emph{~$\mathfrak{C}$}).

\hspace{3mm}$\mathfrak{C}(A):=\{X\in\mathfrak{C}:\, X\subset A\}$~~~(the
set of elements of ~$\mathfrak{C}$ ~contained in ~$A\,$).\medskip{}

\textbf{Definition.} (metric concept on ~$M$) ~Given any two nonvoid
sets ~$X,\, Y\subset M$,\medskip{}

\hspace{3mm}\foreignlanguage{english}{$\big|Y\big|_{X}:=\mbox{sup}\big\{\mbox{dist}(y,X):\,\, y\in Y\big\}\in[0,+\infty]$}\medskip{}

If we naturally identify ~$B(X,+\infty)$ ~with ~$M$, ~then \medskip{}

\hspace{3mm}$\big|Y\big|_{X}=\mbox{inf}\big\{\delta\in[0,+\infty]:\,\, Y\subset B(X,\delta)\big\}$

\bigskip{}
\medskip{}

\noun{3.~PRESENTATION OF THE MAIN RESULT.}\bigskip{}

\noindent ~~~We now give the reader an approximate idea of the
main result of this article (Theorem$\,$1) which is fully stated
and proved in section 6. ~Roughly speaking, it shows that if ~$K$~
is a compact, invariant, proper subset of a ~$C^{\,0}$ ~flow on
a locally compact, connected metric space ~$M$,~ then at least
one, out of twenty eight relevant dynamical phenomena, will necessarily
occur near ~$K$. ~\emph{The reader is warned that the result} \emph{presented
bellow is weaker than Theorem$\,$1}. ~However, this is the best
that can be done without the introduction of certain special denumerable
collections of orbits (called \emph{X}-\emph{trees}, \emph{X}-\emph{$\alpha\,$shells}
and \emph{X}-\emph{$\omega\,$shells}) that will be presented in the
next section. ~In order to simplify the exposition, we shall distribute
the 28 cases among seventeen conditions. ~The 16 cases corresponding
to conditions ~\textbf{A}~ to ~\textbf{G.4}~ will be presented
in a near definitive form. ~The remaining 12 cases, corresponding
to condition ~\textbf{H}, require the introduction of the aforementioned
collections of orbits and are postponed until section 6. ~Notwithstanding,
we shall give bellow an idea of the kind of dynamical phenomenon that
is necessarily encountered if conditions ~\textbf{A}~ to ~\textbf{G.4}~
fail. ~Here we will put emphasis on a {}``geometric flavoured''
and rather descriptive enunciation, in detriment of concision and
elegance. ~Whenever it is possible, Hausdorff ~$d_{_{H}}$ ~metric
properties are reformulated in terms of the simpler and more intuitive
metric ~$d$ ~of ~$M$~ (using lemma 6 of section 5). ~This actually
results in a longer statement than that of the stronger version of
section 6. ~There, using appropriate Hausdorff metric terminology,
a quite sharp, but perhaps at a first sight less illuminating, presentation
is achieved.

\noindent %
\begin{figure}[H]
\noindent \begin{centering}
\includegraphics[scale=1.1]{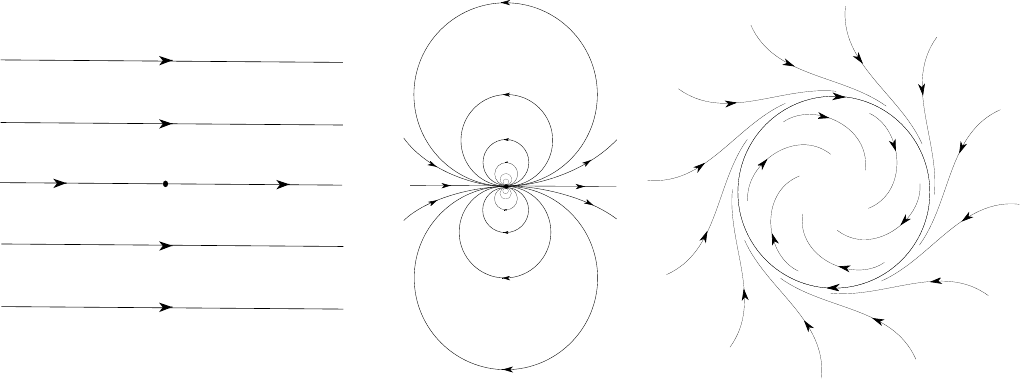}
\par\end{centering}

\caption{{\small Examples of isolated from minimals and stagnant compact invariant
sets on ~$\mathbb{S}^{2}$.~ (Left) fake saddle equilibrium orbit.
~(Right) periodic orbit attracting on one side and repelling on the
other (periodic orbit of saddle-node type).}}
\end{figure}
~~~Let ~$M$ ~be a \emph{locally compact, connected }metric space
with a ~$C^{\,0}$~flow ~$\theta$. ~Consider the following three
propositions where the variable ~$X$ ~assumes values in the set
~$\mathfrak{\mbox{Ci}(}M)$ ~of nonvoid, compact, invariant subsets
of ~$M$:\medskip{}

\hspace{3mm}1\emph{.X} \hspace{6mm}$X$~ \emph{is an} \emph{attractor}

\hspace{3mm}2\emph{.X}\hspace{6mm}$X$~ \emph{is a} \emph{repeller}

\hspace{3mm}3\emph{.X}\hspace{6mm}$X$~ \emph{is} \emph{isolated
from minimals and} \emph{stagnant}\medskip{}

~~~Observe that, due to the \emph{connectedness} of ~$M$, each
compact, invariant,\emph{ proper} subset of ~$M$ ~may satisfy,
at most, one of the above conditions: ~as ~$M$ ~is\emph{ }connected,
if ~$X\in\mbox{Ci}(M)\setminus\{M\}$ ~is an attractor\emph{ }then~
$B^{+}(X)\setminus X\neq\emptyset$ ~since ~$X$ ~is a closed proper
subset of ~$M$ ~and ~$B^{+}(X)$ ~is an open neighbourhood of
~$X$. ~Therefore there are actually points ~$x\in M\setminus X$
~such that ~$\omega(x)\subset X.$ ~An analogue argument shows
that if ~$X$ ~is a repeller then there are points ~$y\in M\setminus X$
~such that ~$\alpha(y)\subset X$. ~Thus an\emph{ }attractor\emph{
}~$X$ ~can neither be a repeller\emph{ }nor stagnant\emph{ }since
both these conditions contradict the \emph{stability }of ~$X.$~
Analogue fact holds if ~$X$ ~is a repeller.\emph{ ~Isolated from
minimals and stagnant} compact, invariant sets play an important role
in the present work. ~In Differentiable Dynamics, typical, dynamically
distinct examples are given by\emph{ }saddle, fake saddle and saddle-node
equilibrium orbits and by\emph{ }hyperbolic saddle periodic orbits\emph{.}%
\footnote{\emph{~}Another instructive example is given by the orbit of the
equilibrium\emph{ }~$(0,\ldots,0,1)\in\mathbb{S}^{n}\subset\mathbb{R}^{n+1}$
~in the compactification of the flow on ~$\mathbb{S}^{\, n}\setminus\{(0,\ldots,0,1)\}$
~induced, via the inverse stereographic projection, by the constant
vector field ~$\frac{\partial}{\partial x_{_{1}}}$ ~on ~$\mathbb{R}^{\, n}$
~(see fig.3, centre, for the case ~$n=2$). ~A more subtle example
is given by the exceptional minimal set\emph{ }on the celebrated Denjoy\emph{
}~$C^{\,1}$ ~flow on ~$\mathbb{T}^{\,2}$.%
}

\bigskip{}

\noindent Weak (provisional) version of \textbf{Theorem 1}:\medskip{}

\noindent \emph{Let} ~$M$ ~\emph{be a locally compact, connected
metric space with a} ~$C^{\,0}$\emph{~flow} ~$\theta$\emph{ ~and}
~$K$~\emph{ a compact, invariant, proper subset of} ~$M$. \emph{~Then:}\medskip{}

\noindent \emph{either}

\noindent \textbf{\emph{A.}}\textbf{\hspace*{0.2in}}$K$~ \emph{is
an attractor}

\emph{or}

\textbf{\emph{B.}}\textbf{\hspace*{0.2in}}$K$~\emph{ is a} \emph{repeller}

\emph{or}

\textbf{\emph{C.}}\textbf{\hspace*{0.2in}}$K$~ \emph{is} \emph{isolated
from minimals and} \emph{stagnant\bigskip{}
}%
\begin{figure}[H]
\noindent \begin{centering}
\includegraphics[scale=0.45]{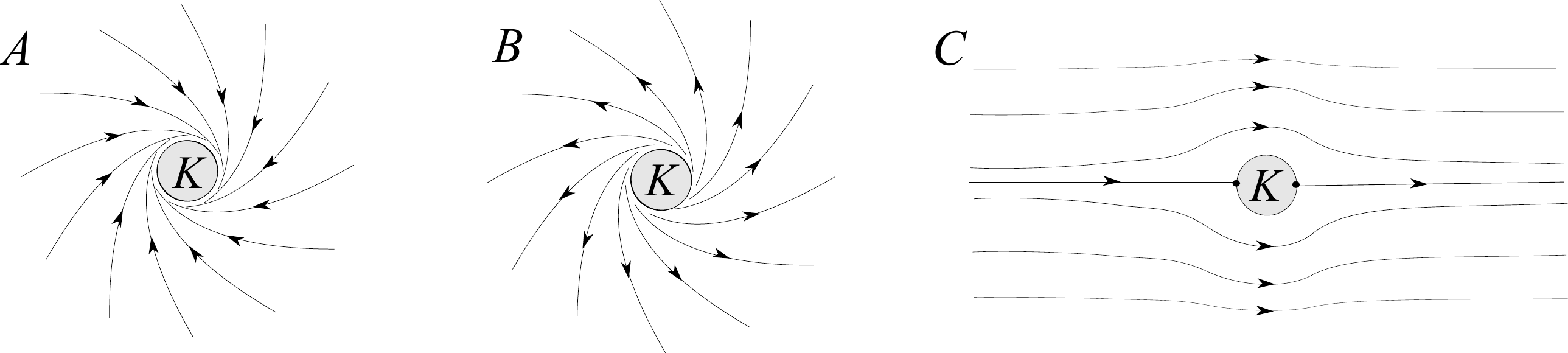}
\par\end{centering}

\caption{}
\end{figure}

\medskip{}

\pagebreak{}

\emph{or at least one of the following nine conditions holds:}

\textbf{\emph{D}}\textbf{.}\textbf{\emph{i}}\textbf{\hspace*{0.2in}}$1\leq$$i\leq3$\\
\emph{there is a sequence of equilibria} ~$e_{_{n}}\in M\setminus K$,\emph{
~converging to an equilibrium ~}$q\in\mbox{bd\,}K$\emph{ ~and}
\emph{such that}\textbf{ }\emph{condition} ~\emph{$i.X$ }~\emph{is
satisfied by all equilibrium orbits ~$\{e_{_{n}}\}$.}

\begin{figure}[H]
\noindent \begin{centering}
\includegraphics[scale=0.21]{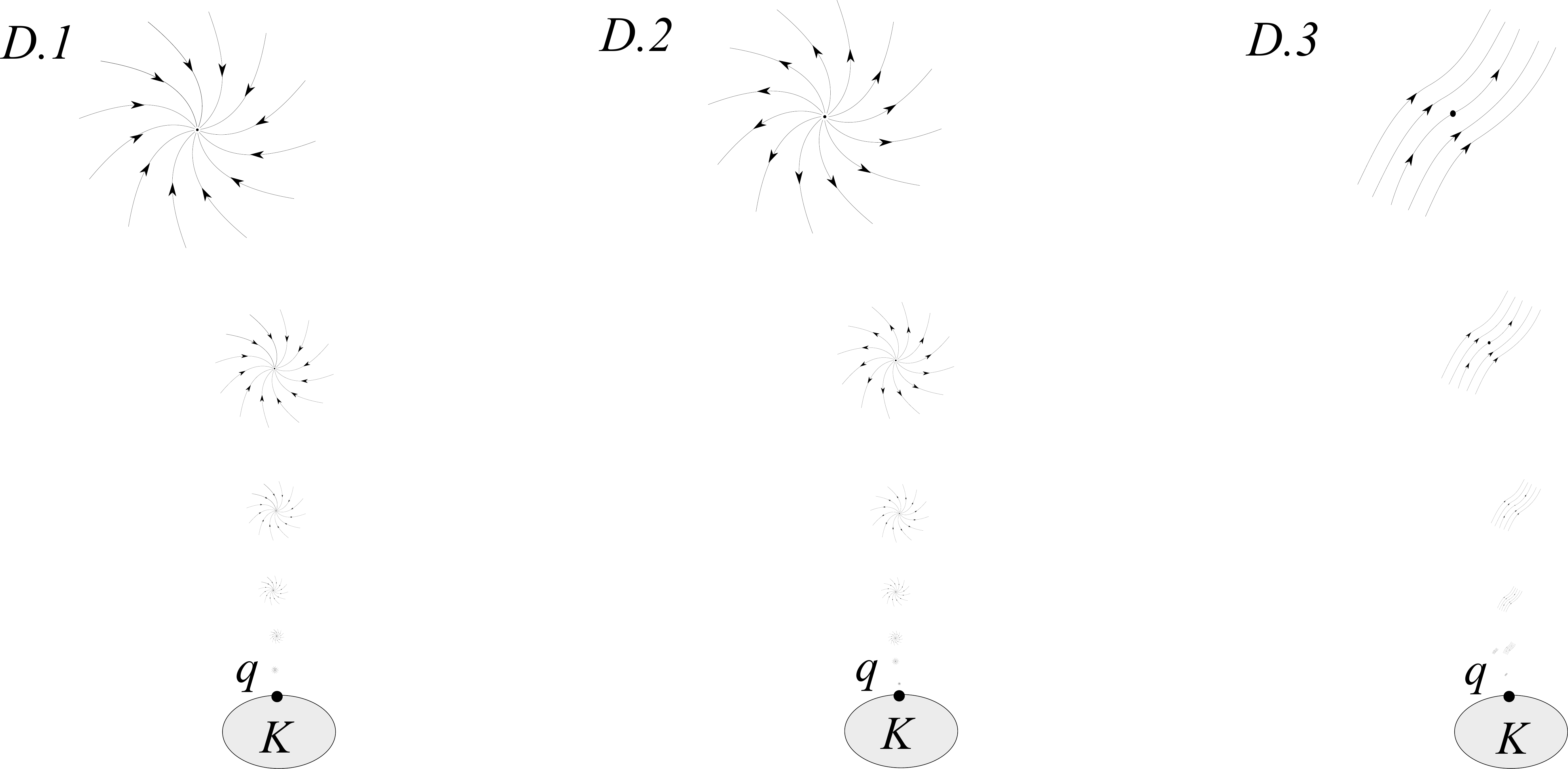}
\par\end{centering}

\caption{}
\end{figure}

\textbf{\emph{E}}\textbf{.}\textbf{\emph{i}}\textbf{\hspace*{0.2in}}$1\leq$$i\leq3$\\
\emph{there is a sequence of periodic orbits ~$\gamma_{_{n}}\subset M\setminus K$,
~$d_{_{H}}$-converging to some nonvoid, compact, connected, invariant
set ~}$Q\subset\mbox{bd\,}K$ ~\emph{and such that}\textbf{ }\emph{condition}
~\emph{$i.X$ }~\emph{is satisfied by all ~$\gamma_{_{n}}$.}

\medskip{}

\begin{figure}[H]
\noindent \begin{centering}
\includegraphics[scale=0.21]{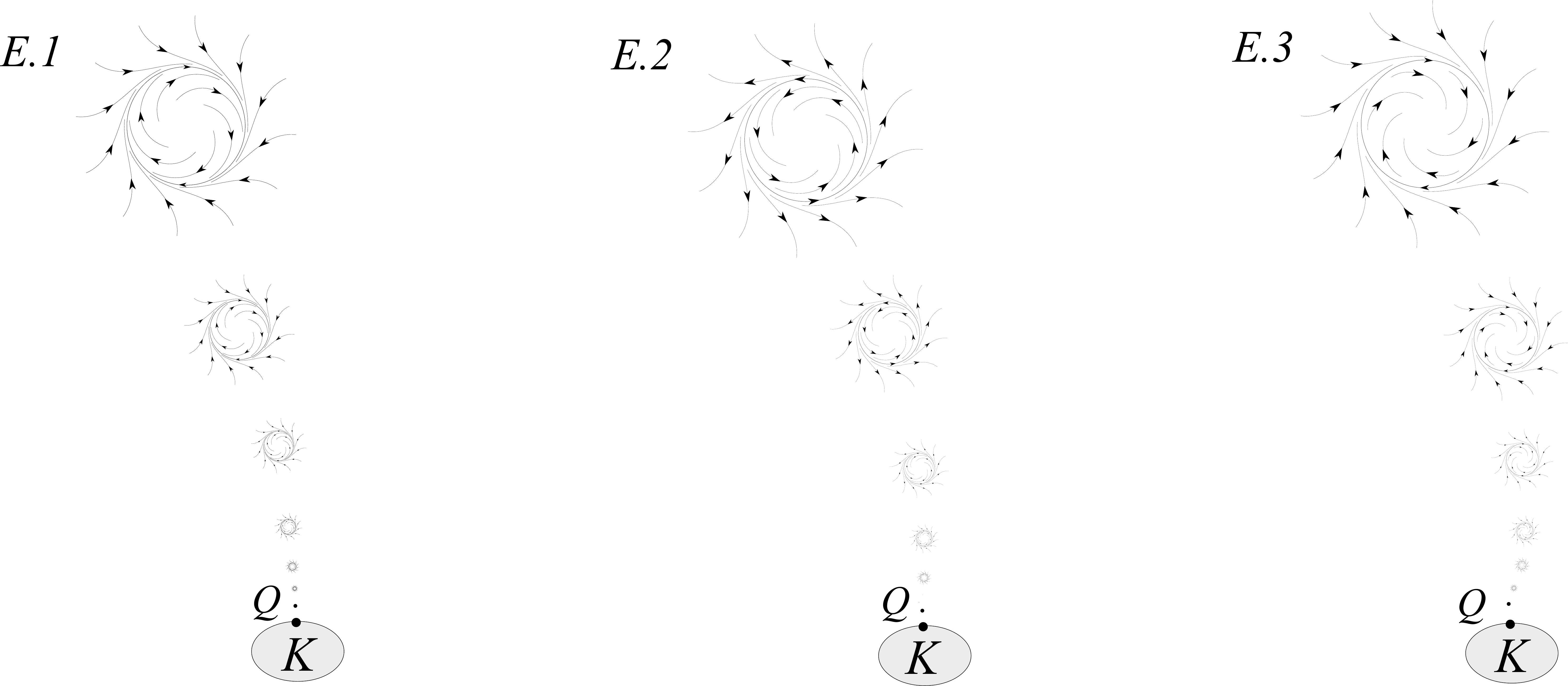}
\par\end{centering}

\caption{}
\end{figure}

\pagebreak{}

\bigskip{}

\begin{figure}[H]
\noindent \begin{centering}
\includegraphics[scale=0.19]{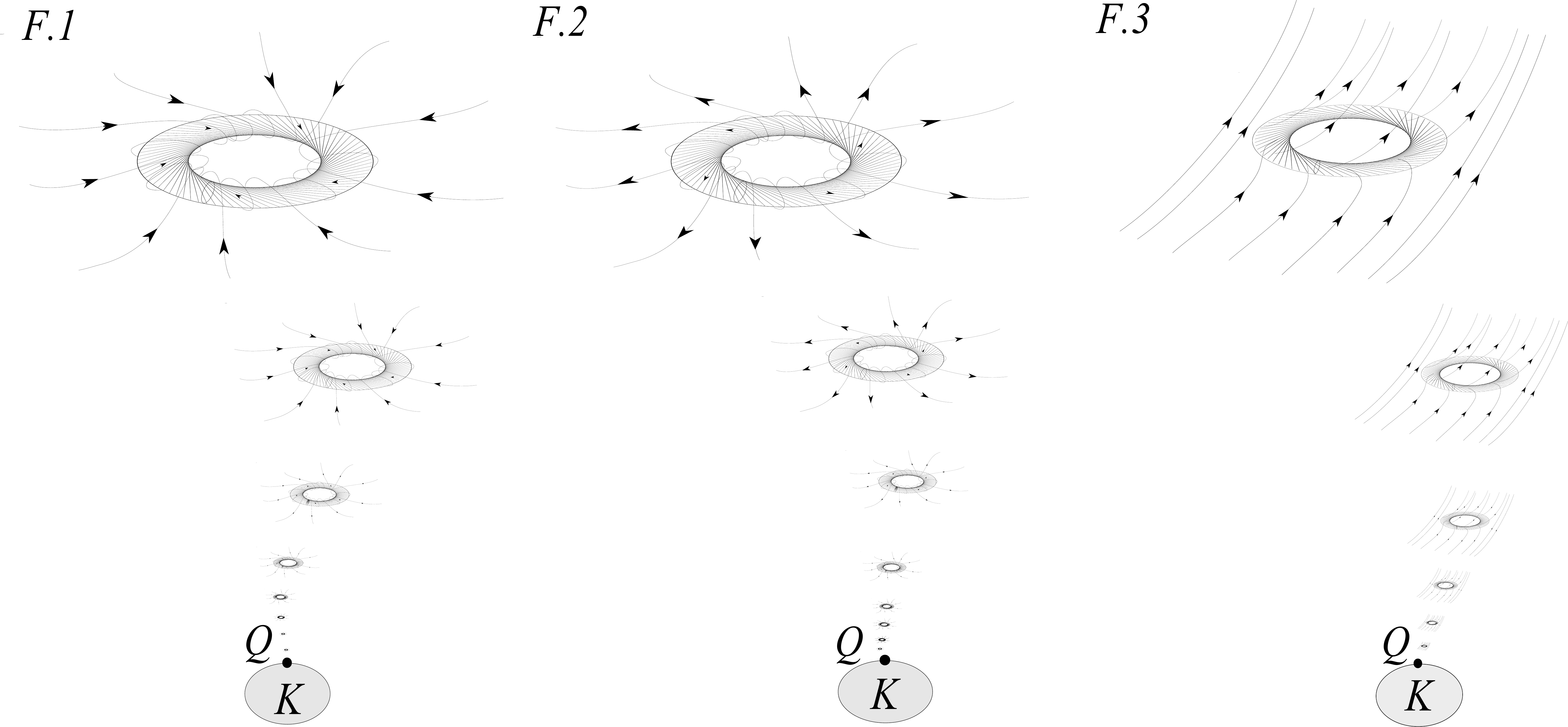}
\par\end{centering}

\caption{{\small this three examples show sequences of 2-tori carrying irrational
linear subflows.}}
\end{figure}

\textbf{\emph{F.i}}\textbf{\hspace*{0.2in}}$1\leq$$i\leq3$\\
\emph{there is a sequence of compact aperiodic minimal sets ~$\varGamma_{_{n}}\subset M\setminus K$,
~$d_{_{H}}$-converging to some nonvoid, compact, connected, invariant
set ~}$Q\subset\mbox{bd\,}K$ ~\emph{and such that}\textbf{ }\emph{condition}
~\emph{$i.X$ }~\emph{is satisfied by all ~$\varGamma_{_{n}}$.}

\begin{figure}[H]
\noindent \begin{centering}
\includegraphics[scale=0.23]{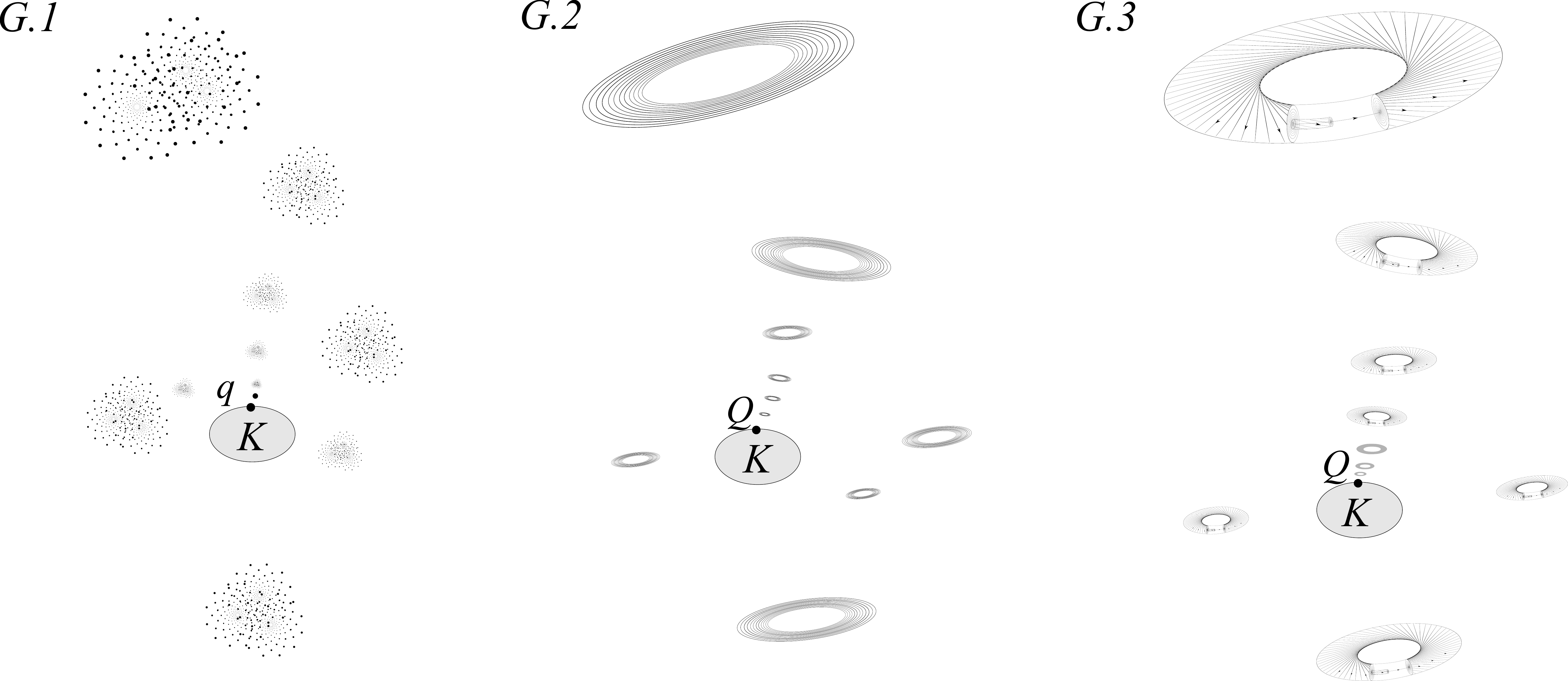}
\par\end{centering}

\caption{{\small figure G.3 shows nested continua of compact aperiodic minimal
sets (nested 2-tori carrying irrational linear subflows).}}
\end{figure}
\emph{or }\\
\textbf{\emph{G.}}\hspace*{0.15in}\emph{there is an open neighbourhood
~$U$ ~of ~$K$ ~such that any neighbourhood of a compact minimal
set ~$\varLambda\subset U\setminus K$ ~contains a} \emph{continuum
of compact minimal sets. ~Moreover, any ~$V\in\mathcal{N}_{_{K}}$
~actually contains a continuum of compact minimal sets disjoint from
~$K$, ~since in addition, at least one of the following four conditions
holds:}

\pagebreak{}

\begin{figure}[H]
\noindent \begin{centering}
\includegraphics[scale=0.3]{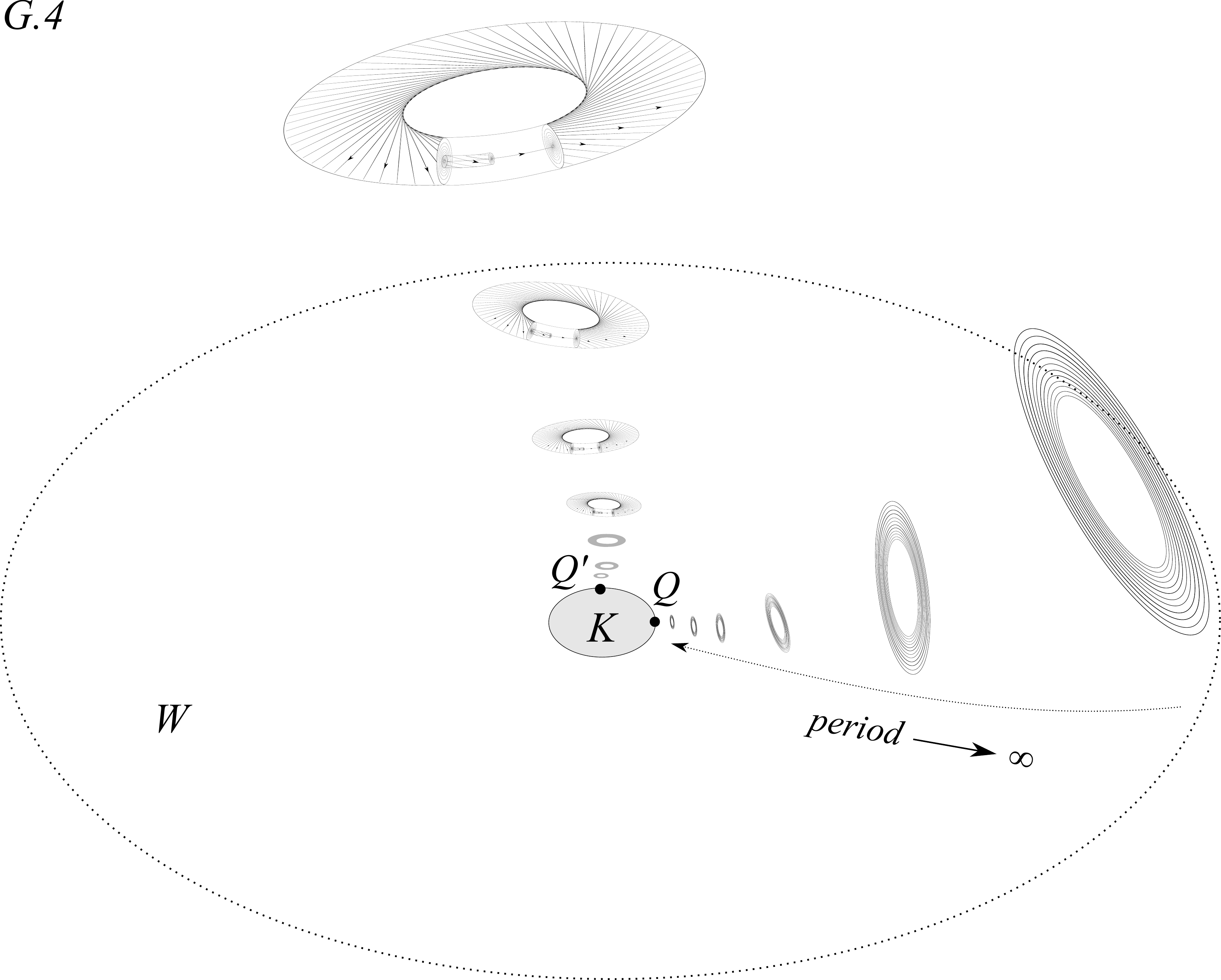}
\par\end{centering}

\caption{}
\end{figure}

\textbf{\hspace{0.6cm}}\textbf{\emph{G}}\textbf{.1}\hspace*{0.15in}\emph{any
neighbourhood of an equilibrium ~$z\in U\setminus K$ ~contains
a continuum of equilibria ~$e\in U\setminus K$, ~and there is a
sequence of these converging to some equilibrium ~}$q\in\mbox{bd}\, K$\emph{.}

\textbf{\hspace{0.6cm}}\textbf{\emph{G.}}\textbf{2}\hspace*{0.15in}\emph{any
neighbourhood of a periodic orbit ~$\gamma\subset U\setminus K$
~contains a continuum of periodic orbits ~$\zeta\subset U\setminus K$,
~and there is a sequence of these ~$d_{_{H}}$-converging to some
nonvoid, compact, connected, invariant set ~}$Q\subset\mbox{bd}\, K$\emph{.}

\textbf{\hspace{0.6cm}}\textbf{\emph{G.}}\textbf{3}\hspace*{0.15in}\emph{any
neighbourhood of a compact aperiodic minimal set ~$\varGamma\subset U\setminus K$
~contains a continuum of compact aperiodic minimal sets ~$\varLambda\subset U\setminus K$,
~and there is a sequence of these $d_{_{H}}$-converging to some
nonvoid, compact, connected, invariant set ~}$Q\subset\mbox{bd}\, K$\emph{.}

\emph{(hence in cases ~G.1, ~G.2 ~and ~G.3 ~ not only every neighbourhood
~$V$ ~of ~$K$ ~contains, respectively, a continuum of equilibria,
a continuum of periodic orbits, a continuum of compact aperiodic minimal
sets, all disjoint from ~$K$, ~but also, if ~$V$ ~is sufficiently
small, then, roughly speaking, we find on ~$V\setminus K$~ a kind
of {}``local super-abundance'' of the respective types of compact
minimal sets).}

\textbf{\hspace{0.6cm}}\textbf{\emph{G.}}\textbf{4}\hspace*{0.15in}\emph{there
is a set ~$P$ ~of periodic orbits contained in ~$U\setminus K$
~and a set ~$A$ ~of compact aperiodic minimal sets contained in
~$U\setminus K$} \emph{~such that the following four conditions
hold:}

\hspace{0.45in}\textbf{.1}\hspace*{0.15in}\emph{there are sequences
~$\gamma_{_{n}}\in P$ ~and ~$\varGamma_{_{n}}\in A$, ~$d_{_{H}}$-converging,
respectively, to some nonvoid, compact, connected, invariant sets
~$Q$,~}$Q'\subset\mbox{bd\,}K$.

\hspace{0.45in}\textbf{.2}\hspace*{0.15in}\emph{ every neighbourhood
~$D$ ~of a periodic orbit ~$\gamma\in P$ ~(resp. of a compact
aperiodic minimal set ~$\varGamma\in A\,$) contains a continuum
of periodic orbits (resp. of compact aperiodic minimal sets) ~and
if ~D ~is sufficiently small, all of them belong to ~$P$ ~(resp.
to ~$A\,$).}

\emph{(hence not only every open neighbourhood ~$V$ ~of ~$K$
~contains both a continuum ~$P(V)$ ~of periodic orbits and a continuum
~$A(V)$ ~of compact aperiodic minimal sets, all disjoint from ~$K$,
~but also, roughly speaking, both ~$P(V)$ ~and ~$A(V)$ ~display
a kind of {}``local super-abundance'' of the respective types of
compact minimal sets).}

\hspace{0.45in}\textbf{.3}\hspace*{0.15in}\emph{K ~is} \emph{bi-stable}
\emph{in relation to ~$P^{*}=\underset{\gamma\in P}{\bigcup}\gamma$
}~\emph{and ~$A^{*}=\underset{\varGamma\in A}{\bigcup}\varGamma$,
~i.e given any neighbourhood ~$W$ ~of ~$K$, ~all periodic orbits$\big/$compact
aperiodic minimal sets ~$\varLambda\in P\,\sqcup\, A$ ~at a sufficiently
small distance from ~$K$ ~are contained in ~$W$.}

\hspace{0.45in}\textbf{.4}\hspace*{0.15in}\emph{given any ~$n\geq1$,
~all periodic orbits ~$\gamma\in P$ ~contained in a sufficiently
small neighbourhood ~$V$ ~of ~$K$ ~have period ~$>n$ ~(note
that by .}1\emph{ and .}2\emph{ there is always a continuum of periodic
orbits~$\gamma\in P$ ~on ~$V\setminus K\,$).}

\emph{(hence, arbitrarily near ~$K$ ~there is a continuum of periodic
orbits disjoint from ~$K$ ~having periods larger than any prescribed
value).}

\begin{figure}[H]
\noindent \begin{centering}
\includegraphics[scale=0.2]{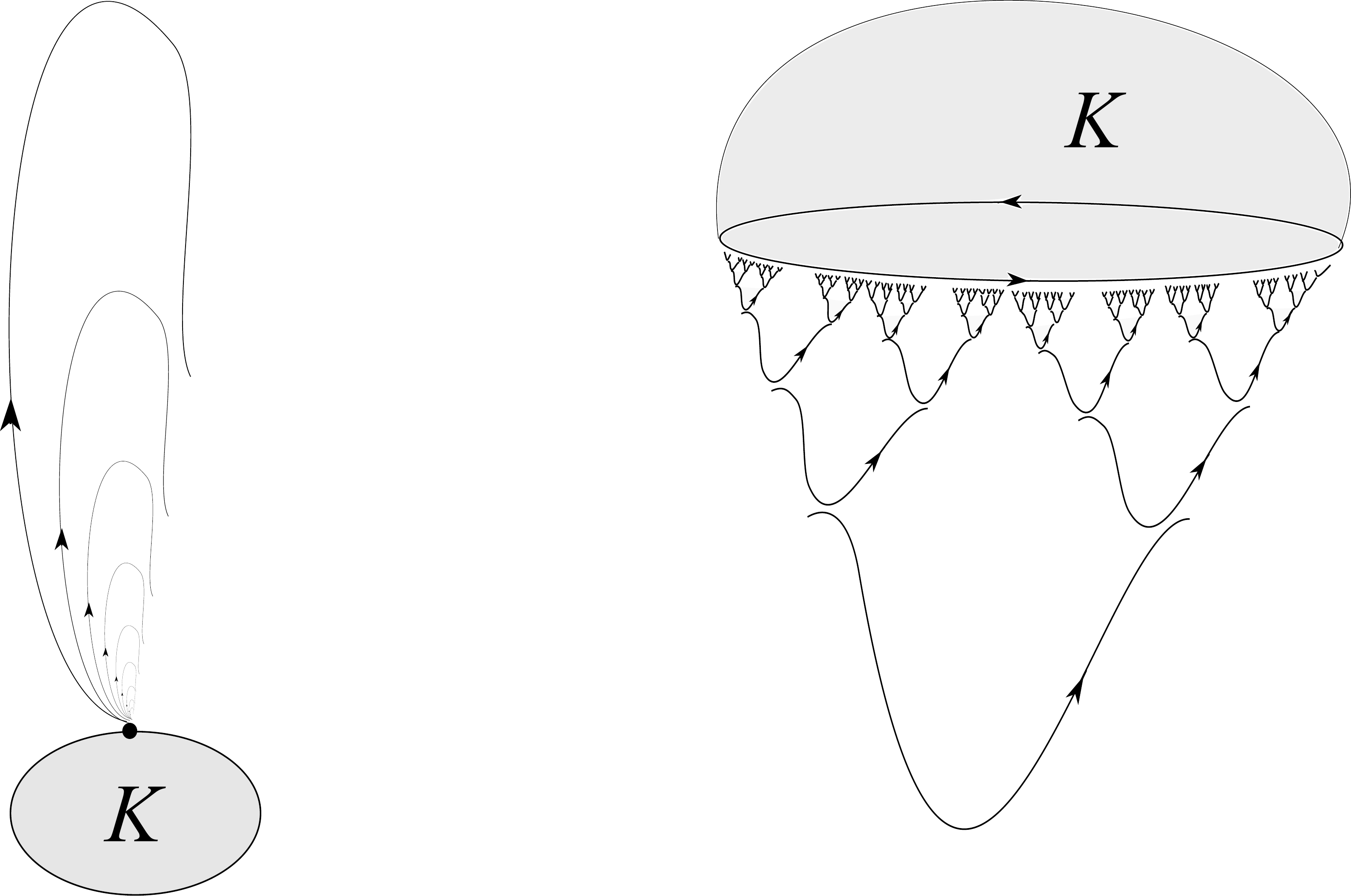}
\par\end{centering}

\caption{{\small for other figures illustrating case ~}\textbf{\small H}{\small{}
see figs. 13 to 15}.}
\end{figure}
\emph{If none of the sixteen conditions }~\textbf{\emph{A}}\emph{~
to }~\textbf{\emph{G.4}}\emph{~ holds, then any neighbourhood of
~$K$~ will contain orbits of infinite height (see fig. 10), more
precisely,}

\textbf{\hspace*{0.2in}}\textbf{\emph{H.}}\hspace*{0.15in}\emph{given
any ~$U\in\mathcal{N}_{_{K}}$, ~there is a sequence of orbits ~$\gamma_{_{n}}\subset U\setminus K$~
such that}\[
\mbox{cl}\,\gamma_{_{1}}\supsetneq\mbox{cl}\,\gamma_{_{2}}\supsetneq\cdots\cdots\supsetneq\mbox{cl}\,\gamma_{_{n}}\supsetneq\cdots\cdots\]

In condition \textbf{~H~} we shall distinguish another 12 relevant
cases, making up a total of $28=16+12$ cases in all. ~Observe that
the occurrence of the above strict inclusion chain is an invulgar
dynamical phenomenon. It is in deep contrast to what happens, for
example, in the case of an aperiodic minimal set ~$\varLambda$ ~of
a flow on a compact metric space $\big($as it is well known, for
any ~$x$ ~belonging to such ~$\varLambda$, ~$\alpha(x)\,=\omega(x)=\mbox{cl\,}\mathcal{O}(x)=\varLambda\,$$\big).$

\emph{Important remark}: ~the reader should keep in mind the following
two elementary facts about locally compact metric spaces ~$M$:

$\bullet$~~any neighbourhood of a compact set ~$K\subset M$ ~contains
a compact ball ~$B[\, K,\delta\,]$, ~for some ~$\delta>0$.

$\bullet$~~every open$\big/$closed subset of ~$M$ ~is also
locally compact.

\emph{Remark }1\emph{.} ~Before entering the proof of Theorem 1,
we make a few brief comments concerning its applicability in the context
of Differentiable Dynamics and also mention some of its consequences
for the topological structure of the set ~$\mbox{CMin}(M)$ ~of
all compact minimal sets of the flow, endowed with the Hausdorff metric
~$d_{_{H}}$.

~~~Let ~$M$ ~be a (2nd countable, Hausdorff) connected ~$C^{\, m}$~
$(0\leq m\leq\infty)$ ~manifold (compact or not) with a ~$C^{\, r}$
~$(0\leq r\leq m)$ ~flow. ~Due to its topological nature, Theorem
1 (see section 6 for the complete statement) gives information, not
only about the possible behaviour of the flow near each one of its
compact, invariant, proper subsets (if there are any), but also, it
{}``illuminates'' the behaviour of the flow \emph{within }each closed,
invariant, connected subset~ $N$, ~provided it contains a compact,
invariant proper subset (this is always the case if ~$N$ ~is, in
addition, a \emph{non-minimal }compact\emph{ }set, see also remark
2). ~Moreover, if ~$M$ ~is \emph{non-compact},\emph{ }then ~$M$~
has an \emph{end-points compactification} ~$M^{\propto}=M\,\sqcup\,\mbox{\ensuremath{E}}(M)$
~that, roughly speaking, captures the different possible ways of
going to infinite on ~$M.$~ As it is well known, besides compact,\emph{
~$M^{\propto}$ }~is connected and metrizable and the flow ~$\phi$
~on ~$M$ ~(uniquely) extends to a ~$C^{\,0}$ ~flow ~$\theta$
~on ~$M^{\propto}$ $\big(\, C^{\, r}$~ on ~$M\,\big)$, ~with
the \emph{end points }~$e\in\mbox{\ensuremath{E}}(M)$ ~becoming
\emph{equilibria}. ~At each end point\emph{ ~$e\in\mbox{\ensuremath{E}}(M)$,
}~not only the differentiable, but also the topological manifold
structure may break \emph{$\big($i.e.} an end point\emph{ }may not
even have ~a neighbourhood $($in ~$M^{\propto})$ homeomorphic
to ~$\mathbb{R}^{n}\,\big).$ ~However as the extended flow is still
continuous at these points, we may apply ~Theorem 1 to the \emph{equilibrium
orbit ~$K=\{e\}$ }of each end point\emph{ }~$e\in E(M)$, ~therefore
obtaining valuable insight about the possible behaviour of the original
flow near each one of its \emph{{}``points at infinite''.}

~~~Theorem 1 has many interesting consequences, some of which can
already be deduced from its weak version. ~Bellow we give a selection
of some simple Corollaries. ~Part II of the present work will be
devoted to the investigation of more subtle implications.

\textbf{Definition.} ~Let ~$M$ ~be a metric space. ~A set ~$\mathfrak{E}\subset2^{M}$
~has elements \emph{arbitrarily near} ~$X\subset M$ ~if ~for
any $\epsilon>0$, ~$B(X,\,\epsilon)$ ~contains an element of ~$\mathfrak{E}$
~\emph{$\big($i.e. }~$\mathfrak{E}\big(B(X,\,\epsilon)\big)\neq\emptyset\big)$.
~In this case we also say that ~$X$ ~has elements of ~$\mathfrak{E}$
~\emph{arbitrarily nearby.}~ More restrictively, $\mathfrak{E}$
~has elements \emph{arbitrarily near (but) outside }~$X$ ~if ~for
any $\epsilon>0$, ~~$\mathfrak{E}\big(B(X,\,\epsilon)\setminus X\big)\neq\emptyset$.
~We also use the expression ~$X$ ~has elements of ~$\mathfrak{E}$
~\emph{outside arbitrarily nearby. }

Observe that the last two concepts are defined using \emph{the metric
of} ~$M$ ~and should not be confused with ~$d_{_{H}}-$\emph{nearness}.%
\footnote{however if ~$M$ ~is a locally compact metric space with a ~$C^{\,0}$
~flow, ~$\mathfrak{E}\subset\mbox{Ci}(M)$ ~and ~$X\in\mbox{CMin}(M)$,
~then ~$\mathfrak{E}$ ~has elements \emph{arbitrarily near} ~$X$
~implies ~$\mathfrak{E}$ ~has elements \emph{arbitrarily ~$d_{_{H}}-$near
}~$X$ ~(see Lemmas 4 and 5, section 5). %
}~ From Theorem 1 we get the following immediate consequences:

\textbf{Corollary 1. ~}\emph{~Let ~$M$ ~be a locally compact,
connected metric space with a ~$C^{\,0}$ ~flow. ~Then} \emph{every
nonvoid, compact invariant set, isolated from minimals and having
no orbits of infinite height arbitrarily nearby is either an attractor
or a repeller or stagnant.\smallskip{}
}

\textbf{Corollary 1}'\textbf{.}\emph{ ~Let ~$M$ ~be a locally
compact, connected metric space with a ~$C^{\,0}$ ~flow and ~$K$
~a (nonvoid) compact, invariant subset of ~$M$, ~isolated from
minimal sets. ~If ~$K$ ~is neither an attractor, nor a repeller,
nor stagnant, then orbits of infinite height occur arbitrarily near
(but) outside ~$K$}.%
\footnote{note that the if ~$M$ ~is compact then ~$M$ ~is always both
an \emph{attractor }and a \emph{repeller }of the flow, hence ~$K$
~must be a proper subset of ~$M$ ~anyway. %
}

\smallskip{}

Ahead we shall see that if the hypothesis of Corollary 1' hold, then
at least one out of three remarkable denumerable collections of orbits
of \emph{infinite height} (called ~$K$-\emph{trees}, ~$K$-$\alpha shells$
~and ~$\mbox{\ensuremath{K}}$-$\omega shells$) will occur in the
flow. ~Any of these implies the existence of orbits of \emph{infinite
height} \emph{arbitrarily near (but) outside} the compact invariant
~$K$. \medskip{}

\textbf{Corollary 2.} \emph{~Let ~$\theta$ ~be a ~$C^{\,0}$
~flow on a locally compact, connected metric space having only a
finite number of compact minimal sets. ~Then any (nonvoid) compact
invariant set that is neither an attractor, nor a repeller, nor stagnant
has orbits of infinite height outside arbitrarily nearby.\medskip{}
}

Again suppose ~$M$ ~is a \emph{locally compact, connected }metric
space with a ~$C^{\,0}$ ~flow ~$\theta$. ~Let\[
\mathfrak{A}:=\big\{ X\in\mbox{CMin}(M):\,\, X\mbox{ \,\ satisfies one of conditions \,}1.X\mbox{ \,\ to \,}3.X\big\}\]
that is, ~$\mathfrak{A}$ ~is the set of \emph{compact }minimal
sets\emph{ }of the flow that are either \emph{attractors}, or \emph{repellers},
or \emph{isolated from minimals and stagnant}. ~The next result shows
that if the compact minimal sets\emph{ }belonging to\emph{ }~$\mathfrak{A}$
~are not ~$d_{_{H}}-$dense in ~$\mbox{CMin}(M)$, ~then ~{}``\emph{$\mathfrak{c}-$abundance}''
~of minimal sets\emph{ }or orbits of \emph{infinite height }will
occur in the flow. ~In the above context,\emph{\medskip{}
}

\textbf{Corollary 3.} \emph{~Let ~$M$ ~be a locally compact, connected
metric space with a ~$C^{\,0}$ ~flow ~$\theta$}. ~\emph{If}
~$\mathfrak{A}$ ~\emph{is not} ~$d_{_{H}}-$\emph{dense} \emph{in}
~$\mbox{CMin}(M)$\emph{,} \emph{~then there is a nonvoid,} ~$\mathfrak{c}-$\emph{dense
in itself, ~$d_{_{H}}-$open subset of} ~$\mbox{CMin}(M)$ \emph{~or
there are orbits of} \emph{infinite height arbitrarily near every
~}$Y\in\mbox{CMin}(M)\setminus\mbox{cl\,}_{_{H}}\mathfrak{A}\neq\emptyset$\emph{.} 

We prove a stronger {}``local'' result. ~Corollary 2 then follows
letting ~$A=M$. \emph{\medskip{}
}

\textbf{Corollary 4.} \emph{~Let} ~$M$\emph{ ~be a locally compact,
connected metric space with a ~$C^{\,0}$ ~flow ~$\theta$ ~and
~$A$ ~an open subset of ~$M$. ~If the set}\[
\mathfrak{A}(A)=\big\{ X\in\mbox{CMin}(A):\,\, X\mbox{ \,\ satisfies one of conditions \,}1.X\mbox{ \,\ to \,}3.X\big\}\]
\emph{is not ~$d_{_{H}}-$dense in ~}$\mbox{CMin}(A)$\emph{, ~then
at least one of the following two situations occurs:}

\emph{a)~}~\emph{there is a nonvoid, ~$\mathfrak{c}-$dense in
itself, $d_{_{H}}-$open subset of ~}$\mbox{CMin}(M)$\\
\hspace*{4mm}\emph{ contained in ~}$\mbox{CMin}(A)$\emph{.}

\emph{b)~}~\emph{there are orbits of infinite height arbitrarily
near every}\\
\hspace*{5mm}$Y\in\mbox{CMin}(A)\setminus\mbox{cl\,}_{_{H}}\mathfrak{A}(A)\neq\emptyset$.

\emph{In particular, if there are only countably many compact minimal
sets in ~$A$, then case b) takes place.}

\medskip{}

\emph{Proof. ~}By hypothesis ~$A\subset M$ ~is open and ~$\varDelta:=\mbox{CMin}(A)\setminus\mbox{cl\,}_{_{H}}\mathfrak{A}(A)\neq\emptyset$,
~hence ~$\mbox{CMin}(A)$ ~is a nonvoid, ~$d_{_{H}}-$open subset
of ~$\mbox{CMin}(M)\subset\mbox{C\ensuremath{(M)}}$ ~(lemma 8,
section 5). ~The set \[
\varUpsilon:=\big\{ Y\in\varDelta:\,\mbox{there are orbits of \emph{infinite height }arbitrarily near \,\ensuremath{Y\,\big\}}}\]

is clearly a ~$d_{_{H}}-$closed subset of ~$\varDelta$. ~Suppose
~$\varTheta:=\varDelta\setminus\varUpsilon\neq\emptyset$ ~\emph{i.e.}
assume there are compact minimal sets in ~$\varDelta$ ~having no
orbits of \emph{infinite height }arbitrarily nearby. ~Then ~$\varTheta$
~is a nonvoid, ~$d_{_{H}}-$open subset of ~$\mbox{CMin}(M)$.
~Let ~$K\in\varTheta\subset\mbox{CMin}(A)$. ~Then since ~$K$
~is a compact minimal set, there is an open ~$U\in\mathcal{N}_{_{K}}$
~with compact closure contained in ~$A$ ~and such that ~$U$
~contains no ~$Y\in\mathfrak{A}(A)$ ~(this follows from lemma
4 (section 5), as ~$K\in\mbox{CMin}(A)\setminus\mbox{cl\,}_{_{H}}\mathfrak{A}(A)\,\big)$,
and also contains no orbit of \emph{infinite height} ~(observe that
these two facts together also imply that any~$\varLambda\in\mbox{CMin}(U)$
~belongs to ~$\varTheta\,\big)$. ~Hence none of the 13 conditions
~\textbf{A} ~to ~\textbf{F.3} ~and ~\textbf{H} ~of the provisional
version of Theorem 1 holds, thus by the same result, at least one
of the 4 conditions ~\textbf{G.1}~ to ~\textbf{G.4} ~must take
place.~ But any of these implies the existence of a \emph{continuum}
of compact minimal sets contained in every ~$B(K,\,\delta)$, ~$\delta>0$
~and thus, by lemma 6 (section 5), of a \emph{continuum }of compact
minimal sets in every ~$B_{_{H}}(K,\,\epsilon)$, ~$\epsilon>0$,
~and for ~$\epsilon$ ~small enough these are contained in ~$U$
~and thus must belong to ~$\varTheta$. ~Therefore ~$\varTheta$
~is ~$\mathfrak{c}-$\emph{dense in itself}\hfill{}$\blacksquare$

Obviously, every ~$X\in\mathfrak{A}(M)$ ~is ~$d_{_{H}}-$isolated
in ~$\mbox{CMin}(M)$ ~$\big($and thus ~$\{X\}$ ~is ~$d_{_{H}}-$open
in ~$\mbox{CMin}(M)\,\big)$, ~hence in the above context,\emph{\medskip{}
}

\textbf{Corollary 5.} ~\emph{If ~$\theta$ ~is a ~$C^{\,0}$ ~
flow on a locally compact, connected metric space ~$M$~ with only
countably many compact minimal sets and displaying no orbits of infinite
height, then the set ~$\mathfrak{A}(M)$ ~is ~$d_{_{H}}-$open
dense in ~}$\mbox{CMin}(M)$\emph{.\medskip{}
}

\emph{Remark 2. ~}Suppose ~$N$ ~is a locally compact, connected
metric space endowed with a ~$C^{\,0}$ ~flow ~$\phi$ ~and ~$M$
~a connected, closed invariant subset of ~$N$, ~containing a compact
invariant proper subset ~$K$. ~Then Theorem 1 applies to the subflow
~$(M,\theta)$~ where ~$\theta:=\phi|_{\mathbb{R}\times M}$, ~$M$
~endowed with the metric of ~$N$. ~In this situation all definitions
must be interpreted {}``within'' ~($M,\theta)$ ~\emph{i.e.} as
concerning this subflow ~(for example a nonvoid, compact invariant
set ~$\varLambda\subset M$ ~may be an \emph{attractor }in ~$(M,\theta)$
~without being one in ~$(N,\phi)\big)$. ~The next result shows
that if addition the phase space ~$N$ ~is \emph{locally connected
}and is separated by the compact invariant set ~$J$, ~then a finner
understanding of the flow behaviour near ~$J$ ~is possible. \emph{\medskip{}
}

\textbf{Corollary 6. ~}\emph{Let ~$N$ ~be a locally compact, connected
and locally connected metric space with a ~$C^{\,0}$~ flow ~and
~$J$ ~a compact, invariant proper subset of ~$N.$ ~Let ~$D$
~be a connected component of ~$N\setminus J$. ~Then ~Theorem
1 applies to ~}$M:=\mbox{cl\,}D$\emph{, ~~$\theta:=\phi|_{\mathbb{R}\times M}$,
~~$K:=M\,\cap\, J$.\smallskip{}
}

Roughly speaking, this result means that \emph{within the closure
of each connected component $D$~ of }~$N\setminus J$, ~at least
one of the 28 phenomena described in Theorem 1 (see section 6 for
the full statement) takes place near the compact invariant ~$(\mbox{cl\,}D)\,\cap\, J$~
(it being possible that within distinct components, different conditions
hold).

\emph{Proof.} ~Consider the collection ~$\varTheta$ ~of all connected
components\emph{ }of~ $N\setminus J$. ~Since ~$N$ ~is locally
connected and ~$N\setminus J$ ~is open, every ~$D\in\varTheta$
~is open in ~$N$, ~hence it cannot be closed as ~$N$ ~is connected.
~On the other hand, every ~$D\in\varTheta$ ~is closed in ~$N\setminus J$,
~hence ~$\emptyset\neq\mbox{bd\,}D=(\mbox{cl\,}D)\setminus D\subset J$.
~The invariance of each ~$D\in\varTheta$ ~now follows from that
of ~$N\setminus J$: ~~the orbit of a point ~$z\in D$ ~cannot
pass from ~$D$ ~to a different ~$D'\in\varTheta$ ~without intercepting
~$\mbox{bd\,}J\subset J$, ~and this is impossible since ~$N\setminus J\supset D$
~is invariant.~ Therefore ~$M:=\mbox{cl\,}D$ ~is a nonvoid, connected,
closed (and hence locally compact) invariant subset of ~$N$ ~and
~$K=(\mbox{cl\,}D)\,\cap\, J$ ~is a nonvoid, compact, invariant
proper subset of ~$M$. ~Define the (sub)flow ~$\theta:=\phi|_{\mathbb{R}\times M}$.
~Now endowed with the metric of ~$N,$ ~~$M$ ~is a compact,
connected metric space with a ~$C^{\,0}$ ~flow ~$\theta$ ~and
~$K$ ~is a compact, invariant (under ~$\theta)$ ~proper subset
of ~$M$. ~Theorem 1 can thus be applied to these ~$M$, ~$\theta$
~and ~$K$.\hfill{}$\blacksquare$

\medskip{}

\emph{Example. ~}Let ~$\phi$ ~be a ~$C^{\, r}$ ~$(r\geq0)$
~flow on ~$N=\mathbb{S}^{\, n}$ ~and ~$K\subset N$ ~an invariant,
codimension one, compact, connected ~$C^{\,0}$~submanifold. ~As
it is well known, by the generalized \emph{Jordan-Brouwer Separation
Theorem,}%
\footnote{J. W. Alexander, {}``A proof and extension of the Jordan-Brouwer
separation theorem'', \emph{Trans. A.M.S.} 23, 333-349 (1922). ~Alexander's
term {}``\emph{immersed''} means ~$C^{\,0}$\emph{-embedded.} ~Recall
that this work is prior to Whitney's foundational papers on the theory
of manifolds.%
} ~$K$ ~separates the flow into three invariant regions, ~$K$,
$B$ ~and ~$A$, ~the last two being the connected components of
~$N\setminus K$, ~with common boundary ~$K$. ~Besides applying
to ~$N$, ~$\phi$, ~$K$,~ Theorem 1\emph{ }also applies to ~$M=A\,\sqcup\, K$,
~~$\theta=\phi|_{\mathbb{R}\times M}$, ~$K$ ~and to ~$M=B\,\sqcup\, K$,
~~$\theta=\phi|_{\mathbb{R}\times M}$, ~$K$. ~Moreover, if ~$K$
~is not a minimal set, ~then it also applies to the (compact, connected,
metric) phase space ~$K$, ~giving, in this case, information about
the possible behaviour of the codimension one subflow ~$\theta=\phi|_{\mathbb{R}\times K}$
~near any compact, invariant, proper subset of ~$K$ ~(there is
at least one). ~This is always the case if, for example, ~$K$ is
the image of a $C^{\,0}$ ~embedding ~$\mathbb{S}^{2m}\hookrightarrow\mathbb{S}^{2m+1},$
~$n=2m+1$ ~$\big($since such ~$K$ must contain an equilibrium,
even if ~$\phi$ is only $C^{\,0}\,\big)$.\emph{\medskip{}
}

\noun{4.~SPECIAL ORBITAL STRUCTURES.}\medskip{}

~~~We will introduce three kinds of {}``orbital structures'',
~$X$\emph{-trees, ~$X$-$\alpha\,$shells} ~and ~$X$-$\omega\,$\emph{shells.
~}The reason for considering these denumerable collections of orbits
lies in the fact that they capture essential features of the {}``dynamical
complexity'' of those flows on which they occur. ~In particular,
their presence implies that arbitrarily near ~$X$~ there are orbits
having limit sets of an outstanding kind. 

~\emph{~~Throughout this section,} ~$X$ \emph{~is a compact,
invariant, proper subset of a ~$C^{\,0}$ ~flow on a locally compact
metric space} ~$M$.

\begin{figure}[H]
\noindent \begin{centering}
\includegraphics[scale=0.8]{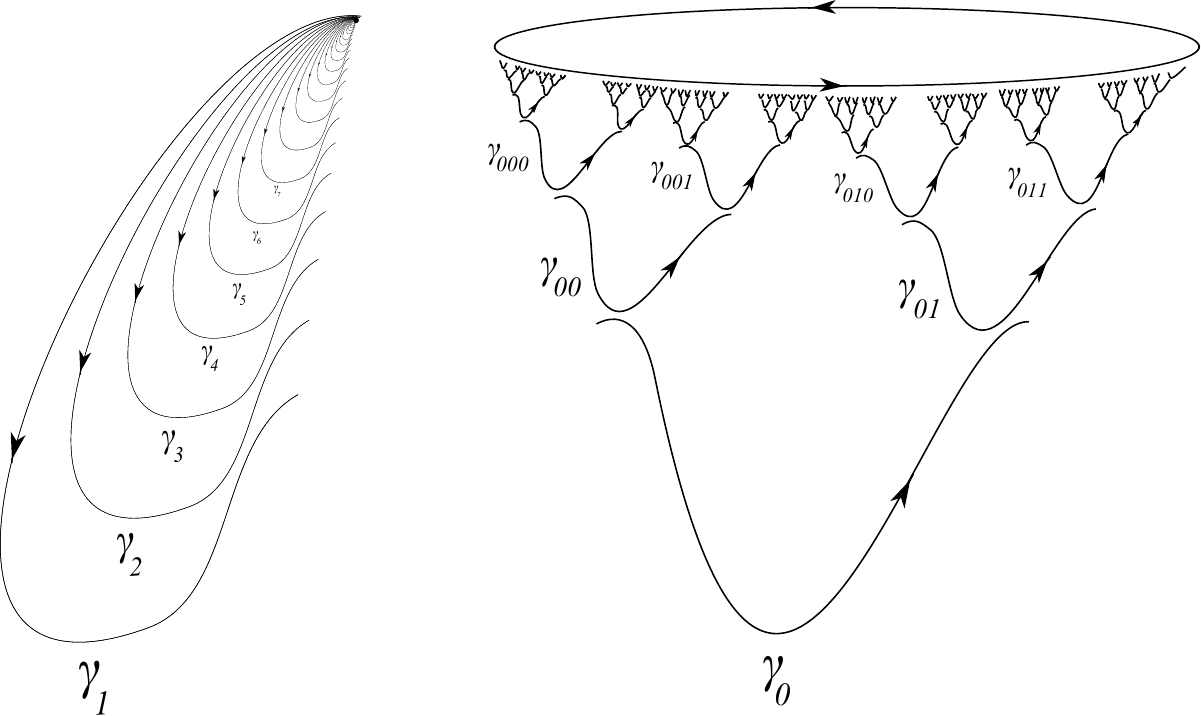}
\par\end{centering}

\caption{{\small (Left) ~a ~$X\mbox{-}\omega shell$ ~with ~$X$ ~an equilibrium
orbit. ~Time reversing the flow ~$X\mbox{-}\alpha shell$ ~is obtained.
~(Right) ~a ~$X\mbox{-}tree$ ~with ~$X$ ~a periodic orbit.}}
\end{figure}
\textbf{4.1 ~}\textbf{\emph{X-trees}}\noun{.}\smallskip{}

Let ~$\mbox{F}:=\{0,1\}$ ~and ~$\mbox{E}_{_{0}}:=\{0\}$. ~~Define\[
\mbox{E}_{_{n}}:=\{0\}\times\mbox{F}^{\, n},\,\, n\geq1\,\,\,\,\,\,\,\,\,\,\,\,\,\,\,\,\,\,\,\,\,\,\,\,\,\,\,\,\,\,\,\,\,\,\,\mathcal{E}:=\underset{n\geq0}{\bigsqcup}\mbox{E}_{_{n}}\,\,\,\,\,\,\,\,\,\,\,\,\,\,\,\,\,\,\,\,\,\,\,\,\,\,\,\,\,\,\,\,\,\,\,\mbox{E}_{_{\infty}}:=\{0\}\times\mbox{F}^{\mathbb{N}}\]
$\mbox{\ensuremath{\big(}}\mathcal{E}\mbox{\,\,\ and\,\,\ E}_{_{\infty}}$
~are, respectively, the set of finite and the set of infinite sequences
of 0's and 1's with first (left) digit $0\,\big)$. ~Since no risk
of ambiguity arises, commas and brackets are omitted in the representation
of both finite and infinite sequences of 0's and 1's ~e.g. we write
~$01$ ~and ~$00\ldots0\ldots\ldots$ ~instead of ~$(0,1)$ ~and
~$(0,0,\ldots,0,\ldots\ldots)$. ~~If ~$a,b\in\mathcal{E},$ ~$ab$
~represents, as usual, the element of ~\foreignlanguage{english}{$\mathcal{E}$}
~obtained by adjoining ~$b$~ to the right side of ~$a$. ~~For
each ~$v\in\mbox{E}_{_{\infty}}\,\,\,\,\,(v=0c_{_{1}}\ldots c_{_{n}}\ldots\ldots,\,\,\,\,\,\, c_{_{n}}\in\{0,1\}\mbox{ }\mbox{ \,\ for all \,}n\geq1)$,
~define ~$v_{_{0}}:=0$~ and ~$v_{_{n}}:=0\ldots c_{_{n}}$ ~~for
all ~$n\geq1$. 

\textbf{Definition. ~}If~$\gamma,\zeta\in\mbox{Orb}(M),$~ we denote
~$\zeta\subset\alpha(\gamma),$ ~~$\zeta\subset\omega(\gamma)$~~
and ~~$\zeta\subset\alpha(\gamma)\,\cup\,\omega(\gamma)$ ~~by~
$\gamma\,\overset{_{0}}{\succ\,\,}\zeta$, ~~$\gamma\,\overset{_{1}}{\succ\,\,}\zeta$
~and ~$\gamma\,\succ\,\,\zeta$, ~respectively. ~Note that all
these three relations are \emph{transitive} and ~$\gamma\,\overset{_{c}}{\succ\,\,}\zeta$
~and ~$\zeta\,\succ\,\,\,\xi$ ~implies ~$\gamma\,\overset{_{c}}{\succ\,\,}\xi$,
~~for ~$c\in\{0,1\}$.\medskip{}

Let ~$U$ ~be a \emph{compact neighbourhood of} ~$X$. ~A ~\emph{X-tree}
~is a pair ~$(\varTheta,\,\psi)$ ~where ~$\varTheta$ ~is a
collection of orbits contained in ~$U\setminus X$ ~and ~$\psi$
~is a surjective map

\[
\begin{array}{lll}
\psi:\mathcal{E} & \longrightarrow & \varTheta\subset\mbox{Orb}(U\setminus X)\\
\,\,\,\,\,\,\,\,\,\,\,\,\, a & \longmapsto & \gamma_{_{a}}\end{array}\]
 such that for any ~$b\in\mathcal{E}$, \begin{equation}
\begin{array}{c}
\gamma_{_{b}}\overset{_{0}}{\,\succ\,\,}\gamma_{_{b0}}\mbox{\,\,\,\,\,\,\,\,\,\ and\,\,\,\,\,\,\,\,}\gamma_{_{b0}}\not\succ\,\,\gamma_{_{b}}\\
\,\,\\
\gamma_{_{b}}\,\overset{_{1}}{\succ\,\,}\gamma_{_{b1}}\mbox{\,\,\,\,\,\,\,\,\,\ and\,\,\,\,\,\,\,\,}\gamma_{_{b1}}\not\succ\,\,\gamma_{_{b}}\end{array}\end{equation}

and for every ~$v\in\mbox{E}_{_{\infty}}$,\begin{equation}
\big|\mbox{cl\,}\gamma_{_{v_{_{n}}}}\big|_{_{X}}\longrightarrow0\end{equation}

$\gamma_{_{0}}$ ~is called the \emph{first orbit }of the $X\mbox{-}tree$
~(see fig.11).~ Observe that (1) implies (because of the transitivity
of relation ~$\succ$) ~that for every~ $v\in\mbox{E}_{_{\infty}}$,
~the sequence ~$\big(\gamma_{_{v_{_{n}}}}\big)$ ~is injective
\emph{i.e. }the ~$\gamma_{_{v_{_{n}}}}$'s ~are distinct and therefore
~$\varTheta$ ~is denumerable (since ~$\mathcal{E}$ ~is).\emph{
~X-trees} ~have significant dynamical properties, some of which
we single out:%
\footnote{recall that ~$A^{+}(X)$~~is the set of points of ~$M$ ~whose
~$\omega$-limit set intercepts both ~$X$ ~and ~$M\setminus X$.
~~$A^{-}(X)$~ is the corresponding negative concept.%
}

~~i)~~~~every ~$z\in\gamma\in\varTheta$ ~belongs to ~$A^{-}(X)\,\cap\, A^{+}(X)$

For each ~$v\in\mbox{E}_{_{\infty}}\,\,\,\,\,\,(v=0c_{_{1}}\cdots c_{_{n}}\cdots\cdots,\,\,\,\,\,\, c_{_{n}}\in\mbox{\{0,1\} }\mbox{for all \,}n\geq1),$

~~ii)\begin{equation}
\gamma_{_{v_{_{0}}}}\overset{_{c_{_{1}}}}{\,\succ\,\,}\gamma_{_{v_{_{1}}}}\,\overset{_{c_{_{2}}}}{\succ\,\,}\gamma_{_{v_{_{2}}}}\,\overset{_{c_{_{3}}}}{\succ}\cdots\cdots\overset{_{c_{_{n}}}}{\succ\,\,}\gamma_{_{v_{_{n}}}}\,\overset{_{c_{_{n+1}}}}{\succ\,\,}\gamma_{_{v_{_{n+1}}}}\overset{_{c_{_{n+2}}}}{\,\succ}\cdots\end{equation}
and\begin{equation}
q>p\,\,\,\implies\,\,\,\gamma_{_{v_{_{q}}}}\not\succ\,\,\gamma_{_{v_{_{p}}}}\end{equation}
thus\begin{equation}
\mbox{cl\,}\gamma_{_{v_{_{n}}}}\supsetneq\mbox{cl\,}\gamma_{_{v_{_{n+1}}}}\,\,\,\,\,\mbox{for all \,\ensuremath{n\geq}0}\end{equation}
~~iii)

\begin{equation}
\mbox{cl}\,\gamma_{_{v_{_{n}}}}\overset{d_{_{H}}}{\longrightarrow}\varLambda_{_{v}}:=\left(\underset{n\geq1}{\bigcap}\mbox{cl}\,\gamma_{_{v_{_{n}}}}\right)\in\mathfrak{S}(X)\end{equation}
\emph{Proof.~ }i)~~~if ~$z\in\gamma_{_{b}}\in\varTheta$, ~$b\in$$\mathcal{E}$
~then ~$\gamma_{_{b}}\,\overset{_{0}}{\succ\,\,}\gamma_{_{b0}}$
, ~~$\gamma_{_{b}}\,\overset{_{1}}{\succ\,\,}\gamma_{_{b1}}$ ~where~
$\gamma_{_{b0}},\gamma_{_{b1}}\in\varTheta\subset\mbox{Orb}(U\setminus X),$
hence both the $\,\alpha\mbox{-limit}$ and $\,\omega\mbox{-limit}$
sets of ~$z$~ have points outside ~$X$. ~On the other hand ~letting
~$k_{_{n}}:=\{0\}^{n}\in\mbox{F}^{\, n}$ and ~$l_{_{n}}:=\{1\}^{n}\in\mbox{F}^{\, n}$
~it follows immediately from (1) that ~$\gamma_{_{b}}\overset{_{0}}{\,\succ\,\,}\gamma_{_{bk_{_{n}}}}$
~and ~$\gamma_{_{b}}\overset{_{1}}{\,\succ\,\,}\gamma_{_{bl_{_{n}}}}$
~for all ~$n\geq1$; ~also (2) implies that both ~$\big|\gamma_{_{bk_{_{n}}}}\big|_{X}$
~and ~$\big|\gamma_{_{bl_{_{n}}}}\big|_{X}$ ~tend to zero when
~$n\longrightarrow+\infty$,~ thus both the ~$\alpha\mbox{-limit}$
and ~$\omega\mbox{-limit}$ sets of ~$z$ ~intercept ~$X$, ~since
these two sets are closed. \\
ii)~~~(3) is trivial; (4) and (5) follow from (1) because of
the transitivity of ~$\succ$.\\
iii)~~~$\mbox{cl}\,\gamma_{_{v_{_{n}}}}\in\mathfrak{S}(U)$
~and ~~$\mbox{cl}\,\gamma_{_{v_{_{n+1}}}}\subset$$\mbox{\,\ cl}\,\gamma_{_{v_{_{n}}}}$
~for all ~$n\geq0,$ ~therefore by Lemma 11 (section 5), ~$\mbox{cl}\,\gamma_{_{v_{_{n}}}}\overset{d_{_{H}}}{\longrightarrow}\varLambda_{_{v}}\in\mathfrak{S}(U)$
~since ~$\mathfrak{S}(U)$ ~is compact $\big($recall that ~$U\in\mathcal{N}_{_{X}}$
~is compact$\big)$; ~on the other hand ~$\big|\mbox{cl}\,\gamma_{_{v_{_{n}}}}\big|_{X}\longrightarrow0$
~hence ~$\varLambda_{_{v}}\subset X$ ~and finally ~$\varLambda_{_{v}}\in\mathfrak{S}(X).$

Observe that if ~$(\varTheta,\psi)$ ~is a ~$X\mbox{-}tree$, ~then
given any\emph{ ~$a\in\mathcal{E}$, ~}letting\[
\varUpsilon=\big\{\gamma_{_{d}}:\, d=a\mbox{\,\,\,\ or\,\,\,}d=ab\mbox{,\,\,\,\,}b\in\mbox{F}^{\, n},\mbox{\,\,\,}n\geq1\big\}\]
and defining the surjective map\[
\begin{array}{llll}
\phi:\mathcal{E} & \longrightarrow & \varUpsilon\\
\,\,\,\,\,\,\,\,\,0 & \longmapsto & \zeta_{_{0}}\,:=\gamma_{_{a}}\,\,\,\,=\psi(a)\\
\,\,\,\,\,\,\,\,\,0b & \longmapsto & \zeta_{_{0b}}:=\gamma_{_{ab}}=\psi(ab) & \mbox{\,\,\,\ for each \,}b\in\underset{n\geq1}{\sqcup}\mbox{F}^{\, n}\end{array}\]
we get a ~$X\mbox{-}tree$ ~with first orbit ~$\gamma_{_{a}}$,
~whose orbits are contained in ~$\varTheta$. ~We call ~$(\varUpsilon,\,\phi)$
~a \emph{sub} $X$-\emph{tree} ~of ~$(\varTheta,\,\psi)$ ~and
commit a safe abuse of expression saying that ~$(\varUpsilon,\,\phi)$
~is contained in ~$(\varTheta,\,\psi)$.~ Note that ~~$\big|\zeta_{_{d}}\big|_{_{X}}\leq\big|\mbox{cl\,}\zeta_{_{0}}\big|_{_{X}}=\big|\zeta_{_{0}}\big|_{_{X}}=\big|\gamma_{_{a}}\big|_{_{X}}$
~for all ~$d\in\mathcal{E}$, ~since ~$\zeta_{_{d}}\subset\zeta_{_{0}}\,\cup\,\alpha(\zeta_{_{0}})\,\cup\,\omega(\zeta_{_{0}})=\mbox{cl\,}\zeta_{_{0}}$.
~Therefore, in virtue of (2), ~given a ~$X$-\emph{tree} ~$(\varTheta,\,\psi)$
~and an ~$\epsilon>0$, ~there is always a \emph{sub} $X$-\emph{tree
}of ~$(\varTheta,\,\psi)$\emph{ ~}with all its orbits contained
in ~$B(X,\,\epsilon)\setminus X$. \medskip{}

\textbf{4.2~~}\textbf{\emph{X-$\mathbf{\alpha\,}$shells ~and~
X-$\mathbf{\omega}\,$shells.}}\smallskip{}

~~~We define ~$X\mbox{-}\omega\, shells$. ~~$X\mbox{-}\alpha\, shells$
~are the time symmetric concept \emph{i.e.} a sequence of orbits
~$(\gamma_{_{n}})_{_{n\geq1}}$ ~is a ~$X\mbox{-}\alpha\, shell$
~if it is a ~$X\mbox{-}\omega\, shell$ ~in the time reversal flow
~$\phi(t,x)=\theta(-t,x).$ ~Let ~$U$ ~be a \emph{compact neighbourhood
of }~$X$. ~A ~$X\mbox{-}\omega\, shell$ ~is a sequence ~of
orbits ~$\gamma_{_{n}}\subset U\setminus X$ ~satisfying the following
three conditions:

$\bullet$~~~$\gamma_{_{n}}\subset B^{-}(X)$ ~for every ~$n\geq1$

$\bullet$~~~$\gamma_{_{n}}\,\overset{_{1}}{\succ\,\,}\gamma_{_{n+1}}$
~~$\mbox{and}$~~$\gamma_{_{n+1}}\,\not\succ\,\,\gamma_{_{n}}$,
~~$\mbox{for all}$ ~$n\geq1$

$\bullet$~~~$\big|\mbox{cl\,}\gamma_{_{n}}\big|_{X}\longrightarrow0$\smallskip{}

These imply%
\footnote{clearly ~$\gamma_{_{n}}\,\overset{_{1}}{\succ\,\,}\gamma_{_{m}}$
~for every ~$1\leq n<m$, ~thus ~$\omega(\gamma_{_{n}})\,\cap\, X\neq\emptyset$
~since ~$\big|\mbox{cl\,}\gamma_{_{m}}\big|_{X}\longrightarrow0$
~and ~$\omega(\gamma_{_{n}})$ ~is closed. ~On the other hand,
~$\omega(\gamma_{_{n}})\not\subset X$ ~because ~$\gamma_{_{n+1}}\subset\omega(\gamma_{_{n}})$
~and ~$\gamma_{_{n+1}}\subset M\setminus X$. ~Hence ~$\gamma_{_{n}}\subset\, A^{+}(X)$. %
} that ~$\gamma_{_{n}}\subset A^{+}(X)$ ~for every ~$n\geq1$ ~and
hence\begin{equation}
\gamma_{_{n}}\subset B^{-}(X)\,\cap\, A^{+}(X)\mbox{\,\,\,\,\,\,\,\ for every \,}n\geq1\end{equation}
Also, the sequence ~$(\gamma_{_{n}})$ ~is necessarily injective
\emph{i.e.} the ~$\gamma_{_{n}}$'s ~are distinct (see fig.11).
~~Again, as in the case of ~\emph{X-trees, }it is easily seen that:\[
\gamma_{_{1}}\,\overset{_{1}}{\succ\,\,}\gamma_{_{2}}\overset{_{1}}{\succ\,}\cdots\cdots\overset{_{1}}{\succ\,\,}\gamma_{_{n}}\,\overset{_{1}}{\succ\,\,}\gamma_{_{n+1}}\overset{_{1}}{\succ\,}\cdots\]
\[
q>p\,\,\,\implies\,\,\,\gamma_{_{q}}\not\succ\,\,\gamma_{_{p}}\]
\[
\mbox{cl\,}\gamma_{_{n}}\supsetneq\mbox{cl\,}\gamma_{_{n+1}}\mbox{\,\,\,\,\ for all \,\,\ensuremath{n\geq}1}\]
\[
\mbox{cl}\,\gamma_{_{n}}\overset{d_{_{H}}}{\longrightarrow}\varLambda:=\left(\underset{n\geq1}{\bigcap}\mbox{cl}\,\gamma_{_{n}}\right)\in\mathfrak{S}(X)\]

$X\mbox{-}\alpha\, shells$ ~have exactly the same properties, interchanging
~$\alpha$ ~with ~$\omega$, ~$+$ ~with ~$-$~ and changing
~$\overset{_{1}}{\succ\,\,}$~ to ~$\overset{_{0}}{\succ\,\,}$~everywhere.~
Obviously, ~if ~$(\gamma_{_{n}})_{_{n\geq1}}$ ~is ~a ~$X\mbox{-}\omega\, shell$
~then ~any subsequence ~$\big(\gamma_{_{n_{i}}}\big)_{_{i\geq1}}$
~is also a ~$X\mbox{-}\omega\, shell$ ~and we call it a ~\emph{sub
$X\mbox{-}\omega shell$} ~of ~$(\gamma_{_{n}})_{_{n\geq1}}$. ~Therefore,
since ~$\big|\mbox{cl\,}\gamma_{_{n}}\big|_{X}\longrightarrow0$,
~~given any ~$\epsilon>0,$ ~a ~$X\mbox{-}\omega\, shell$ ~always
has a \emph{sub }$X\mbox{-}\omega\, shell$\emph{ ~}with all its
orbits contained in ~$B(X,\,\epsilon)\setminus X$. ~Analogue fact
holds for ~$X\mbox{-}\alpha\, shells$.\medskip{}

\noun{5.~LEMMAS.}\smallskip{}

~~~The following result gives an unusual characterisation of \emph{attractors}
in terms of the behaviour of the negative orbits of points outside
the compact invariant in question. ~It illustrates a topological-dynamical
phenomenon that plays a key role in the present work.\medskip{}

\textbf{Lemma 1.} ~\emph{Let} ~~$M$ ~\emph{be a locally compact
metric space with a} $C^{\,0}$ ~\emph{flow} ~$\theta$ \emph{~and}
~$K$~ \emph{a compact, invariant, proper subset of} ~$M$\emph{.
~Then ~$K$ ~is an attractor iff ~there is a neighbourhood ~$U$
~of ~$K$~ such that no point ~$z\in U\setminus K$ ~has its
negative orbit ~$\mathcal{O}^{-}(z)$ ~entirely contained in ~$U$.
~Analogously, ~$K$ ~is a repeller iff there is a ~$U\in\mathcal{N}_{_{K}}$
~such that ~$z\in U\setminus K\implies\mathcal{O}^{+}(z)\not\subset U$.\medskip{}
}

We need only to prove the characterisation of attractors\emph{ }in
Lemma 1 since a compact, invariant set is a repeller\emph{ }iff it
is an\emph{ }attractor\emph{ }in the time reversal flow.

\begin{figure}[H]
\noindent \begin{centering}
\includegraphics[scale=0.8]{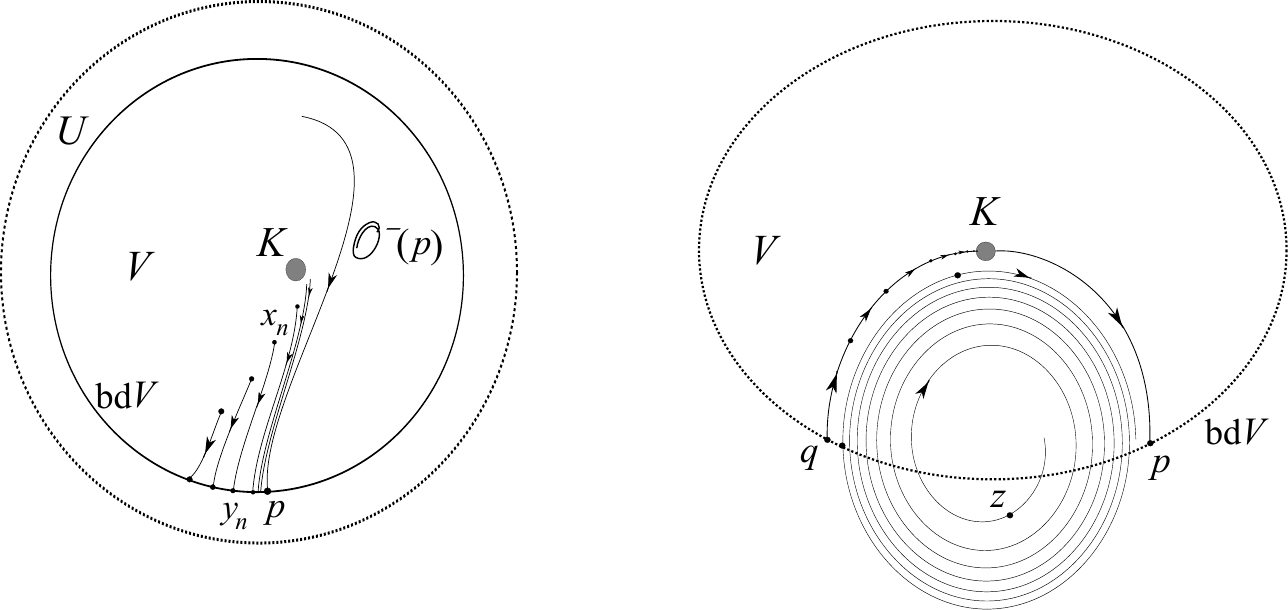}
\par\end{centering}

\caption{{\small (Left) Lemma 1, the non-stability}\emph{\small{} }{\small of
~$K$ ~implies the existence of a negative orbit trapped in ~$U\setminus K$,
~for every ~$U\in\mathcal{N}_{_{K}}\,$. (Right) Lemma 3.}}
\end{figure}

~Let ~$M$, ~$\theta,$ ~$K$ ~be as in Lemma 1. ~The following
elementary result (and the analogue fact for repellers)\emph{ }will
be used:

Claim: ~\emph{in the context of lemma }1\emph{, ~if ~$K$ ~is
an attractor and ~$z\in B^{+}(K)\setminus K$ ~then ~$\alpha(z)\subset\mbox{bd}\, B^{+}(K)\subset M\setminus B^{+}(K)$.}%
\footnote{~Suppose ~$z\in B^{+}(K)\setminus K$. ~It is easily seen that
~$B^{+}(K)$ ~is an open, invariant set hence ~$\alpha(z)\subset\mbox{cl\,}\mathcal{O}(z)\subset\mbox{cl\,}B^{+}(K)=B^{+}(K)\,\sqcup\,\mbox{bd\,}B^{+}(K)$.
~Necessarily, ~$\alpha(z)\,\cap\, B^{+}(K)=\emptyset\,$:~ otherwise
there is a ~$y\in\alpha(z)\,\cap\, B^{+}(K)$ ~implying both ~$\emptyset\neq\omega(y)\subset K$
~and ~$\omega(y)\subset\alpha(z)$ ~$\big(\omega(y)\subset\mbox{\ensuremath{\mbox{cl}\mathcal{O}(y)\subset}}\alpha(z)$,
~as ~$\alpha(z)$ ~is closed invariant), which by their turn imply
that ~$\alpha(z)\,\cap\, K\neq\emptyset$, ~and this contradicts
the \emph{stability }of ~$K$, ~since $z\in M\setminus K$. ~Therefore
~$\alpha(z)\subset\mbox{bd\,}B^{+}(K)\subset M\setminus B^{+}(K)$.
~Note however that we may have ~$\alpha(z)=\emptyset$. ~A time-symmetric
argument proves the analogue fact for repellers.%
}

\emph{Proof of Lemma 1.} ~$(\implies)$ ~Suppose ~$K$ ~is an
\emph{attractor}. ~Let ~$U$ ~be a compact neighbourhood of ~$K$
~contained in the open set ~$B^{+}(K)\in\mathcal{N}_{_{K}}$. ~Clearly
~$U$ ~and the compact ~$\mbox{bd\,}B^{+}(K)$ ~are disjoint,
hence, ~$z\in U\setminus K\subset B^{+}(K)\setminus K\implies\mathcal{O}^{-}(z)\not\subset U$
~since ~$\mathcal{O}^{-}(z)\subset U\implies\emptyset\neq\alpha(z)\subset U\subset B^{+}(K)$,
~which contradicts ~$\alpha(z)\subset\mbox{bd}\, B^{+}(K)\subset M\setminus B^{+}(K)$.

$($$\Longleftarrow)$ ~Suppose ~$U\in\mathcal{N}_{_{K}}$ ~is
such that ~$z\in U\setminus K\implies\mathcal{O}^{-}(z)\not\subset U$.
~Without loss of generality we may assume ~$U$ ~is compact.%
\footnote{if ~$U$ ~is not compact, take a compact ~$V\in\mathcal{N}_{_{K}}$
~contained in ~$U$ ~(this is always possible since ~$K$ ~is
compact ~and ~$M$ ~is locally compact). ~Clearly ~$z\in V\setminus K\implies\mathcal{O}^{-}(z)\not\subset V$.
~Let ~$U:=V$. %
}

Claim I.~~$K$ ~\emph{is} \emph{stable: ~}Suppose the contrary.\emph{~
}Then there is a ~$V_{_{0}}\subset\mathcal{N}_{_{K}}$ ~such that
every neighbourhood of ~$K$~ contains a point ~$x$ ~for which
~$\mathcal{O}^{+}(x)\not\subset V_{_{0}}$~ (observe that ~$x\in M\setminus K$
~since ~$K$ ~is invariant). ~Let ~$V$ ~be a compact neighbourhood
of ~$K$ ~such that ~$V\subset U\,\cap\, V_{_{0}}$. ~Since ~$M\setminus K$
~is invariant, the above remark implies that there are sequences
~$x_{_{n}}\in V\setminus K$~ and ~$t_{_{n}}\in\mathbb{R}^{+}$~
such that:

a)\textbf{\hspace*{3mm}}$\mbox{dist}(x_{_{n}},\, K)\longrightarrow0$\smallskip{}

b)\textbf{\hspace*{3mm}}$y_{_{n}}:=\theta(t_{_{n}},\, x_{_{n}})\in\mbox{bd\,}V$\smallskip{}

c)\textbf{\hspace*{3mm}}$\theta\big(\,[\,0,t_{_{n}}[\times\{x_{_{n}}\}\,\big)\subset V\setminus K\mbox{ \,\ which is equivalent to\,}\theta\big([\,-t_{_{n}},0\,[\times\{y_{_{n}}\}\,\big)\subset V\setminus K$\medskip{}

Since ~$\mbox{bd}\, V$ ~is compact, replacing ~$y_{_{n}}$ ~by
a convergent subsequence we may suppose ~$y_{_{n}}\longrightarrow p$
~for some ~$p\in\mbox{bd}\, V\subset U\setminus K\subset M\setminus K$.~
We claim that ~$\mathcal{O}^{-}(p)\subset V\subset U$.~ Assuming
the contrary \emph{i.e. }~$\mathcal{O}^{-}(p)\not\subset V$, ~there
is a ~$T<0$ ~such that ~$\lambda:=\mbox{dist}(\, p_{_{T}},\, V)>0$.~
By the continuity of the flow, there is a ~$\delta>0$ ~such that\[
z\in B(\, p,\,\delta)\implies z_{_{T}}\in B(\, p_{_{T}},\,\lambda/2)\]
hence for some ~$n_{_{0}}\geq1$,\[
n>n_{_{0}}\implies y_{_{n}}\in B(\, p,\,\delta)\implies\theta(T,\, y_{_{n}})\in B(\, p_{_{T}},\,\lambda/2)\subset M\setminus V\]
But necessarily~ $t_{_{n}}\longrightarrow+\,\infty$~ since ~$\mbox{dist}(x_{_{n}},\, K)\longrightarrow0$
~and ~$y_{_{n}}\in\mbox{bd}\, V$ ~$\big($observe that if ~$\big(t_{_{n_{i}}}\big)$~
is a bounded subsequence of ~$(t_{_{n}})$ ~then ~$\big\{ x_{_{n_{i}}}:\, i\geq1\big\}\subset\varTheta:=\theta\big([-t,0\,]\times\mbox{bd\,}V\,\big)$~
where~ $t:=\sup\big\{ t_{_{n_{i}}}:\, i\geq1\big\}$, ~$\varTheta$
~being a compact disjoint from the compact invariant ~$K,$~ hence~~$\mbox{dist}(x_{_{n_{i}}},\, K)\centernot\longrightarrow0,$
~a contradiction$\big)$. ~Thus for some $n_{_{1}}>n_{_{0}}$,\[
n>n_{_{1}}\implies\big(-t_{_{n}}<T<0\,\,\,\,\,\,\mbox{and}\,\,\,\,\,\,\theta(T,\, y_{_{n}})\not\in V\big)\]
in contradiction with c). ~Hence ~$p\in U\setminus K$ ~and ~$\mathcal{O}^{-}(p)\subset V\subset U$
~in contradiction with the initial hypothesis. ~Therefore ~$K$
~is \emph{stable}.

Claim II.~~$B^{+}(K)$ \emph{~is a neighbourhood of} ~$K$\emph{:
~}Recall that by hypothesis ~$U$ ~is a compact neighbourhood of
~$K$. ~Since ~$K$ ~is \emph{stable,} there is a ~$V\in\mathcal{N}_{_{K}}$
~such that ~$\mathcal{O}^{+}(V)\subset U$, ~which implies ~$\emptyset\neq\omega(x)\subset U$
~for every ~$x\in U$. ~We claim that ~$V\subset B^{+}(K)$. ~By
\emph{reductio ad absurdum,} ~suppose there are ~$x\in V$ ~and
~$y\in M\setminus K$ ~such that ~$y\in\omega(x)$.~ Clearly\[
\mbox{cl}\,\mathcal{O}^{+}(x)=\mathcal{O}^{+}(x)\,\cup\,\omega(x)\subset\mbox{cl\,}U=U\]
hence ~$y\in U\setminus K$.~ But by hypothesis, ~$\mathcal{O}^{-}(y)\not\subset U$
~which by its turn implies ~$\omega(x)\not\subset U$, ~since ~$\mathcal{O}(y)\subset\omega(x)$
~(as ~$\omega(x)$~ is an invariant set). ~We have reached contradiction.~
The proof of Lemma 1 is complete. \hfill{}$\blacksquare$ 

\smallskip{}

The following result is essentially a version of the \emph{Ura-Kimura-Bhatia
Theorem} for locally compact, connected metric spaces (see \emph{e.g.}
{[}1{]}, p.114). ~Our proof is given within the spirit of the present
work.\medskip{}

\textbf{Lemma 2.} ~\emph{Let} ~$M$ \emph{~be a locally compact,
connected metric space with a} $C^{\,0}$ ~\emph{flow} ~$\theta$
~\emph{and} $K$ \emph{a compact, invariant, proper subset of }~$M$.~\emph{
Then:}\textbf{\emph{ }}

\emph{either}\hspace*{3mm}I.\hspace*{3mm}\emph{K is an attractor.}

\emph{or}\hspace*{8.5mm}II.\hspace*{3mm}\emph{K is a repeller.}

\emph{or at least one of the following two situations occurs:}

\hspace*{12mm}III.\hspace*{3mm}\emph{K is isolated from minimals
and stagnant.}

\hspace*{12mm}IV.\hspace*{3mm}\emph{given any neighbourhood ~U
~of ~K,~~$U\setminus K$ ~contains an (entire)}\\
\hspace*{20mm}\emph{orbit. }\smallskip{}

\emph{Proof.} ~If ~$K$ ~is an attractor\emph{ }then it is easily
seen that none of conditions ~II, ~III~ and ~IV~ holds: clearly
~$B^{+}(K)\setminus K\neq\emptyset$ ~since ~$M$ ~is connected,
~$B^{+}(K)$ ~is an open neighbourhood of ~$K$ ~and ~$K$ ~is
a closed, proper subset of ~$M.$ ~Hence there is necessarily a
~$x\in M\setminus K$ ~such that ~$\emptyset\neq\omega(x)\subset K$
~and thus ~$K$ ~cannot be negatively stable\emph{,} therefore
~II~ does not hold; ~III~ implies the existence of a ~$y\in M\setminus K$
~such that ~$\emptyset\neq\alpha(y)\subset K$ ~which contradicts
the stability of ~$K$; ~IV contradicts the existence of neighbourhood
~$U$ ~in Lemma 1. ~Analogously if ~$K$ ~is a repeller\emph{
}then none of conditions I, III and IV holds. ~Suppose ~$K$ ~is
neither an attractor nor a repeller\emph{. ~}If ~III~ does not
hold, then

a)~~~\emph{K ~}is non-stagnant\smallskip{}

or\smallskip{}

b)~~~for every ~$W\in\mathcal{N}_{_{K}}$, ~~$\mbox{Min}(W\setminus K)\neq\emptyset$

If ~b) ~holds then condition ~IV ~is clearly satisfied. ~Suppose
~a) ~holds.~ Since ~$K$ ~is neither an attractor\emph{ }nor
a repeller, by lemma 1, given any compact ~$V\in\mathcal{N}_{_{K}}$,
~there are ~$x,\, y\in V\setminus K$ ~such that ~$\mathcal{O}^{-}(x)\subset V$
~and ~$\mathcal{O}^{+}(y)\subset V$ ~$\big($which implies ~$\emptyset\neq\alpha(x)\subset V$
~and ~$\emptyset\neq\omega(y)\subset V\,\big)$. ~Since ~$K$
~is non-stagnant we cannot have both ~$\alpha(x)\subset K$ ~and
~$\omega(y)\subset K$. ~If ~$\alpha(x)\not\subset K$~ then there
is a ~$z\in\alpha(x)\,\cap\,(V\setminus K)$ ~and\[
\mathcal{O}(z)\subset\alpha(x)\subset V\]
As ~$M\setminus K$ ~is invariant and ~$z\in M\setminus K$, ~$\mathcal{O}(z)\subset V\setminus K$.
~Since every neighbourhood of ~$K$ ~contains a compact ~$V\in\mathcal{N}_{_{K}}$,
~condition ~IV ~necessarily holds. ~If ~$\omega(y)\not\subset K$
~a time-symmetric argument conducts to the same conclusion.\hspace*{\fill}$\blacksquare$

\medskip{}

The next result is, in a certain sense, a counterpart to the \emph{Butler-Mcgehee
Lemma }(Butler and Waltman {[}3{]}, p.259).\medskip{}

\textbf{Lemma 3.}\emph{ ~Let} ~$M$ \emph{~be a locally compact
metric space with a} $C^{\,0}$ ~\emph{flow} ~$\theta$ ~\emph{and
~X~ a nonvoid, compact, invariant proper subset of ~M. ~If ~X~
is non-stagnant and ~$z\in A^{+}(X)$ ~}$\big($\emph{resp.} \emph{$z\in A^{-}(X)\,\big)$
~then given any ~$U\in\mathcal{N}_{_{K}}$ ~there is a ~$y\in\omega(z)$
~$\big($resp. $\, y\in\alpha(z)\,\big)$ ~such that ~$\mathcal{O}(y)\subset U\setminus X$}.\medskip{}

\emph{Proof.} ~The result will be proved in the case ~$z\in A^{+}(X),$
~a time symmetric argument disposes of case ~$z\in A^{-}(X)$. ~Assume
the hypothesis of Lemma 3 hold. ~Given any ~$U\in\mathcal{N}_{_{K}}$,
~let ~$0<\lambda<\mbox{min\ensuremath{\big(}1,}\,|\omega(z)|_{_{X}}\big)$
~be such that ~$V:=B[\, X,\lambda\,]\subset U$ ~is compact. ~We
will first show that there are ~$p,\, q\in\omega(z$) ~with ~$\mathcal{O}^{-}(p)$
~and ~$\mathcal{O}^{+}(q)$ ~both nonvoid and contained in ~$V\setminus X$.~
Since ~$\omega(z)\,\cap\, X\neq\emptyset$~ and ~$\omega(z)\not\subset V$
~it is easily seen that there are sequences ~$t_{_{n}},\,\, T_{_{n}}\in\mathbb{R}^{+}$
~such that:\smallskip{}

1)\textbf{\hspace*{3mm}}$t_{_{n}}<T_{_{n}}<t_{_{n+1}}$ ~~for all
~$n\geq1$\smallskip{}

2)\textbf{\hspace*{3mm}}$t_{_{n}}\longrightarrow+\infty$\smallskip{}

3)\textbf{\hspace*{3mm}}$\mbox{dist}\big(z_{_{t_{_{n}}}},\, K\big)\longrightarrow0$
\smallskip{}

4)\textbf{\hspace*{3mm}}$a_{_{n}}:=z_{_{T_{_{n}}}}\in\mbox{bd\,}V$\smallskip{}

5)\textbf{\hspace*{3mm}}$\theta\big([\, t_{_{n}},T_{_{n}}]\times\{z\}\big)\subset V\,\,\,\,\,\,\mbox{which is equivalent to}\,\,\,\,\,\,\theta\big([-(T_{_{n}}-t_{_{n}}),0\,]\times\{a_{_{n}}\}\big)\subset V$\smallskip{}

Let ~$\big(a_{_{n_{i}}}\big)$ ~be a subsequence of ~$(a_{_{n}})$
~converging to some point ~$p$ ~of the compact ~$\mbox{bd\,}V.$
~~Note that ~$a\in\omega(z)$~ since ~$T_{_{n_{i}}}\longrightarrow+\infty$.~
We claim that ~$\mathcal{O}^{-}(p)\subset V\setminus X$: ~~$\mathcal{O}(p)\subset M\setminus X$
~since ~\emph{$p$} ~belongs to the invariant set ~$M\setminus X$.~
To prove ~$\mathcal{O}^{-}(p)\subset V$ ~observe that together,
~$X$ ~is a compact invariant set, ~$a_{_{n_{i}}}\in\mbox{bd\,}V\subset M\setminus X$,
~$\mbox{bd\,}V$ ~is compact, ~$z_{_{t_{_{n_{i}}}}}=\theta\big(-(T_{_{n_{i}}}-t_{_{n_{i}}}),\, a_{_{n_{i}}}\big)$
~and ~$\mbox{dist}\big(z_{_{t_{_{n_{i}}}}},\, X\big)\longrightarrow0$
~imply, by an argument identical to the one given in the proof of
the Lemma 1, that we must have ~$-\big(T_{_{n_{i}}}-t_{_{n_{i}}}\big)\longrightarrow-\infty.$~
Now together,\[
\theta\big(\big[-\big(T_{_{n_{i}}}-t_{_{n_{i}}}\big),0\,\big]\times\{a_{_{n_{i}}}\}\big)\subset V,\]
\[
a_{_{n_{i}}}\longrightarrow p\mbox{\,\,\,\,\,\,\,\,\ and\,\,\,\,\,\,\,\,\,}-\big(T_{_{n_{i}}}-t_{_{n_{i}}}\big)\longrightarrow-\infty\]
imply in virtue of the continuity of the flow (again as in Lemma 1),
~that ~$\mathcal{O}^{-}(p)\subset V$~ and therefore ~$\mbox{cl\,}\mathcal{O}^{-}(p)=\mathcal{O}^{-}(p)\,\cup\,\alpha(p)\subset V$.~
A similar argument%
\footnote{observe that there are sequences ~$t_{_{n}},\, T_{_{n}}\in\mathbb{R}^{+}$
~satisfying conditions ~1) ~to ~4) ~above plus\[
5')\,\,\,\,\,\,\,\,\,\,\theta\big([\, T_{_{n}},t_{_{n+1}}]\times\{z\}\big)\subset V\,\,\,\,\,\,\mbox{which is equivalent to}\,\,\,\,\,\,\theta\big([\,0,(t_{_{n+1}}-T_{_{n}})\,]\times\{a_{_{n}}\}\big)\subset V\]
} (see fig.12) shows the existence of a ~$q\in\omega(z)$ ~such that
~~$\mathcal{O}^{+}(q)\subset V\setminus X$, ~which by its turn
implies ~$\mbox{cl\,}\mathcal{O}^{+}(q)=\mathcal{O}^{+}(q)\,\cup\,\omega(q)\subset V.$
~Note that both ~$\alpha(p)\neq\emptyset$ ~and ~$\omega(q)\neq\emptyset$,
~as ~$V$ ~is a nonvoid compact.~ Now since ~$X$ ~is \emph{non-stagnant}
and ~$p,\, q\in M\setminus X$, ~we cannot have both ~$\emptyset\neq\alpha(p)\subset X$
~and ~$\emptyset\neq\omega(q)\subset X$, ~therefore as ~$M\setminus X$
~is invariant, there is necessarily a ~$y\in\alpha(p)$ ~such that
~$\mathcal{O}(y)\subset V\setminus X$ ~or there is a ~$y\in\omega(q)$
~such that ~$\mathcal{O}(y)\subset V\setminus X.$~ In both cases
~$\mathcal{O}(y)\subset\omega(z)$ ~since ~$p,\, q\in\omega(z)$
~and ~$\omega(z)$ ~is closed and invariant.~ As ~$V\setminus X\subset U\setminus X$
~the proof is complete.~\hspace*{\fill}$\blacksquare$\smallskip{}

\textbf{Definition.}~ Let ~$M$ ~be a metric space with a ~$C^{\,0}$
~flow. ~A sequence ~$A_{_{n}}\subset M$ \emph{~approaches ~$X\subset M$}
~if for every ~$\epsilon>0$, ~$\, A_{_{n}}\subset B(\, X,\,\epsilon)$
~for all sufficiently large $n$ ~\emph{i.e. if ~$\big|A_{_{n}}\big|_{_{X}}\longrightarrow0$}.\smallskip{}

~In the next proof we use the following elementary result:

\emph{{}``In a metric space a sequence converges to a point ~z ~if
every subsequence contains a (sub)subsequence converging to ~z''.}\smallskip{}

\textbf{Lemma 4.}\emph{ ~Let ~$M$ ~be a locally compact metric
space} \emph{with a ~$C^{\,0}$} ~\emph{flow.} ~\emph{If a sequence}
~$\varLambda_{_{n}}\in\mbox{Ci(}M)$ ~\emph{approaches a compact
minimal set ~$S$ ~then ~$(\varLambda_{_{n}})$ ~actually ~$d_{_{H}}-$converges
to ~$S$ ~i.e.\[
\big|\varLambda_{_{_{n}}}\big|_{S}\longrightarrow0\,\implies\,\varLambda_{_{n}}\overset{d_{_{H}}}{\longrightarrow}S\]
Proof.} ~Let ~$U$ ~be a compact neighbourhood of ~$K$. ~Since
~$\big|\varLambda_{_{_{n}}}\big|_{S}\longrightarrow0$, ~given any
subsequence ~$\big(\varLambda_{_{n_{i}}}\big)$, ~there is a ~$i_{_{0}}\geq1$
~such that ~$\varLambda_{_{n_{i}}}\subset U$ ~for all ~$i>i_{_{0}}$.
~Now ~$\mbox{\ensuremath{\big[}Ci}(U),\, d_{_{H}}\big]$ ~is a
compact metric space by \emph{Blaschke Theorem} (see section 2), hence
by \emph{Blaschke Principle} there is a (sub)subsequence $\big(\varLambda_{_{n_{i_{k}}}}\big)$
~$d_{_{H}}-$converging to some nonvoid, compact, invariant set ~$Q\subset U.$
~But $\big|\varLambda_{_{n_{i_{k}}}}\big|_{S}\longrightarrow0$ ~implies
~$Q\subset S$ ~and since $S$ ~is a minimal set, ~$Q=S.$ ~Hence
the above convergence criterion is satisfied.\hfill{}$\blacksquare$\smallskip{}

Therefore, ~$X\in\mbox{CMin}(M)$ ~is an \emph{isolated minimal
set}%
\footnote{recall (section 2), that $X\in\mbox{CMin}(M)$ \emph{~}is an isolated
minimal (set) if there is a neighbourhood ~$U\subset M$ ~of ~$X$
~such that ~$U\setminus X$ ~contains no minimal set of the flow.%
} ~iff ~$X$ ~is ~$d_{_{H}}-$\emph{isolated in ~}$\mbox{CMin}(M)$
~(lemma 4 establishes $(\Longleftarrow)$; ~$(\Longrightarrow)$
~follows from the definition of the Hausdorff metric).\smallskip{}

\textbf{Lemma 5.}\emph{ ~Let ~$M$ ~be a locally compact metric
space} \emph{with a ~$C^{\,0}$} ~\emph{flow. ~If ~}$Q\in\mbox{CMin}(M)$
~\emph{then for each ~$\epsilon>0$ ~there is a ~$\delta>0$ ~such
that}\[
\mbox{Ci}\big(B(Q,\,\delta)\big)\subset B_{_{H}}(Q,\,\epsilon).\]
\emph{Proof.} ~Suppose the contrary. ~Then there is an ~$\epsilon>0$
~and sequences ~$\lambda_{_{n}}>0$ ~and ~$\varLambda_{_{n}}\in\mbox{Ci}\big(B(Q,\,\lambda_{_{n}})\big)$~
such that ~$\lambda_{_{n}}\longrightarrow0$ ~and ~$d_{_{H}}(\varLambda_{_{n}},\, Q)\geq\epsilon.$
~But this contradicts Lemma 4 since ~$\varLambda_{_{n}}\in\mbox{Ci}(M)$,
~~$Q\in\mbox{CMin}(M)$ ~and ~$\big|\varLambda_{_{n}}\big|_{Q}\longrightarrow0$.~\hfill{}$\blacksquare$\smallskip{}

The above useful relation between the metric ~$d$ ~of ~$M$ ~and
the \emph{Hausdorff metric} $d_{_{H}}$ ~of ~$\mbox{C}(M)$ ~will
be repeatedly used. ~Recall that if ~$A\subset M$ ~and ~$\mathfrak{M}\subset2^{M}$
~then ~$\mathfrak{M}(A)$ ~is the set of all ~$X\in\mathfrak{M}$
~contained in ~$A$ ~(see section 2). ~Note lemma 5 and Hausdorff
metric's definition together imply that if ~$\mathfrak{M}\subset\mbox{CMin}(M)$,
~then for every ~$X\in\mathfrak{M}$ ~and ~$\epsilon>0$, ~there
is a ~$\delta>0$ ~such that\[
B_{_{H}}(X,\,\delta)\,\cap\,\mathfrak{M}\,\subset\,\mathfrak{M}\big(B(X,\,\delta)\big)\,\subset\, B_{_{H}}(X,\,\epsilon)\,\cap\,\mathfrak{M}\]
hence there are arbitrarily small ~$d_{_{H}}-$neighbourhoods of
~$X\in\mathfrak{M}$ ~in the form ~$\mathfrak{M}\big(B(\, X,\,\delta\,)\big)$,
~a most useful fact.\smallskip{}

\textbf{Lemma 6.}\emph{ ~Let ~$M$ ~be a locally compact metric
space} \emph{with a ~$C^{\,0}$} ~\emph{flow. ~If ~}$\mathfrak{M}\subset\mbox{CMin}\big(M\big)$
~\emph{and for every} ~$X\in\mathfrak{M}$, ~$\epsilon>0$,\[
\#\,\mathfrak{M}\big(B(X,\,\epsilon)\big)=\mathfrak{c}\,\,\,\,\,\,\,\,\,\,\,\,\,\big(\mbox{\emph{resp.}}\mathcal{\mbox{ \,\,\,}}\mathfrak{M}\big(B(X,\,\epsilon)\setminus X\big)\neq\emptyset\big)\]
\emph{then} \emph{~$\mathfrak{M}$ ~is ~$\mathfrak{c}-$dense in
itself} \emph{~(resp. ~}$\mathfrak{M}$ ~\emph{is} \emph{~$d_{_{H}}-$dense
in itself).}\smallskip{}

Hence a set of compact minimal sets ~$\mathfrak{M}$ ~is ~$\mathfrak{c}-$\emph{dense
in itself} ~(resp. ~$d_{_{H}}-$dense in itself)\emph{ }~iff every
neighbourhood ~$U_{_{X}}\subset M$ ~of each ~$X\in\mathfrak{M}$
~contains a \emph{continuum }of elements of ~$\mathfrak{M}$ ~(resp.
an element of ~$\mathfrak{M}$ ~distinct from ~$X\,$). \smallskip{}

\emph{Proof of lemma 6. }~1) Suppose that for every ~$X\in\mathfrak{M}$
~and $\epsilon>0,$~ $\#\,\mathfrak{M}\big(B(X,\,\lambda)\big)=\mathfrak{c}$.
~Let ~$\epsilon>0$~ and ~$X\in\mathfrak{M}$ ~be given. ~Since
~$\mathfrak{M}\subset\mbox{CMin}(M)\subset\mbox{Ci}(M)$, ~by lemma
5, ~for ~$\delta>0$ ~sufficiently small,\[
\mathfrak{M}\big(B(X,\,\delta)\big)\,\subset\,\mbox{Ci}\big(B(X,\,\delta)\big)\,\subset\, B_{_{H}}(X,\,\epsilon)\]
therefore ~$\#\,\big(B_{_{H}}(X,\,\epsilon)\,\cap\,\mathfrak{M}\big)=\mathfrak{c}$
~since by hypothesis ~$\#\,\mathfrak{M}\big(B(X,\,\delta)\big)=\mathfrak{c}$.
~Hence ~$\mathfrak{M}$ ~is ~$\mathfrak{c}-$\emph{dense in itself}.
~2) ~Suppose now that for every ~$X\in\mathfrak{M}$ ~and ~$\epsilon>0,$
~~$\mathfrak{M}\big(B(X,\,\epsilon)\setminus X\big)\neq\emptyset$.
~Then, given any ~$X\in\mathfrak{M}$, ~there is a sequence ~$\varLambda_{_{n}}\in\mathfrak{M}\setminus\{X\}\subset\mbox{Ci}(M)$
~such that ~$\big|\varLambda_{_{n}}\big|_{X}\longrightarrow0$ ~and
since ~$X\in\mbox{CMin}(M)$, ~by lemma 4~ it follows that ~$\varLambda_{_{n}}\overset{d_{_{H}}}{\longrightarrow}X$.
~As all the $\varLambda_{_{n}}'s$ are distinct from $X,$ ~$X$
is not $d_{_{H}}-$isolated in ~$\mathfrak{M}$\emph{, }hence as
~$X\in\mathfrak{M}$ is arbitrary,\emph{ }~$\mathfrak{M}$ is $d_{_{H}}-$dense
in itself.\hfill{}$\blacksquare$\smallskip{}

\textbf{Lemma 7.} \emph{~Let ~$M$ ~be a locally compact metric
space} \emph{with a ~$C^{\,0}$} ~\emph{flow.} ~\emph{If ~$A\subset M$
~is open and ~}$\mbox{CMin}(A)$\emph{ ~is $d_{_{H}}-$dense in
itself ~then ~}$\mbox{CMin}(A)$\emph{ ~is ~$\mathfrak{c}-$dense
in itself. ~If ~$\mathfrak{M}$ ~is a ~$d_{_{H}}-$open and dense
in itself subset of ~}$\mbox{CMin}(M)$\emph{ ~then ~$\mathfrak{M}$
~is ~$\mathfrak{c}-$dense in itself}\smallskip{}

The proof of lemma 7 is presented on section 8.\smallskip{}

\textbf{Lemma 8.} \emph{~Let ~$M$ ~be a metric space. If} ~$\mathfrak{C}\subset\mbox{C}(M)$
\emph{~and ~$A\subset M$ ~is open, then ~$\mathfrak{C}(A):=\{X\in\mathfrak{C}:\,\, X\subset A\}$
~is $d_{_{H}}-$open in ~$\mathfrak{C}.$}\smallskip{}

\emph{Proof. ~}Let ~$X\in\mathfrak{C}(A).$~ Since ~$X$ ~is
compact and ~$A$ ~is open, ~there is a ~$\lambda>0$ ~such that
~$B(X,\,\lambda)\subset A$.~ On the other hand for any ~$\epsilon>0,$
~$Y\in B_{_{H}}(X,\,\epsilon)\implies Y\subset B(X,\,\epsilon)$.~
Therefore ~$Y\in\mathfrak{C\,}\cap B_{_{H}}(X,\,\epsilon)\implies Y\in\mathfrak{C}\big(B(X,\,\epsilon)\big)\implies Y\in\mathfrak{C}(A)$.~\hfill{}$\blacksquare$\smallskip{}

\textbf{Lemma} \textbf{9.} ~\emph{~Let ~$M$ ~be a metric space
with a ~$C^{\,0}$ ~flow and ~$x_{_{n}}$ ~a sequence of periodic
points with (minimal) period ~$\lambda_{_{n}}.$ ~If ~$x_{_{n}}\longrightarrow x$
~and} ~$\lambda_{_{n}}$ ~\emph{is convergent} ~\emph{then ~x
~is either a periodic point or an equilibrium}.\smallskip{}

\emph{Proof.} ~Assume ~$x$ ~is neither a periodic point nor an
equilibrium. ~Suppose ~$\lambda_{_{n}}\longrightarrow T\in[0,+\infty\,[$.~
Let ~$d:=\mbox{dist}(x,\,\, x_{_{T}})/3>0.$ ~By the continuity
of the flow, there is a ~$0<\delta<d$ ~such that\[
\big(\,\,\, z\in B(x,\,\delta)\,\,\,\,\,\,\mbox{and}\,\,\,\,\, t\in]\, T-\delta,T+\delta\,[\,\,\,\big)\,\,\,\implies\,\,\theta(t,\, z)\in B(x_{_{T}},\, d)\]
Since ~$x_{_{n}}\longrightarrow x$ ~and ~$\lambda_{_{n}}\longrightarrow T$,
~there is a ~$n_{_{0}}\in\mathbb{N}$ ~such that \[
n>n_{_{0}}\,\,\,\implies\,\,\,\big(\,\,\, x_{_{n}}\in B(x,\,\delta)\subset B(x,\, d)\,\,\,\,\mbox{and\,\,\,\,}\lambda_{_{n}}\in]\, T-\delta,T+\delta\,[\,\,\,\big)\]
and therefore, as ~$x_{_{n}}$ ~is periodic with period ~$\lambda_{_{n}}$,
~it follows that for ~$n>n_{_{0}},$ ~~$x_{_{n}}=\theta(\lambda_{_{n}},\, x_{_{n}})\in B(x_{_{T}},\, d)$
~and ~$x_{_{n}}\in B(x,\, d)$, which is absurd since ~$B(x,\, d)\,\cap\, B(x_{_{T}},\, d)=\emptyset$.
~This actually shows that ~$x$ ~is either an equilibrium or ~$x$
~is a periodic point with period ~$T/n$, ~for some $n\in$$\mathbb{N}.$
\hfill{}$\blacksquare$ \smallskip{}

\textbf{Lemma 10.}\emph{ ~Let ~$M$ ~be a locally compact metric
space} \emph{with a ~$C^{\,0}$} ~\emph{flow.} ~\emph{If ~Q ~is
a compact aperiodic minimal set, then for any ~$m\geq1$ ~there
is an~ $\epsilon>0$ ~such that}\[
\gamma\in\mbox{Per}(M)\,\,\mbox{\,\ and\,\,\,}\,\mbox{dist}(\gamma\,,\, Q)<\epsilon\,\,\,\,\,\implies\,\,\,\,\,\mbox{period}(\gamma)>m\]

\emph{Proof. ~}Suppose the contrary. ~Then there is a ~$m\geq1$
~and there are sequences ~$\gamma_{_{n}}\in\mbox{Per}(M)$ ~and
~$x_{_{n}}\in\gamma_{_{n}}$ ~such that ~$\mbox{dist}(x_{_{n}},\, Q)\longrightarrow0$
~and ~$0<\mbox{period}(x_{_{n}})=\mbox{period}(\gamma_{_{n}})\leq m.$
~Since ~$Q$ ~has a compact neighbourhood of the form ~$B[\, Q,\,\delta\,]$,
~for some ~$\delta>0$, ~it is easily seen (taking ~$K_{_{n}}:=B[\, Q,\,\delta/n\,]\,\big)$
that applying Lemma 11.b)\emph{ }followed by \emph{Bolzano-Weierstrass
Theorem, }we may select from ~$(x_{_{n}})$ ~a sub-sequence ~$\big(x_{_{n_{i}}}\big)$
~such that ~$x_{_{n_{i}}}\longrightarrow x\in Q$ ~and ~$\mbox{period}\big(x_{_{n_{i}}}\big)\longrightarrow\lambda\in[\,0,m\,]$.
~By lemma 9 this implies ~$x\in Q$ ~is either a periodic point
or an equilibrium, which is absurd since $Q\in\mbox{Am}(M).$~\hfill{}$\blacksquare$\smallskip{}

\textbf{Lemma~11}(\emph{Nested Compacts Lemma})\emph{. ~Let} ~$M$
~\emph{be a locally compact metric space and} ~$\mathfrak{C}$ \emph{~a
~$d_{_{H}}$-closed subset of} ~$\mbox{C}(M).$ \emph{~If} ~$K_{_{n}}\in\mathfrak{C}$
\emph{~and} ~$K_{_{n}}\supset K_{_{n+1}}$ ~\emph{for all} ~$n\geq1$
\emph{~then:}

A)~~~~$K_{_{n}}\overset{d_{_{H}}}{\longrightarrow}\left(\underset{n\geq1}{\bigcap}K_{_{n}}\right)\in\mathfrak{C}$,

B)~~~\emph{~every sequence} ~$x_{_{n}}\in K_{_{n}}$ \emph{~
has a subsequence converging to some} ~ $x\in\underset{n\geq1}{\bigcap}K_{_{n}}.$

C)~~~~\emph{every sequence} ~$\varLambda_{_{n}}\subset K_{_{n}}$
\emph{where} ~$\varLambda_{_{n}}\in\mathfrak{B}$,\emph{ }~$\mathfrak{B}$
~\emph{ a $d_{_{H}}$-closed subset of ~}$\mbox{C}(M)$\emph{$,$
~has a subsequence $d_{_{H}}$-converging to some ~$X\in\mathfrak{B}\left(\underset{n\geq1}{\bigcap}K_{_{n}}\right)$.}\smallskip{}

\emph{Proof.} ~~A)~~By \emph{Blaschske Theorem}, ~$\big[\,\mbox{C}(K_{_{1}}),d_{_{H}}\,\big]$
~is a compact metric space (see section 2) and in addition, ~$K_{_{n}}\in\mbox{C}(K_{_{1}})\subset\mbox{C}(M)$
~for all ~$n\geq1$, ~hence by \emph{Blaschke Principle} there
is a sub-sequence ~$\big(K_{_{n_{i}}}\big)$ ~converging to some
~$K\in\mbox{C}(K_{_{1}})\,\cap\,\mathfrak{C}$. ~This actually implies%
\footnote{using the elementary fact that if ~$A,B,C,Y\in\mbox{C}(M)$,~ $A\supset B\supset C$,~
$d_{_{H}}(A,\, Y)<\epsilon$ and ~$d_{_{H}}(C,\, Y)<\epsilon$ then
$d_{_{H}}(B,\, Y)<\epsilon$. %
} the whole sequence ~$(K_{_{n}})$ ~$d_{_{H}}$-converges to ~$K,$~
since $K_{_{n}}\supset K_{_{n+1}}$ ~for all ~$n\geq$1. ~Now ~$\underset{n\geq1}{\bigcap}K_{_{n}}\subset K$
~since for each ~$z\in\underset{n\geq1}{\bigcap}K_{_{n}}$, ~the
sequence ~$x_{_{n}}:=z\in K_{_{n}}$ ~converges to ~$z,$~ which
implies ~$z\in\mbox{lim}\, K_{_{n}}=K.$~ On the other hand, ~$K_{_{m}}$~
is compact for each ~$m\geq1$ ~and ~\foreignlanguage{english}{$K_{_{m}}\supset K_{_{m+k}}$}
~for all ~$k\geq$1, ~hence since each ~$\mbox{C}(K_{_{m}})$~
is ~$d_{_{H}}$-closed in ~$\mbox{C}(K_{_{1}})$ ~(by \emph{Blaschke
Theorem}) ~it follows that ~$\mbox{lim\,}K_{_{n}}=K\subset K_{_{m}}$
~for each ~$m\geq1,$~ therefore ~$K\subset\underset{n\geq1}{\bigcap}K_{_{n}}$
~and finally ~$K=\underset{n\geq1}{\bigcap}K_{_{n}}$. ~~B)~~Select
from ~$x_{_{n}}\in K_{_{n}}\subset K_{_{1}}$~ a subsequence ~$\big(x_{_{n_{i}}}\big)$
~converging to some ~$x\in K_{_{1}}$ ~(this is possible since
~$K_{_{1}}$ ~is compact). ~Then ~$x_{_{n_{i}}}\in K_{_{n_{i}}}$
~and ~$K_{_{n_{i}}}\overset{d_{_{H}}}{\longrightarrow}\underset{n\geq1}{\bigcap}K_{_{n}}$,
~thus necessarily ~$x\in\underset{n\geq1}{\bigcap}K_{_{n}}.$~~
C)~~By ~\emph{Blaschke Principle,} ~$\varLambda_{_{n}}\in\mathfrak{B}(K_{_{n}})\subset\mathfrak{B}(K_{_{1}})\subset\mbox{C}(K_{_{1}})$
~has a subsequence ~$(\varLambda_{_{n_{i}}})$ ~$d_{_{H}}$-converging
to some ~$X\in\mbox{C}(K_{_{1}})\,\cap\,\mathfrak{B}$.~ But ~$\varLambda_{_{n_{i}}}\subset K_{_{n_{i}}}\overset{d_{_{H}}}{\longrightarrow}\underset{n\geq1}{\bigcap}K_{_{n}}$
~implies ~$\varLambda_{_{n_{i}}}$ approaches ~$K=\underset{n\geq1}{\bigcap}K_{_{n}}$
~i.e. ~$\mbox{\ensuremath{\big|}}\varLambda_{_{n_{i}}}\big|_{K}\longrightarrow0,$
~thus necessarily ~$X\subset\underset{n\geq1}{\bigcap}K_{_{n}}$
~and finally ~$X\in\mathfrak{B}\left(\underset{n\geq1}{\bigcap}K_{_{n}}\right)$.\hfill{}$\blacksquare$

\medskip{}

\textbf{Lemma 12}(\emph{Cantor-Dirichlet Principle}).\emph{ ~If ~$A=\underset{1\leq n\leq n_{_{0}}}{\bigcup}A_{_{n}}$
~and ~A ~is infinite then~~$\#A_{_{n}}=\#A$ ~for some ~$1\leq n\leq n_{_{0}}.$}\medskip{}

\emph{Proof. }~Recall that if ~$D=B\,\cup\, C$, ~and ~$D$ ~is
infinite then ~$B$ ~or ~$C$ ~is infinite (otherwise ~$\#D\leq\#\, B+\#\, C<\infty\,\big)$,
~and also that using the Axiom of Choice, if ~\emph{$B$ ~}and
~\emph{$C$~ }are\emph{ }two sets\emph{, }at least one of which
is infinite, then ~$\#(B\,\cup\, C)=\mbox{max}\{\#B,\#C\}.$ ~Therefore
by finite induction over ~$1\leq n\leq n_{_{0}}$~ it follows that
~$\#A=\mbox{max}\{\#A_{_{n}}:\,1\leq n\leq n_{_{0}}\}.$~\hfill{}$\blacksquare$

\bigskip{}

\noun{6.1~THE MAIN THEOREM. TOPOLOGICAL STRUCTURE OF }~$\mbox{CMin}(M)$.

\medskip{}
~~~Again let ~$M$ ~be a locally compact, connected metric space
with a ~$C^{\,0}$ ~flow. ~Consider the following six propositions
where the variable ~$X$ ~assumes values in the set ~$\mathfrak{\mbox{Ci}(}M)$
~of nonvoid, compact, invariant subsets of ~$M$:\bigskip{}

\hspace{3mm}\emph{$1.X$} \hspace{6mm}$X$~ \emph{is an attractor.}

\hspace{3mm}$2.X$\hspace{6mm}$X$~ \emph{is a} \emph{repeller.}

\hspace{3mm}$3.X$\hspace{6mm}$X$~ \emph{is} \emph{isolated from
minimals and} \emph{stagnant. }

\hspace{3mm}$4.X$\hspace{6mm}$X$~ \emph{is} \emph{isolated from
minimals} \emph{and there is a} ~$X\mbox{-}\alpha\, shell$.

\hspace{3mm}$5.X$\hspace{6mm}$X$~ \emph{is} \emph{isolated from
minimals} \emph{and there is a} ~$X\mbox{-}\omega\, shell$.

\hspace{3mm}$6.X$\hspace{6mm}$X$~ \emph{is} \emph{isolated from
minimals and there is a} ~$X\mbox{-}tree$.

\pagebreak{}

\textbf{Theorem 1.}\emph{~ Let ~M ~be a locally compact, connected
metric space with a ~$C^{\,0}$ ~flow ~$\theta$ ~and ~K ~a
compact, invariant proper subset of ~M.~ Then:}

\emph{either}\textbf{}\\
\textbf{1.\hspace*{5mm}}$K$~ \emph{is an attractor.}

\emph{or}\textbf{}\\
\textbf{2.\hspace*{5mm}}$K$~\emph{ is a} \emph{repeller.\smallskip{}
}

\emph{or at least one of the following four conditions holds}:

\textbf{3.\hspace*{5mm}}$K$~ \emph{is} \emph{isolated from minimals
and} \emph{stagnant.}\smallskip{}

\textbf{4.\hspace*{5mm}}\emph{K}~ \emph{is} \emph{isolated from
minimals and there is a} ~$K\mbox{-}\alpha\, shell$.\smallskip{}

\textbf{5.\hspace*{5mm}}\emph{K}~\emph{ is} \emph{isolated from
minimals} \emph{and there is a} ~$K\mbox{-}\omega\, shell$.\smallskip{}

\textbf{6.\hspace*{5mm}}\emph{K}~ \emph{is} \emph{isolated from
minimals} \emph{and there is a} ~\emph{$K\mbox{-}tree$.}\smallskip{}

\emph{or}\textbf{\emph{ }}\emph{at least one of the following eighteen
conditions hold}

\textbf{7.i\hspace*{4mm}}$1\leq i\leq6$\textbf{}\\
\emph{there is a sequence} ~$\{e_{_{n}}\}\in\mbox{Eq}(M\setminus K),$
$d_{_{H}}-$\emph{converging to some}~$\{q\}\in\mbox{Eq}(\mbox{bd\,}K)$\\
\emph{and such that}\textbf{ }\emph{condition ~$i.X$ }~\emph{is
satisfied by all equilibrium orbits~$\{e_{_{n}}\}$.}

\textbf{8.i\hspace*{4mm}}$1\leq i\leq6$\textbf{}\\
\emph{there is a sequence} ~$\gamma_{_{n}}\in\mbox{Per}(M\setminus K),$
~$d_{_{H}}-$\emph{converging to some} ~$Q\in\mathfrak{S}(\mbox{bd\,}K)$\\
\emph{and such that condition ~$i.X$ ~is satisfied by all}
\emph{periodic orbits ~$\gamma_{_{n}}$.}

\textbf{9.i\hspace*{4mm}}$1\leq i\leq6$\textbf{}\\
\emph{there is a sequence} ~$\varGamma_{_{n}}\in\mbox{Am}(M\setminus K),$
~\emph{$d_{_{H}}-$converging to some} ~$Q\in\mathfrak{S}(\mbox{bd\,}K)$\\
\emph{and such that condition ~$i.X$ ~is satisfied by all compact
aperiodic minimals ~$\varGamma{}_{_{n}}$.}

\emph{or }\\
\textbf{10.\hspace*{3mm}}\emph{there is an open neighbourhood}
\emph{~$U$ ~}of\emph{ ~$K$ ~such that} ~$\mbox{CMin}(U\setminus K)$\emph{
~is} ~$\mathfrak{c}-$\emph{dense in itself }\textbf{~ }\emph{and}\textbf{
}\emph{at least one of the following four conditions holds:}

\textbf{\hspace*{5mm}10.1\hspace*{4mm}}$\mbox{Eq}(U\setminus K)$
\emph{~is a ~$\mathfrak{c}-$dense in itself set, }~$d_{_{H}}-$\emph{accumulating}
\emph{in ~}$\mbox{Eq}(\mbox{bd\,}K)$.

\textbf{\hspace*{5mm}10.2\hspace*{4mm}}$\mbox{Per}(U\setminus K)$
~\emph{is a ~$\mathfrak{c}-$dense in itself} \emph{set,} ~$d_{_{H}}-$\emph{accumulating
in} \emph{~}$\mathfrak{S}(\mbox{bd\,}K)$.

\textbf{\hspace*{5mm}10.3\hspace*{4mm}}$\mbox{Am}(U\setminus K)$
~\emph{is a} \emph{~$\mathfrak{c}-$dense in itself set,} ~$d_{_{H}}-$\emph{accumulating
in ~}$\mathfrak{S}(\mbox{bd\,}K)$.

\textbf{\hspace*{5mm}10.4\hspace*{4mm}}\emph{there are} $\mathfrak{c}-$\emph{dense
in itself sets }~$P\subset\mbox{Per(}U\setminus K)$\emph{ ~and
~}$A\subset\mbox{Am}(U\setminus K)$\emph{,}\\
\textbf{\hspace*{5.8mm}}$d_{_{H}}-$\emph{open}\textbf{ }\emph{in
}~$\mbox{Per}(M)$\emph{ ~and in ~}$\mbox{Am}(M)$, \emph{~respectively,
and such that}:

\hspace*{10mm}$\bullet$\textbf{\hspace*{3mm}}\emph{both} \emph{~P
~and ~A }~~$d_{_{H}}-$\emph{accumulate in ~}$\mathfrak{S}(\mbox{bd\,}K)$\medskip{}

\hspace*{10mm}$\bullet$\textbf{\hspace*{3mm}}\emph{K ~is} \emph{bi-stable}
\emph{in relation to ~$P^{*}=\underset{\gamma\in P}{\bigcup}\gamma$
}~\emph{and ~$A^{*}=\underset{\varGamma\in A}{\bigcup}\varGamma$}

\hspace*{10mm}$\bullet$\textbf{\hspace*{3mm}}\emph{for any sequence
$\gamma_{_{n}}\in P$, ~}$\mbox{dist}$\emph{$(\gamma_{_{n}},K)\rightarrow0$
~implies~ }$\mbox{period}$\emph{$(\gamma_{_{n}})\rightarrow+\infty\,$}.\medskip{}

\begin{figure}[H]
\noindent \begin{centering}
\includegraphics[scale=0.22]{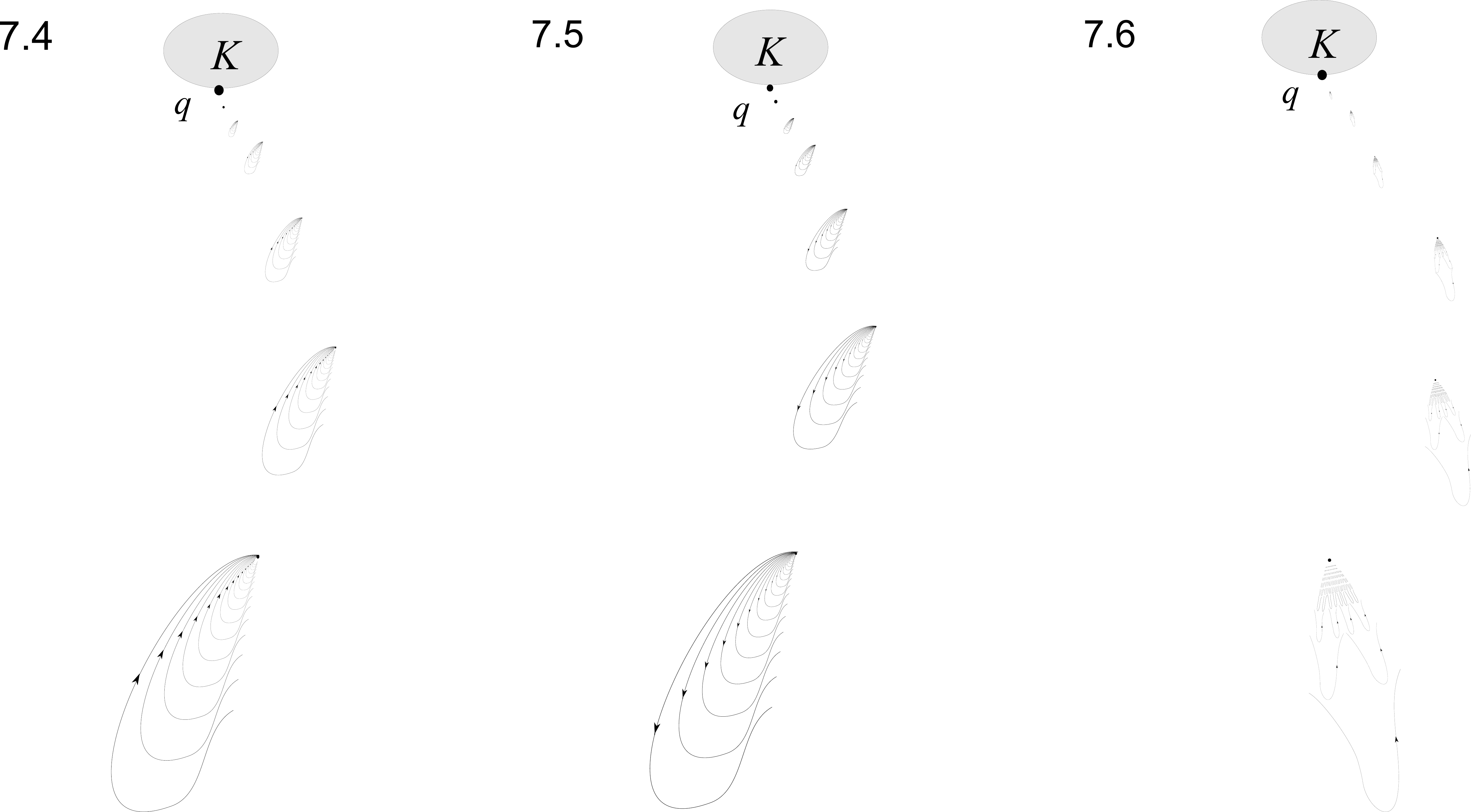}
\par\end{centering}

\caption{}
\end{figure}

\begin{figure}[H]
\noindent \begin{centering}
\includegraphics[scale=0.21]{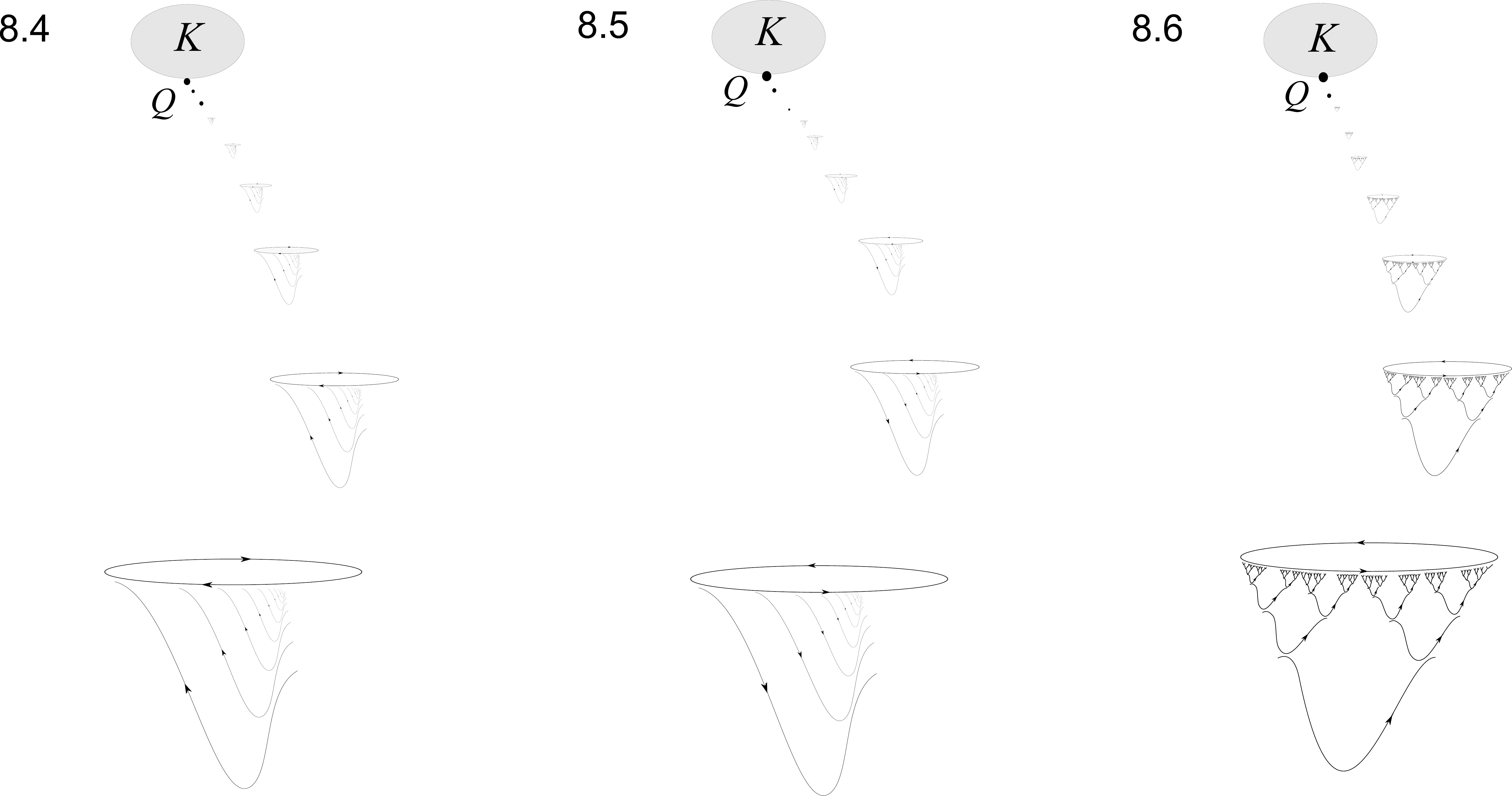}
\par\end{centering}

\caption{}
\end{figure}

Observe that, \emph{$\gamma_{_{n}}\in P$ ~}and\emph{ }$\mbox{dist}$\emph{$(\gamma_{_{n}},K)\rightarrow0$
~}together actually imply that ~$|\gamma_{_{n}}|_{_{K}}\rightarrow0$,
~since ~$P^{*}$ ~is \emph{bi-stable }in relation to ~$K$. ~On
the other hand, by the 1st point of ~\textbf{10.4} ~there is a sequence
~$P\ni\gamma_{_{n}}\overset{d_{_{H}}}{\longrightarrow}Q\in\mathfrak{S}(\mbox{bd\,}K)$.
~As ~$P$ ~is ~$\mathfrak{c}-$\emph{dense in itself,} ~the 3rd
point can thus be replaced by the following condition:\medskip{}

$\bullet$\textbf{\hspace*{3mm}}\emph{given any ~$n\geq1$, ~all
periodic orbits ~$\gamma\in P$ ~contained in a sufficiently small
neighbourhood ~$V$ ~of ~$K$ ~have period ~$>n$, ~and there
is always a continuum of these.}\medskip{}

\emph{Remark. ~}Due to the fact that, for any given ~$\epsilon>0$,
~a ~$X$-\emph{tree} ~always has a \emph{sub} $X$-\emph{tree ~}with
all its orbits contained in ~$B(X,\epsilon)\setminus X$, and analogue
fact holds for ~$X$-$\alpha shells$ ~and for ~$X$-$\omega shells$
~(see section 4), it is immediate to verify that any of the twelve
conditions ~\textbf{4}, ~\textbf{5}, ~\textbf{6}, ~\textbf{7.4}
~to ~\textbf{7.6},~\textbf{ 8.4} ~to ~\textbf{8.6} ~and ~\textbf{9.4}
~to ~\textbf{9.6 }~implies that ~\emph{for every neighbourhood
of ~$U$ ~of ~$K$, ~there is a sequence of orbits ~$\gamma_{_{n}}\subset U\setminus K$~
such that}\[
\mbox{cl}\,\gamma_{_{1}}\supsetneq\mbox{cl}\,\gamma_{_{2}}\supsetneq\cdots\cdots\supsetneq\mbox{cl}\,\gamma_{_{n}}\supsetneq\cdots\cdots\]
Thus the claim of condition~ \textbf{H}\textbf{\emph{ ~}}in section
3 is justified.\medskip{}

\begin{figure}[H]
\noindent \begin{centering}
\includegraphics[scale=0.2]{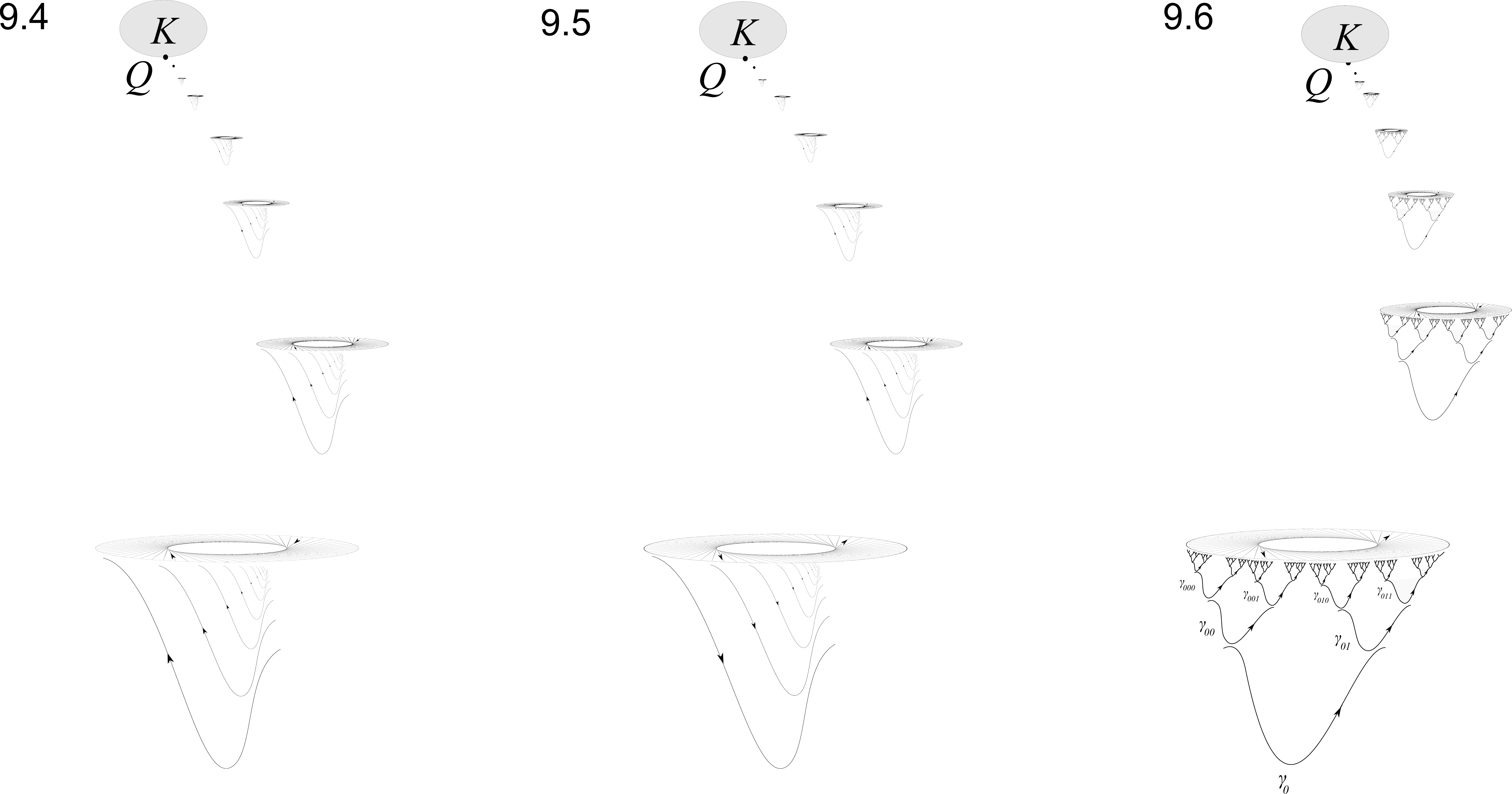}
\par\end{centering}

\caption{{\small examples 9.4 to 9.6 show sequences of 2-tori carrying irrational
linear subflows.} }
\end{figure}

\noun{interlude: ~topological structure of ~}$\mbox{CMin}(M)$\noun{.}

~~~Theorem 1 brings to light the importance of compact minimal
sets to the characterization of the possible {}``dynamical landscapes''
in the vicinity of a compact invariant proper subset of a flow. ~Obviously
there is a close relation between the dynamical behaviour of a flow
near a compact minimal set ~$X$ ~and the topological Hausdorff
structure of ~$\mbox{CMin}(M)$ ~near ~$X$. ~Actually, from Theorem
1 we easily obtain an elegant characterization of the set ~$\mbox{CMin}(M)$
~of all compact minimal sets of the flow, endowed with the Hausdorff
metric.

\emph{~}~~Let ~$M$ ~be a locally compact, connected metric space
with a ~$C^{\,0}$ ~flow.~ Consider the following seven propositions,
where the variable ~$X$ ~now takes values in the set ~$\mbox{CMin}(M)$
~of all \emph{compact minimal subsets} of the flow:

\textbf{\hspace*{4mm}}\emph{$1.X$} \textbf{\hspace*{5.5mm}}$X$~
\emph{is an attractor.}

\textbf{\hspace*{4mm}}$2.X$\textbf{\hspace*{6mm}}$X$~ \emph{is
a} \emph{repeller.}

\textbf{\hspace*{4mm}}$3.X$\textbf{\hspace*{6mm}}$X$~ \emph{is}
\emph{an isolated minimal set and stagnant. }

\textbf{\hspace*{4mm}}$4.X$\textbf{\hspace*{6mm}}$X$~ \emph{is}
\emph{an isolated minimal set} \emph{and there is a} ~$X\mbox{-}\alpha\, shell$.

\textbf{\hspace*{4mm}}$5.X$\textbf{\hspace*{6mm}}$X$~ \emph{is
an isolated minimal set} \emph{and there is a} ~$X\mbox{-}\omega\, shell$.

\textbf{\hspace*{4mm}}$6.X$\textbf{\hspace*{6mm}}$X$~ \emph{is}
\emph{an isolated minimal set and there is a} ~$X\mbox{-}tree$.

\textbf{\hspace*{4mm}}$10.X$\textbf{\hspace*{5mm}}\emph{there is
an ~$\epsilon>0$ ~such that the compact minimal sets contained
in}\\
\textbf{\hspace*{17mm}}\emph{$B(X,\,\epsilon)$ ~form a ~$\mathfrak{c}-$dense
in itself subset of }~$\mbox{CMin}(M)$\emph{.}\smallskip{}

or equivalently (by lemma 6, section 5),\smallskip{}

\textbf{\hspace*{4mm}}$10'.X$\textbf{\hspace*{5mm}}\emph{there
is an} ~$\epsilon>0$ \emph{~such that every neighbourhood} ~$U_{_{Y}}\subset M$
\emph{~of each}\textbf{}\\
\textbf{\hspace*{17mm}} $Y\in\mbox{CMin}\big(B(X,\,\epsilon)\big)$\emph{
~contains a continuum of compact minimal sets.}

If ~$X$ ~satisfies ~$10.X$, ~then as ~$X$ ~is itself a compact
minimal set, every neighbourhood of ~$X$ ~actually contains a \emph{continuum
}of compact minimal sets.

\emph{Important remark: ~}Recall that,\emph{ }the definition of the
Hausdorff metric and lemma 4 (section 5), together imply that

\textbf{\hspace*{4mm}}{}``$X$ ~\emph{is an isolated compact minimal
set iff it is ~$d_{_{H}}-$isolated in ~}$\mbox{CMin}(M)$.''

(by definition (section 2), ~$X$ ~is an \emph{isolated minimal
(set)} if ~for some ~$U\in\mathcal{N}_{_{X}}$, ~$U\setminus X$
~contains no minimal set of the flow). 

\textbf{Definition. }\textbf{\emph{~}}\emph{denote by ~$\mathfrak{M}_{_{i}}$,
~$1\leq i\leq6$ ~or ~$i=10$ ~the set of all ~}$X\in\mbox{CMin}(M)$\emph{
~satisfying condition ~$i.X$ ~and by ~$\mathfrak{M}_{_{i-j}}$,
~$1\leq i<j\leq6$ ~the set of all ~}$X\in\mbox{CMin}(M)$\emph{
~satisfying (at least) one of condition ~$i.X$ ~to ~$j.X$.\medskip{}
}

\textbf{Theorem 2.} ~\emph{Let ~M ~be a locally compact, connected
metric space with a ~$C^{\,0}$ ~flow. ~Then:}

1)\textbf{\hspace*{4mm}}\emph{$\mathfrak{M}_{_{1-6}}$ ~is the set
of isolated compact minimal sets and thus a countable,}\\
\textbf{\hspace*{7.5mm}}\emph{ $d_{_{H}}-$open subset of} ~$\mbox{CMin}(M)$.

2)\textbf{\hspace*{4mm}}\emph{$\mathfrak{M}_{_{10}}$ ~is a ~$d_{_{H}}-$open
and ~$\mathfrak{c}-$dense in itself subset of ~}$\mbox{CMin}(M)$.
\emph{~It is either}\\
\textbf{\hspace*{7.5mm}}\emph{empty or has the cardinal of the
continuum.}

3)\textbf{\hspace*{4mm}}\emph{$\mathfrak{M}_{_{1-6}}$ ~is ~$d_{_{H}}-$dense
in ~}$\mbox{CMin}(M)\setminus\mathfrak{M}_{_{10}}$.\emph{\medskip{}
}

\emph{Proof. ~}In first place, note the following trivial fact that
will be implicitly used in several instances bellow: if ~$X$ ~is
a compact minimal set and ~$Y_{_{n}}$ ~is a sequence of compact
minimal sets ~$d_{_{H}}-$converging to ~$Q\in\mathfrak{S}(X)$
~then ~$Q=X$.

\emph{~}1) ~Clearly every ~$X\in\mathfrak{M}_{_{1-6}}$ ~is an
isolated compact minimal set; on the other hand, by Theorem 1, any
compact minimal set ~$X$ ~satisfying none of the six conditions
~$1.X$ ~to ~$6.X$ ~is \emph{not} an isolated compact minimal
set,%
\footnote{note that if ~$X\in\mbox{CMin}(M)$ ~satisfies none of conditions
~$1.X$ ~to ~$6.X$, ~then ~$X$ ~is necessarily a proper subset
of ~$M,$ ~as ~$M$ ~compact implies ~$M$ ~is both an attractor
and a repeller in the flow, hence Theorem 1 can be applied.%
} ~hence ~$\mathfrak{M}_{_{1-6}}$ ~is the set of isolated compact
minimal sets of the flow. ~By a remark above, ~every ~$X\in\mathfrak{M}_{_{1-6}}$
~is thus ~$d_{_{H}}-$isolated in $\mbox{CMin}(M)$, ~hence ~$\mathfrak{M}_{_{1-6}}$
~is ~$d_{_{H}}-$open in ~$\mbox{CMin}(M)$. ~~$\mathfrak{M}_{_{1-6}}$
~is countable since ~it is a ~$d_{_{H}}-$discrete and separable
metric space.%
\footnote{recall (section 8.2) that ~$\mbox{C}(M)\supset\mathfrak{M}_{_{1-6}}$
~is ~$d_{_{H}}-$separable.%
}

2) ~If ~$X\in\mathfrak{M}_{_{10}}$ ~then there is an ~$\epsilon_{_{X}}>0$
~such that ~$\mbox{CMin}\big(B(X,\,\epsilon_{_{X}})\big)$ ~is
a ~$\mathfrak{c}-$\emph{dense in itself}, ~$d_{_{H}}-$open subset
of ~$\mbox{CMin}(M)$ ~(lemma 8, section 5). ~It is immediate to
verify that ~$\mbox{CMin}\big(B(X,\,\epsilon_{_{X}})\big)\subset\mathfrak{M}_{_{10}}$:
~if ~$Y\in\mbox{CMin}\big(B(X,\,\epsilon_{_{X}})\big)$ ~then ~$Y\subset B(X,\,\epsilon_{_{X}})$
~and since ~$Y$ ~is compact, ~$\epsilon_{_{X}}-|Y|_{_{X}}>0$,
~thus ~$B(Y,\,\epsilon_{_{X}}-|Y|_{_{X}})\subset B(X,\,\epsilon_{_{X}})$,
~hence as ~$B(Y,\,\epsilon_{_{X}}-|Y|_{_{X}})$ ~is open in $M$,
by lemma 8, ~$\mbox{CMin}\big(B(Y,\,\epsilon_{_{X}}-|Y|_{_{X}})\big)$
~is a ~$d_{_{H}}-$open subset of ~$\mbox{CMin}\big(B(X,\,\epsilon_{_{X}})\big)$
~and hence ~$\mathfrak{c}-$\emph{dense in itself. ~}Thus ~$Y\in\mathfrak{M}_{_{10}}$
~and therefore~ $\mathfrak{M}_{_{10}}=\underset{X\in\mathfrak{M}_{_{10}}}{\bigcup}\mbox{CMin}\big(B(X,\,\epsilon_{_{X}})\big)$
~is a ~$\mathfrak{c}-$\emph{dense in itself},\emph{ }~$d_{_{H}}-$open
subset of ~$\mbox{CMin}(M)$. ~Since the phase space ~$M$ ~is
separable, there is at most a \emph{continuum} of compact minimal
sets in the flow, ~therefore ~$\mathfrak{M}_{_{10}}$ ~is either
empty or has cardinal~ $\mathfrak{c}$, ~since it is ~$\mathfrak{c}-$\emph{dense
in itself.}

3) ~We claim that ~$\mathfrak{M}_{_{10}}$ ~is the set of compact
minimal sets of the flow satisfying condition ~\textbf{10} ~of Theorem~1
~(taking ~$K:=X\,$): ~if ~$X\in\mathfrak{M}_{_{10}}$ ~then
it clearly satisfies none of the 24 conditions ~\textbf{1 ~}to ~\textbf{9.6}
~(since ~$X$~ is not the ~$d_{_{H}}-$limit of a sequence of
isolated compact minimal sets), ~hence by the same theorem, it must
satisfy condition ~\textbf{10}. ~On the other hand, if ~$X$ ~satisfies
condition ~\textbf{10}, then there is an open ~$U\in\mathcal{N}_{_{X}}$
~such that ~$\mbox{CMin}(U\setminus X)$ ~is ~$\mathfrak{c}-$\emph{dense
in itself}, ~hence taking ~$\epsilon>0$ ~such that ~$B(\, X,\,\epsilon)\subset U$,
~it follows (by lemma 8) that ~$\mbox{CMin}\big(B(X,\,\epsilon)\setminus X\,\big)$
~is a ~$d_{_{H}}-$open subset of ~$\mbox{CMin}(U\setminus X)$
~and hence ~$\mathfrak{c}-$\emph{dense in itself}. ~Also by Theorem
1, at least one of the 4 conditions ~\textbf{10.1~} to ~\textbf{10.4}
~is satisfied, hence ~$X$ ~is the ~$d_{_{H}}-$limit of a sequence
of compact minimal sets ~$Y_{_{n}}\in\mbox{CMin}\big(B(X,\,\epsilon)\setminus X\,\big)$.
~As ~$\mbox{CMin}\big(B(X,\,\epsilon)\setminus X\,\big)$ ~is ~$\mathfrak{c}-$\emph{dense
in itself},\emph{ }this implies that in every neighbourhood of ~$X$
~there is a \emph{continuum }of compact minimal sets contained in
~ $B(X,\,\epsilon)$. ~Hence ~$\mbox{CMin}\big(B(X,\,\epsilon)\big)$
~is ~$\mathfrak{c}-$\emph{dense in itself}, ~thus ~$X$ ~satisfies
condition ~$10.X$. ~Therefore, by Theorem 1, if ~$X\in\mbox{CMin}(M)\setminus\mathfrak{M}_{_{10}}$
~then ~$X$ ~must satisfy (at least) one of the 24 conditions ~\textbf{1
}to \textbf{9.6 ~}and this clearly implies ~$X\in\mbox{cl\ensuremath{_{_{H}}}}\mathfrak{M}_{_{1-6}}$
~since ~$\mathfrak{M}_{_{i}}$, ~$1\leq i\leq6$ ~is the set of
compact minimal sets of the flow satisfying condition $1\leq\mathbf{i}\leq6$
~of Theorem 1.\hfill{}$\blacksquare$\emph{ \medskip{}
}

~~~Note, however, that ~$10.X$ ~is indeed a very strong condition,
essentially due to its ~$d_{_{H}}-$\emph{openness:} even when ~$\mbox{CMin}(M)\setminus\mathfrak{M}_{_{1-6}}$
~is nonvoid and ~\emph{$\mathfrak{c}-$dense in itself}, ~it can
happen that ~$\mathfrak{M}_{_{7}}$ ~is empty, since it is still
possible that ~$\mathfrak{M}_{_{1-6}}$ ~is ~$d_{_{H}}-$dense
in the whole ~$\mbox{CMin}(M)$ ~(simple examples of ~$C^{\,\infty}$
~flows exhibiting this phenomenon already occur on ~$\mathbb{S}^{1}$
~and ~$\mathbb{S}^{2}\big)$. ~However, the next result shows that
a nonvoid ~$\mathfrak{c}-$\emph{dense in itself} ~set of compact
minimal sets always occurs, whenever there are uncountably many compact
minimal sets in the flow. ~More precisely, if ~$\mbox{CMin}(M)$
~is \emph{uncountable}, then removing from this set a suitable \emph{countable}
(possibly empty) set we obtain a \emph{nonvoid} ~$\mathfrak{c}-$\emph{dense
in itself} ~set of compact minimal sets. ~This decomposition theorem
is analogous to the celebrated \emph{Cantor-Bendixson} \emph{Theorem
}for\emph{ }separable, complete metric spaces (\emph{Polish spaces})\emph{.
~}Note, however, that although ~$d_{_{H}}-$separable,%
\footnote{since ~$\mbox{C}(M)\supset\mbox{CMin}(M)$ ~is, see section 8.2.%
} in general, ~$\mbox{CMin}(M)$ ~is neither ~$d_{_{H}}-$complete
nor ~$d_{_{H}}-$locally compact. ~Also observe that, since there
is at most a \emph{continuum }of compact minimal sets in the flow
(see section 2), the above result implies that ~$\mbox{CMin}(M)$
~obeys, in a certain sense, to the \emph{Continuum Hypothesis}: its
cardinal is either finite (possibly null), denumerable ~$(\aleph_{_{0}})$
~or the \emph{continuum} ~$\mathfrak{c}=2^{\,\aleph_{_{0}}}.$\emph{\medskip{}
}

\textbf{Theorem 3.} ~\emph{Let ~$\theta$ ~be a ~$C^{\,0}$ ~flow
on a locally compact, separable metric space ~$M$,~ displaying
uncountably many compact minimal sets. ~Then there is a countable
(possibly empty) set ~}$\mathfrak{I}\subset\mbox{CMin}(M)$\emph{
~such that:}

\textbf{\hspace*{2mm}}I.\textbf{\hspace*{4.5mm}}$\mathfrak{D}:=\mbox{CMin}(M)\setminus\mathfrak{I}$
~\emph{is a} ~$\mathfrak{c}-$\emph{dense in itself} \emph{and ~$d_{_{H}}-$closed
subset of ~}$\mbox{CMin}(M)$\emph{, ~having the cardinal of the
continuum.}

\textbf{\hspace*{2mm}}II.\textbf{\hspace*{4mm}}$\mathfrak{I}$ \emph{~is
the set of all} ~$X\in\mbox{CMin}(M)$ \emph{~having a neighbourhood
containing only countably many (possibly one) compact minimal sets,
hence} ~$\mathfrak{D}$ \emph{~is the largest} ~$\mathfrak{c}-$\emph{dense
in itself} \emph{~subset of }~$\mbox{CMin}(M).$\emph{\medskip{}
}

Hence, \emph{~if ~}$\mbox{CMin}(M)$\emph{ ~is uncountable, then
all but a countable number of compact minimal sets of the flow have
a continuum of compact minimal sets on each of their neighbourhoods
}or, equivalently, ~$\mbox{CMin}(M)$\emph{ ~is the union of a countable
(possibly empty) set and a ~$\mathfrak{c}-$dense in itself set.} 

The proof uses in an essential way a \emph{{}``}Cantor's ternary
set - like'' construction that constitutes the core of the proof
of lemma 7 (section 8.1).

\emph{Proof.} ~Suppose\emph{ ~}$\mbox{CMin}(M)$\emph{ ~}is uncountable.
~Let ~$\mathfrak{I}$ ~be the set of all ~$X\in\mbox{CMin}(M)$~
for which there is an ~$\epsilon>0$ ~such that\[
B_{_{H}}(X,\,\epsilon)\mbox{\,\ \emph{contains} \emph{only} \emph{countably} \emph{many} \emph{(possibly} \emph{only} \emph{one)} \emph{compact} \emph{minimals.}}\]
(by lemma 5 this is actually equivalent to: there is a ~$\delta>0$~
such that there are only countably many compact minimal sets contained
in ~$B(X,\,\delta)\subset M\,\big)$.~ For each ~$X\in\mathfrak{I}$
~define\[
\epsilon_{_{X}}:=\mbox{sup}\big\{\epsilon>0:\,\, B_{_{H}}(X,\,\epsilon)\,\cap\,\mbox{CMin}(M)\mbox{ \,\ is countable\,\ensuremath{\big\}}}\]
Note that ~$0<\epsilon_{_{X}}<+\infty$ ~since ~$\mbox{CMin}(M)$
~is, by hypothesis, uncountable and ~$\mbox{CMin}(M)=\underset{n\geq1}{\bigcup}\big(B_{_{H}}(X,\, n)\,\cap\,\mbox{CMin}(M)\big)$.
~Also, observe that for each ~$X\in\mathfrak{I}$, ~$B_{_{H}}(X,\,\epsilon_{_{X}})\,\cap\,\mbox{CMin}(M)$
~is a countable ~$d_{_{H}}-$open subset of ~$\mbox{CMin}(M)$
~and thus is contained in ~$\mathfrak{I}$ ~(by definition of ~$\epsilon_{_{X}}$,
~$B_{_{H}}\big(X,\,\epsilon_{_{X}}(1-1/n)\big)\,\cap\,\mbox{CMin}(M)$
~is countable for each ~$n\in\mathbb{N},$~ hence ~$B_{_{H}}(X,\,\epsilon_{_{X}})\,\cap\,\mbox{CMin}(M)=\underset{n\geq1}{\bigcup}\big(B_{_{H}}\big(X,\,(1-1/n)\epsilon_{_{X}}\big)\,\cap\,\mbox{CMin}(M)\big)$\emph{
~}is countable).

Claim 1: ~$\mathfrak{I}$ \emph{~is countable:} ~assume ~$\mathfrak{I}$
~is infinite (otherwise the claim is proved). ~$\mbox{CMin}(M)$
~is ~$d_{_{H}}-$separable ~(since ~$\mbox{C}(M)\supset\mbox{CMin}(M)$
~is),%
\footnote{see section 8.2%
} ~hence ~$\mathfrak{I}$ ~ is ~$d_{_{H}}-$separable. ~Let ~$I=\{X_{_{1}},\, X_{_{2}},\ldots,X_{_{n}},\ldots\}$
~be a denumerable ~$d_{_{H}}-$dense subset of ~$\mathfrak{I}$.
~We claim that \[
\mathfrak{I}=\underset{n\geq1}{\bigcup}\Big(B_{_{H}}\big(X_{_{n}},\,\epsilon_{_{X_{_{n}}}}\big)\,\cap\,\mbox{CMin}(M)\Big)\]
therefore proving the countability of ~$\mathfrak{I}$. ~The inclusion
~$\supset$ ~is already established by a remark above. ~To prove
~$\subset$ ~observe that given any ~$Y\in\mathfrak{I}$, ~there
is a ~$X_{_{n}}\in I$ ~such that ~$X_{_{n}}\in B_{_{H}}(Y,\,\epsilon_{_{Y}}/2)$.
~Thus, by the triangle inequality for the ~$d_{_{H}}$ ~metric,
~$B_{_{H}}(X_{_{n}},\,\epsilon_{_{Y}}/2)\subset B_{_{H}}(Y,\,\epsilon_{_{Y}})$,
~hence ~$B_{_{H}}(X_{_{n}},\,\epsilon_{_{Y}}/2)\,\cap\,\mbox{CMin}(M)\subset B_{_{H}}(Y,\,\epsilon_{_{Y}})\,\cap\,\mbox{CMin}(M)$
~is countable. ~Therefore ~$\epsilon_{_{X_{_{n}}}}\geq\epsilon_{_{Y}}/2$.
~Hence ~$Y\in B_{_{H}}\big(X_{_{n}},\,\epsilon_{_{X_{_{n}}}}\big)\,\cap\,\mbox{CMin}(M)$
~since ~$d_{_{H}}(X_{_{n}},\, Y)<\epsilon_{_{Y}}/2$. ~The claim
is proved. ~Note that the identity above also proves that~ $\mathfrak{I}$
~is ~$d_{_{H}}-$open in ~$\mbox{CMin}(M)$.

Now let ~$\mathfrak{D}:=\mbox{CMin}(M)\setminus\mathfrak{I}$. ~Clearly
~$\mathfrak{D}$~ is nonvoid since ~$\mbox{CMin}(M)$ ~is uncountable
(by hypothesis) ~and ~$\mathfrak{I}$ ~is countable. ~Actually,
by definition of ~$\mathfrak{I}$, ~given any ~$X\in\mathfrak{D}$
~and ~$\epsilon>0$, ~$B_{_{H}}(X,\,\epsilon)\,\cap\,\mbox{CMin}(M)$
~is uncountable\emph{, }hence ~$B_{_{H}}(X,\,\epsilon)\,\cap\,\mathfrak{D}$
~is also ~uncountable and in particular, ~$\mathfrak{D}$ ~is
a nonvoid ~$d_{_{H}}-$\emph{dense in itself} ~subset of ~$\mbox{CMin}(M)$.

Claim 2: $\mathfrak{D}$ ~\emph{is} ~$\mathfrak{c}-$\emph{dense
in itself:} ~As ~$\mathfrak{D}\subset\mbox{CMin}(M)$, ~in virtue
of lemma 6 (section 5), we need only to prove that given any ~$X\in\mathfrak{D}$
~and ~$\epsilon>0$, ~there is a \emph{continuum }of compact minimal
sets ~$Y\in\mathfrak{D}$ ~contained in ~$B(X,\,\epsilon)\subset M$.
~Taking ~$\epsilon$ ~sufficiently small we may assume ~$B[\, X,\,\epsilon]$
~is compact ~($X$ ~is compact and ~$M$ ~is locally compact).
~Let ~$A:=B(X,\,\epsilon)$, ~$\varLambda_{_{0}}:=X$\emph{ ~}and
$\epsilon_{_{0}}:=\epsilon/2$. ~Now since ~$\mathfrak{D}$ ~is
~$d_{_{H}}-$\emph{dense in itself}, ~we may carry the construction
of the proof of lemma 7 (section 8.1) within ~$\mathfrak{D}(A)=\big\{ Z\in\mathfrak{D}:\,\, Z\subset A\big\}$
~\emph{i.e.} we may select each ~$\varLambda_{_{a}}$, ~$a\in$$\mathcal{F}$~
in ~$\mathfrak{D}(A)$ ~instead of in ~$\mbox{CMin}(A)$. ~As
in the proof of lemma 7 we get a \emph{continuum }of ~$d_{_{H}}-$Cauchy
sequences, ~$d_{_{H}}-$converging to a \emph{continuum }of mutually
disjoint, nonvoid, compact invariant sets contained in ~$A$, ~therefore
proving the existence of a \emph{continuum }of compact minimal sets
contained in in this open set ~(as each ~$K\in\mbox{Ci}(A)$ ~contains
at least one compact minimal set). ~Now since ~$\mathfrak{I}$ ~is
countable, a \emph{continuum}~of these compact minimal sets ~$\varGamma\in\mbox{CMin}(A)$
~actually belongs to ~$\mathfrak{D}=\mbox{CMin}(M)\setminus\mathfrak{I}$.
~Therefore ~$\mathfrak{D}$ ~is ~$\mathfrak{c}-$\emph{dense in
itself}. ~ $\mathfrak{D}$ ~is ~$d_{_{H}}-$closed in ~$\mbox{CMin}(M)$
~since ~$\mathfrak{I}$ ~is ~$d_{_{H}}-$open in the same set.\hfill{}$\blacksquare$

It is simple to see that the set ~$E$ ~of \emph{equilibria }satisfies
the following stronger analogue property to that expressed on Theorem
3: 

{}``~\emph{if ~$E$ ~is uncountable, then ~$E$ ~is the union
of a countable set and a perfect subset ~$\mathfrak{E}$~ of ~$M$,
~with the cardinal of the continuum. ~For each ~$z\in\mathfrak{E}$
~and ~$\epsilon>0$ ~there is an embedding $h$ ~of ~Cantor's
ternary set into ~$B(z,\,\epsilon)\,\cap\,\mathfrak{E}$ ~with }~$z\in\mbox{im}\, h$\emph{.
~Hence ~$\mathfrak{E}\subset E$ ~is a ~$\mathfrak{c}-$dense
in itself closed subset of ~$M$.''}%
\footnote{This follows immediately from the following observation: ~since the
phase space ~$M$ ~is locally compact and separable it can be endowed
with an equivalent \emph{boundedly compact} metric (on which every
closed bounded set is compact, see \emph{e.g.} \cite{lima}, p.278),
thus becoming a complete, separable metric space (\emph{Polish space)}\@.
~$E$ ~is closed in $M$, ~hence is also a \emph{Polish space }in
this equivalent metric, and this proposition is well known to hold
on such spaces (see \emph{e.g. }\cite{levy}\emph{, }chap. VII.2).%
}

The question now arises as whether the corresponding propositions
analogue to Theorem 3, for the set ~$\mbox{Per}(M)$ ~of all \emph{periodic
orbits }and for the set ~$\mbox{Am}(M)$ ~of all \emph{compact aperiodic
minimal sets }of the flow, also hold\smallskip{}

~~I.~~\emph{If ~}$\mbox{Per}(M)$ \emph{~is uncountable then
all but a countable number of periodic orbits of the flow have a continuum
of periodic orbits on each of their neighbourhoods. }

~~II.~\emph{~If ~}$\mbox{Am}(M)$\emph{ ~is uncountable then
all but a countable number of compact aperiodic minimal sets of the
flow have a continuum of compact aperiodic minimal sets on each of
their neighbourhoods.}

As already mentioned in the introduction, it is, in a certain sense,%
\footnote{\emph{i.e.} working within \emph{Zermelo-Fraenkel Set Theory }and
provided this {}``standard'' axiomatic\emph{ }is consistent.%
}~ useless to look for counterexamples to any of these two propositions
within {}``standard'' Dynamical Systems Theory: ~both I. and II.
are provable in \emph{~ZFC~ }set theory under the additional assumption
of the \emph{Continuum Hypothesis ~CH. ~}Hence each turns out to
be either demonstrable in \emph{~ZFC ~}or independent of this standard
axiomatic (due to Gödel's result, 1938).\emph{ ~}The proof that ~$CH\Longrightarrow\mbox{I}\,\wedge\,\mbox{II}$
~is simple and actually depends only on the fact ~$\mbox{Per}(M)$
~is a separable (Hausdorff) metric space: ~I. and II. are particular
cases of the following proposition, which is equivalent to the \emph{Continuum
Hypothesis:\medskip{}
}

$\mathfrak{c}-$\emph{Denseness Hypothesis }(\emph{$\mathfrak{c}DH\,$}):
~\emph{If ~$L$ ~is an uncountable separable metric space, then
a ~$\mathfrak{c}-$dense in itself set is obtained removing from
~$L$ ~a suitable countable set (possibly empty).}%
\footnote{here ~$X\subset M$ ~is ~$\mathfrak{c}-$\emph{dense in itself
}~means, as for the ~$d_{_{H}}$ ~metric, ~that for every ~$x\in X$
~and ~$\epsilon>0$, ~~$B(x,\,\epsilon)\,\cap\, X$ ~has the
cardinal ~$\mathfrak{c}$ ~of the \emph{continuum.}%
}

As we could locate no reference for the equivalence ~$CH\Longleftrightarrow\mathfrak{c}DH$,
~a short proof is included in section 8.3, for the sake of completeness.

\noun{proof of theorem 1.}

\emph{Synopsis:} ~Assume neither \textbf{1} nor \textbf{2} hold.
~A) ~If ~$K$ ~is isolated from minimal sets\emph{ }then it is
shown that at least one of conditions \textbf{3} to \textbf{6 }necessarily
holds. ~B) ~If $K$ ~is not isolated from minimals, then we consider
two possible cases: B.1) if for every neighbourhood $U$ ~of ~$K,$
~$U\setminus K$ ~contains a compact minimal set of ~$(M,\,\theta)$
~satisfying (at least) one of the six conditions ~1.\emph{X} to
6.\emph{X, ~}then it is proved that at least one of the eighteen
cases \textbf{7.1} to \textbf{9.6} holds; B.2)~if the contrary is
true, then there is a neighbourhood ~$U$ ~of ~$K$ ~such that
~$\mbox{CMin}(U\setminus K)$ is a $\mathfrak{c}-$\emph{dense in
itself},\emph{ }$d_{_{H}}-$open subset of ~$\mbox{CMin}(M)$,\emph{
}~$d_{_{H}}-$accumulating in ~$\mathfrak{S}(K)$ ~and at least
one of the four conditions \textbf{10.1} to \textbf{10.4} necessarily
holds.\emph{\bigskip{}
}

\emph{Proof of Theorem 1.} ~It easily seen that condition \textbf{~1~}
excludes the remaining 27 conditions and the same holds with condition
~\textbf{2}.%
\footnote{Actually, if $K$~ is an \emph{~attractor }~then ~$B^{+}(K)\setminus K\neq\emptyset$
~and ~~$x\in B^{+}(K)\setminus K\implies\alpha(x)\subset\mbox{bd}\, B^{+}(K)\subset M\setminus B^{+}(K)$
~(see footnotes 13 and 25).~ A time-symmetric fact holds if ~$K$~is
a~ \emph{repeller.} ~This immediately implies that if ~\textbf{1}
(resp. \textbf{2}) holds , then\textbf{ }none of the remaining 27
conditions can take place.%
}~ \emph{Assume, throughout the remaining of this proof, that neither}
~\textbf{1} ~\emph{nor} ~\textbf{2} ~\emph{holds}.~~In this
situation we distinguish the two possible cases:\medskip{}

A$)$~~~$K$ \emph{~is} \emph{isolated from minimals }~\emph{i.e.}
\emph{~for some} ~$U\in\mathcal{\mathcal{N}}_{_{K}}$, ~~$\mbox{CMin(}U\setminus K)=\emptyset.$

B$)$~~~\emph{for every} ~$V\in\mathcal{N}_{_{K}},$ ~~$\mbox{CMin}(V\setminus K)\neq\emptyset$.

We recall an important elementary that will be implicitly used in
several instances bellow:\emph{ on a locally compact metric space,
every sufficiently small neighbourhood of a compact set has compact
closure, and thus may only contain} \emph{compact minimal sets}.

\medskip{}

Case A): 

Since ~$K$ ~is compact and ~$M$ ~is locally compact. we may
assume, without loss of generality, ~that ~$U$ ~is compact.~
Then for any ~$z\in U,$\[
\mathcal{O}^{+}(z)\subset U\implies\omega(z)\,\cap\, K\neq\emptyset\implies z\in A^{+}(K)\,\sqcup\, B^{+}(K)\]
\[
\mathcal{O}^{-}(z)\subset U\implies\alpha(z)\,\cap\, K\neq\emptyset\implies z\in A^{-}(K)\,\sqcup\, B^{-}(K)\]

since otherwise we would have ~$\mbox{CMin(}U\setminus K)\neq\emptyset$
~$\big(\,$clearly ~$\mathcal{O}^{+}(z)\subset U\implies$$\mbox{cl\,}\mathcal{O}^{+}(z)=\mathcal{O}^{+}(z)\,\cup\,\omega(z)\subset U,$
~thus if ~$\omega(z)\,\cap\, K=\emptyset$ ~then ~$\omega(z)\subset U\setminus K.$~
But ~$\omega(z)$ ~is a nonvoid, compact, invariant set (since ~$\mathcal{O}^{+}(z)\subset U$
~is compact, hence it must contain a compact minimal set of ~$(M,\theta),$
~in contradiction with ~$\mbox{CMin}(U\setminus K)=\emptyset$.~
If ~$\mathcal{O}^{-}(z)\subset U,\,$ then assuming ~$\alpha(z)\,\cap\, K=\emptyset$
~we arrive at the same contradiction$\big)$. ~Suppose now that
condition ~\textbf{3} ~does not hold.~ Since ~$K$ ~is \emph{isolated
from minimals }it follows that \emph{ }~$K$ ~is \emph{non-stagnant,
}therefore for every orbit\emph{ ~}$\mathcal{O}(z)\subset U\setminus K,$\smallskip{}

~~~either~~~0.~~~~$z\in A^{-}(K)\,\cap\, B^{+}(K)$

~~~or~~~~~~~~~I.~~~~\foreignlanguage{english}{$z\in B^{-}(K)\,\cap\, A^{+}(K)$}

\selectlanguage{english}%
~~~or~~~~~~~~~II.~~~$z\in A^{-}(K)\,\cap\, A^{+}(K)$

\selectlanguage{british}%
Depending on which of these three cases is satisfied, we say ~$\mathcal{O}(z)$
~is an orbit of ~\emph{type ~}0, ~I ~or~ II. ~More generally,
the fact ~$K$ ~is \emph{non-stagnant} implies that orbits of \emph{type
~}0 ~and ~I ~cannot coexist in ~$U\setminus K$. ~This implies
that exactly one of the following three conditions holds:\medskip{}

i)~~$\,\,$there is an orbit ~$\mathcal{O}(x)\subset U\setminus K$
such that $\big(\mbox{cl\,}\mathcal{O}(x)\big)\setminus K$ contains
only orbits\\
\hspace*{5mm}of \emph{type} 0

ii)~~there is an orbit ~$\mathcal{O}(y)\subset U\setminus K$ such
that $\big(\mbox{cl\,}\mathcal{O}(y)\big)\setminus K$ contains only
orbits\\
\hspace*{6mm}of \emph{type} I

iii)~for every orbit ~$\mathcal{O}(z)\subset U\setminus K$, ~$\,\big(\mbox{cl\,}\mathcal{O}(z)\big)\setminus K$
~contains an orbit of \emph{type} II

\medskip{}

$\big($observe that ~$\mathcal{O}(w)\subset U\setminus K\implies\big(\mbox{cl\,}\mathcal{O}(w)\big)\setminus K\subset U\setminus K\,\big).$
~~We claim that 

a)\hspace*{6mm}i)~~~~implies there is a \emph{~K-$\alpha\, shell$}

b)\hspace*{6mm}ii)~~~implies there is a ~$K\mbox{-}\omega\, shell$

c)\hspace*{6mm}iii)~~implies there is a ~$K\mbox{-}tree$

Suppose there is an orbit ~$\mathcal{O}(y)$~ satisfying condition
ii). ~Since ~$\mathcal{O}(y)$ ~is of \emph{type} I, ~by Lemma
3 (recall that ~$K$ ~is, by hypothesis, non-stagnant), given any
neighbourhood ~$V$ ~of ~$K$, ~there is a ~$p\in\omega(y)\setminus K\subset\mbox{cl\,}\mathcal{O}(y)$
~such that ~$\mathcal{O}(p)\subset V\setminus K$. ~Clearly ~$\mathcal{O}(p)$
~is also of \emph{type }I ~since ~$\mathcal{O}(p)\subset\mbox{cl\,}\mathcal{O}(y)$.
~The existence of a ~$K\mbox{-}\omega shell$\emph{ ~}with first
orbit ~$\mathcal{O}(y)$ ~is now a straightforward inductive consequence
of Lemma 3. ~Analogously, if ~$\mathcal{O}(x)$ ~is an orbit satisfying
condition i) then there is a ~$K\mbox{-}\alpha shell$ ~with first
orbit ~$\mathcal{O}(x)$. ~We now prove c). ~Recall that by hypothesis,
~\emph{K} ~satisfies none of conditions \textbf{~1}, \textbf{~2}
and\emph{ }\textbf{~3}, ~therefore by Lemma 2, there is necessarily
an orbit ~$\mathcal{O}(z)\subset U\setminus K$.~ By iii), ~$\big(\mbox{cl}\,\mathcal{O}(z)\big)\setminus K$
~contains an orbit ~$\gamma_{_{0}}$ ~of \emph{type} II. ~We will
inductively define a map\[
\begin{array}{lll}
\psi:\mathcal{E} & \longrightarrow & \mbox{Orb}\big((\mbox{cl}\,\gamma_{_{0}})\setminus K\big)\subset\mbox{Orb}\big(U\setminus K\big)\subset\mbox{Orb}\big(M\setminus K\big)\\
\,\,\,\,\,\,\,\,\,\,\,\,\,\, a & \longmapsto & \gamma_{_{a}}\end{array}\]

so that ~$\varTheta:=\mbox{im\,}\psi$ is a \emph{K-tree.} Adopt
the following lexicographic order on ~$\mathcal{E}:$\[
0<00<01<000<001<010<011<0000<0001<0010<\cdots\cdots\]
Suppose ~$a\in\mathcal{E}$ ~is such that for all ~$\mathcal{E}\ni d<a,$
~$\gamma_{_{d}}$ ~is an already defined orbit of \emph{type }II
contained in\emph{ ~}$(\mbox{cl\,}\gamma_{_{0}})\setminus K\subset U\setminus K$.~
We define ~$\gamma_{_{a}}$:

Evidently, ~$a=bc$ ~for some ~$b\in\mathcal{E}$ ~and ~$c\in\{0,1\}.$~
By Lemma 3%
\footnote{note that ~$K\in\mbox{Ci}(M)\setminus\{M\}$, ~$K$ ~is \emph{non-stagnant}
~and by hypothesis,~$\gamma_{_{b}}$ ~is of\emph{ type} II i.e.
~$\gamma_{_{b}}\subset A^{-}(K)\,\cap\, A^{+}(K)$.%
} there is an orbit ~$\zeta_{_{bc}}$ ~such that:

~~$\bullet$~~~~$\gamma_{_{b}}\overset{_{_{c}}}{\succ\,\,\,}\zeta_{_{bc}}$\smallskip{}

~~$\bullet$~~~~$0<\big|\zeta_{_{bc}}\big|_{K}<\big|\gamma_{_{b}}\big|_{K}\big/2$

hence~ $\zeta_{_{bc}}\not\succ\gamma_{_{b}}$~ and~ $\zeta_{_{bc}}\subset(\mbox{cl\,}\gamma_{_{b}})\setminus K\subset(\mbox{cl\,}\gamma_{_{0}})\setminus K\subset U\setminus K$.
~~By hypothesis iii), ~$\big(\mbox{cl\,}\zeta_{_{bc}}\big)\setminus K$
~contains an orbit of \emph{type} II and we identify ~$\gamma_{_{a}}$
~with it.~ Clearly ~$\gamma_{_{b}}\overset{_{c}}{\succ\,\,\,}\gamma_{_{bc}}=\gamma_{_{a}}$
~for every ~$b\in\mathcal{E}$, ~$c\in\{0,1\}$ ~since ~$\gamma_{_{a}}\subset\mbox{cl\,}\zeta_{_{bc}}$
~and ~$\gamma_{_{b}}\overset{_{c}}{\succ\,\,\,}\zeta_{_{bc}}$.
~Note that ~inequality ~$|\gamma_{_{a}}|_{_{K}}=|\gamma_{_{bc}}|_{_{K}}\leq|\zeta_{_{bc}}|_{_{K}}<|\gamma_{_{b}}|/2$
~guarantees ~$\gamma_{_{bc}}\not\succ\gamma_{_{b}}$ ~for every
~$b\in\mathcal{E}$, ~$c\in\{0,1\}$ ~and ~$\big|\mbox{cl\,}\gamma_{_{v_{_{n}}}}\big|_{_{K}}\longrightarrow0$
~for every ~$v\in\mbox{E}_{_{\infty}}$.~ It is now immediate to
verify that ~$(\varTheta,\,\psi)$, ~where ~$\varTheta$ ~is the
inductively defined set ~$\big\{\gamma_{_{a}}:\, a\in\mathcal{E}\big\}$,~
is indeed a\emph{ K-tree}.

We thus conclude that \emph{case} A) \emph{implies that at least one
of conditions }~\textbf{3}, \textbf{~4}, \textbf{~5~}\emph{ or}
\textbf{~6~ }\emph{necessarily holds,} therefore if conditions ~\textbf{1
}~to ~\textbf{6}\emph{~ }(i.e. 1.\emph{K} to 6.\emph{K}) all fail
then condition ~B) holds (recall we assumed, in the beginning of
the proof, that both ~\textbf{1} ~and \textbf{~2~} are false).~~Note
that since ~$K$ ~is an arbitrary compact, invariant, proper subset
of ~$M$, the above observation is true for all ~$X\in\mbox{Ci}(M)\setminus\{M\}$\emph{~
i.e}.\medskip{}

\textbf{Lemma 13. ~}\emph{If ~X~ is a compact, invariant, proper
subset of} ~\emph{$M$ ~and all conditions ~$1.X$ }to \emph{$6.X$~}
\emph{fail, then arbitrarily near} ~$X$ \emph{~there is always
a compact minimal set} \emph{of ~$(M,\,\theta)$ ~disjoint from
~$X$ ~i.e.} \emph{for any} ~$\epsilon>0,$ ~$\mbox{CMin\ensuremath{\big(}}B(X,\,\epsilon)\setminus X\big)\neq\emptyset.$\medskip{}

Case B):~~~for every ~$V\in\mathcal{N}_{_{K}},$ ~~$\mbox{CMin}(V\setminus K)\neq\emptyset.$\medskip{}

We distinguish the two (sub)cases:

B.1)~~~for every ~$V\in\mathcal{N}_{_{K}},$ ~~$V\setminus K$~
contains a compact minimal set ~$X$ ~satisfying (at least) one
of the six conditions \emph{$1.X$} to $6.X$.

B.2)~~~there is an open ~$U\in\mathcal{N}_{_{K}}$ ~such that
no ~$X\in\mbox{CMin}(U\setminus K)$ ~satisfies any of the six conditions
\emph{$1.X$} to \emph{$6.X$. \medskip{}
}

Case B.1):

We will show that in this case at least one of the eighteen cases
\textbf{~7.1~} to \textbf{~9.6~} will necessarily hold.~ Take
~$\epsilon_{_{1}}>0$~ such that ~$U:=B[\, K,\,\epsilon_{_{1}}]$
~is compact. ~Let ~$\varLambda_{_{1}}$ ~be a compact minimal
set contained in ~$U\setminus K$, ~satisfying at least one of the
six conditions ~1.\emph{X} to 6.\emph{X.~} Since ~$K$ ~and ~$\varLambda_{_{1}}$
~are disjoint compacts, ~$B(K,\,\epsilon_{_{2}})\,\cap\,\varLambda_{_{1}}=\emptyset$~
where ~$\epsilon_{_{2}}:=\mbox{dist}(\varLambda_{_{1}},\, K)/2$.~
We may obviously define two sequences ~$\varLambda_{_{n}}\in\mbox{CMin}(U\setminus K)$~
and ~$\epsilon_{_{n}}>0$~ such that for every ~$n\geq1,$

$\bullet$~~~$\varLambda_{_{n}}$ ~satisfies at least one of the
six conditions \emph{~$1.X$~} to $6.X$

$\bullet$~~~$\varLambda_{_{n}}\subset B[\, K,\epsilon_{_{n}}]$

$\bullet$~~~$\epsilon_{_{n+1}}:=\mbox{dist}(\varLambda_{_{n}},\,\, K)/2$

Since ~$\varLambda_{_{n}}\in\mbox{CMin}(U)\subset\mathfrak{S}(U)$
~for all ~$n\geq1$~ and ~$|\varLambda_{_{n}}|_{_{K}}\longrightarrow0$,
~by Lemma 11, we may select from ~$(\varLambda_{_{n}})$~ a subsequence
$d_{_{H}}-$converging to some ~$Q\in\mathfrak{S}(K)$.%
\footnote{n$\varLambda_{_{n}}\in\mathfrak{S}(M)$, ~~$\varLambda_{_{n}}\subset B[\, K,\,\epsilon_{_{n}}]\in\mbox{C}(M)$
~and ~$\epsilon_{_{n+1}}<\epsilon_{_{n}}/2$, ~thus by Lemma 11.c),
~$\varLambda_{_{n}}\overset{d_{_{H}}}{\longrightarrow}Q\in\mathfrak{S}\big(\underset{n\geq1}{\cap}B[\, K,\,\epsilon_{_{n}}]\big)=\mathfrak{S}(K)$.%
}~~Obviously ~$Q\subset\mbox{bd\,}K$, ~since ~$\varLambda_{_{n}}\subset M\setminus K$,
~thus in fact $Q\in\mathfrak{S}(\mbox{bd\,}K)$.~ By the \emph{Cantor-Dirichlet
Principle }we may\emph{ }select from this subsequence another subsequence
consisting of \emph{compact minimal sets all belonging to the same
one of the following three classes:} \emph{equilibrium orbits }$\mbox{Eq}(M),$
\emph{periodic orbits} $\mbox{Per}(M),$ \emph{compact aperiodic  minimals}
$\mbox{Am}(M).$ ~Finally since each term of ~$(\varLambda_{_{n}}$)
~satisfies at least one of the six conditions 1.\emph{X} to 6.\emph{X,
}using the \emph{Cantor-Dirichlet Principle} again, we select from
the last obtained subsequence another subsequence such that (at least)
\emph{one of the six conditions} 1.\emph{X} \emph{to} 6.\emph{X is
satisfied by all its terms}, therefore obtaining a sequence of \emph{compact
minimal sets }contained in ~~$M\setminus K$ ~satisfying at least
one of the eighteen conditions \textbf{~7.1~} to \textbf{~9.6}.

To complete the proof of Theorem 1 we show that in case B.2) \emph{i.e}.
if ~ 

$\mathfrak{a}$)~~$\mbox{CMin}(V\setminus K)\neq\emptyset$ ~for
all ~$V\in\mathcal{N}_{_{K}}$

and

$\mathfrak{b})$~~for some open ~$U\in\mathcal{N}_{_{K}}$, no
~$X\in\mbox{CMin}(U\setminus K)$ satisfies any of the six conditions
\emph{$1.X$} to \emph{$6.X$}

then at least one of the four cases \textbf{~10.1} ~to \textbf{~10.4~}
necessarily holds.

Case B.2):

We may obviously assume, without loss of generality, that ~$U$ ~has
compact closure. ~Observe that:

$-$~~~$\mbox{CMin}(U\setminus K)$~ is $d_{_{H}}-$open in ~$\mbox{CMin}(M)$
~since ~$U\setminus K$ ~is open and ~$\mbox{CMin}(M)\subset\mbox{C}(M)$
(Lemma 8).

$-$~~~~In virtue of Lemma 13, ~$\mathfrak{b}$) implies that
for any ~$X\in\mbox{CMin}(U\setminus K)$ ~and ~$\epsilon>0,$
~$\mbox{CMin\ensuremath{\big(}}B(\, X,\,\epsilon)\setminus X\big)\neq\emptyset,$
~thus by Lemma 5, ~every ~$X\in\mbox{CMin}(U\setminus K)$ ~is
\emph{~}non\emph{ $d_{_{H}}-$}isolated\emph{~} in ~$\mbox{CMin}(M)$
~and since ~$\mbox{CMin}(U\setminus K)$ ~is ~$d_{_{H}}$-open
in ~$\mbox{CMin}(M)$ ~(lemma 8), ~it follows that every ~$X\in\mbox{CMin}(U\setminus K)$
~is non\emph{ $d_{_{H}}-$}isolated in ~$\mbox{CMin}(U\setminus K)$
~\emph{i.e. }~$\mbox{CMin}(U\setminus K)$ ~is ~$d_{_{H}}-$dense
in itself. ~By Lemma 7, ~$\mbox{CMin}(U\setminus K)$~ is in fact
~$\mathfrak{c}-$\emph{dense in itself.} 

$-$~~~Using Lemma 11 (\emph{Nested Compacts Lemma})%
\footnote{Let ~$n_{_{0}}\geq1$ ~be such that ~$B[\, K,\,1/2^{n_{_{0}}}]\subset U$
~and define ~$K_{_{n}}:=B[\, K,1/2^{n}\,]\in\mbox{C}(\mbox{cl\,}U)$
~for all ~$n\geq n_{_{0}}$. ~Select a sequence ~$\varLambda_{_{n}}\in\mbox{CMin}\big(K_{_{n}}\setminus K\big)$,
~$n\geq n_{_{0}}.$~ Note that ~$K_{_{n}}\supset K_{_{n+1}},$
~~$\underset{n\geq n_{_{0}}}{\bigcap}K_{_{n}}=K$ ~and all ~$\varLambda_{_{n}}$,
~$n\geq n_{_{0}}$ ~belong to the $d_{_{H}}-$compact ~$\mathfrak{S}(\mbox{cl\,}U)$,
~therefore ~$(\varLambda_{_{n}})$ ~has a subsequence $d_{_{H}}-$converging
to some ~$Q\in\mathfrak{S}\big(\underset{n\geq1}{\bigcap}K_{_{n}}\big)=\mathfrak{S}(K)$. %
} we infer from $\mathfrak{a}$) ~that ~$\mbox{CMin}(U\setminus K)$
~~$d_{_{H}}-$accumulates in ~$\mathfrak{S}(K)$.

Therefore, in case\textbf{ }~B.2), ~there is an open ~$U\in\mathcal{N}_{_{K}}$
~with compact closure such that\smallskip{}

$\bullet$~~$\mbox{CMin}(U\setminus K)$ \emph{is a} $\mathfrak{c}-$\emph{dense
in itself, $d_{_{H}}-$open subset of }$\mbox{CMin}(M)$\emph{, $d_{_{H}}-$
accumulating in ~$\mathfrak{S}(K)$. }

\selectlanguage{english}%
\smallskip{}

\selectlanguage{british}%
Note that in particular, ~$\mbox{Per}(U\setminus K)\,\sqcup\,\mbox{Am}(U\setminus K)$
~is a ~$\mathfrak{c}-$\emph{dense in itself, }~$d_{_{H}}-$ open
~subset of ~$\mbox{CMin}(M)$~ since~ $\mbox{Eq}(M)$ ~is a ~$d_{_{H}}-$closed~
subset of ~$\mbox{CMin}(M)$. ~Now the above remark concerning ~$\mbox{CMin}(U\setminus K)$
implies there is a sequence ~$\varLambda_{_{n}}\in\mbox{CMin}(U\setminus K)$,
~$d_{_{H}}-$accumulating in ~$\mathfrak{S}(K).$ ~By the \emph{Cantor-Dirichlet
Principle} we may suppose this sequence is such that\emph{ all} ~$\varLambda_{_{n}}$~
\emph{belong to the same one of the following three classes: equilibrium
orbits }$\mbox{Eq}(M)$, \emph{periodic orbits }$\mbox{Per}(M),$
\emph{compact aperiodic minimals} $\mbox{Am}(M).$

Suppose now that conditions \textbf{~10.1}, \textbf{~10.2} ~and
\textbf{~10.3} ~all fail.~ We will show that condition \textbf{~10.4}
is necessarily true. The equality ~$\mbox{CMin}(N)=\mathcal{\mbox{Eq}}(N)\,\,\sqcup\,\,\mbox{Per}(N)\,\,\sqcup\,\,\mbox{Am}(N)$,
~valid for all ~$N\subset M$,~ will be repeatedly used $\big(\,\sqcup$
~denotes disjoint union$\big)$. ~Three possible cases are distinguished:\medskip{}

1st case:~~~There is a sequence~~$\varLambda_{_{n}}\in\mbox{Am}(U\setminus K)$~~
$d_{_{H}}$-accumulating in ~\foreignlanguage{english}{$\mathfrak{S}(K)$.}\medskip{}

~~~Since ~$\mbox{Am}(U\setminus K)$~~$d_{_{H}}-$accumulates
in~ $\mathfrak{S}(K)$~ but ~\textbf{10.3~ }is false, for every
open ~$V\in\mathcal{N}_{_{K}}$ ~there are ~$\varGamma\in\mbox{Am}(V\setminus K)$
~and ~$\epsilon>0$ ~such that\[
\#\big(B_{_{H}}(\varGamma,\,\epsilon)\,\cap\,\mbox{Am}(V\setminus K)\big)<\mathfrak{c}\]
Now since ~$V\setminus K$ ~is open and ~$\varGamma$ ~is closed,
~$B[\,\varGamma,\,\delta\,]\subset V\setminus K$ ~for a sufficiently
small ~$\delta>0$. ~Moreover by Lemma 5 $\big($$\varGamma\in\mbox{CMin}(M)$
~and ~$\mbox{Am}(M)\subset\mbox{Ci}(M)\big)$, taking ~$\delta$
~even smaller if necessary, we may further guarantee that \[
\mbox{Am}\big(B(\varGamma,\,\delta)\big)\subset B_{_{H}}(\varGamma,\,\epsilon)\,\cap\,\mbox{Am}(V\setminus K)\]
which implies\[
\#\,\mbox{Am}\big(B(\varGamma,\,\delta)\big)<\mathfrak{c}\]
Therefore it is easily seen that there are sequences ~$\varGamma_{_{n}}\in\mbox{Am}(U\setminus K)$
~and~~$\delta_{_{n}}>0$ ~such that:

\textbf{\hspace*{3mm}}1)\textbf{\hspace*{6mm}}$\big|\varGamma_{_{n}}\big|_{_{K}}\longrightarrow0$\smallskip{}

\textbf{\hspace*{3mm}}2)\textbf{\hspace*{6mm}}$B[\,\varGamma_{_{n}},\,\delta_{_{n}}\,]\subset U\setminus K$
\smallskip{}

\textbf{\hspace*{3mm}}3)\textbf{\hspace*{6mm}}$\mbox{\#\ensuremath{\,}Am}\big(B(\varGamma_{_{n}},\,\delta_{_{n}})\big)<\mathfrak{c}$\smallskip{}

Again, by~Lemma 11 (recall that ~$\varGamma_{_{n}}\subset U\setminus K$
~and ~$\mbox{cl\,}U$ ~is compact) we may replace condition 1)
above by\smallskip{}

\textbf{\hspace*{3mm}}1)\textbf{\hspace*{6mm}}$\varGamma_{_{n}}\overset{d_{_{H}}}{\longrightarrow}Q$
~for some ~$Q\in\mathfrak{S}(\mbox{bd\,}K)$

Clearly\smallskip{}

\textbf{\hspace*{3mm}}4)\textbf{\hspace*{6mm}}$\delta_{_{n}}<d_{_{H}}(\varGamma_{_{n}},\, Q)$\smallskip{}

since ~$B[\,\varGamma_{_{n}},\,\delta_{_{n}}]\,\cap\, Q=\emptyset$
~(recall that ~$\emptyset\neq Q\subset K\,$). 

Taking a smaller ~$\delta_{_{n}}$ ~if necessary, we may further
require that:

\smallskip{}

\textbf{\hspace*{3mm}}5)\textbf{\hspace*{6mm}}$\mbox{Eq}\big(B(\varGamma_{_{n}},\delta_{_{n}})\big)=\emptyset$
~~~~%
\footnote{the set ~$\mbox{Eq}(M)\subset\mbox{C}(M)$~~is $d_{_{H}}-$closed
and ~$\varGamma_{n}\notin\mbox{Eq}(M)$, ~(see Lemma 4).%
}\smallskip{}

\textbf{\hspace*{3mm}}6)\textbf{\hspace*{6mm}}$\gamma\in\mbox{Per}\big(B(\varGamma_{_{n}},\,\delta_{_{n}})\big)\implies\mbox{period}(\gamma)>n$~~~~~%
\footnote{By Lemma 10, this condition is necessarily satisfied for ~$\delta{}_{n}$
~small enough since ~$\varGamma_{_{n}}\in\mbox{Am}(M).$%
}

\smallskip{}

$\mbox{CMin}(U\setminus K)$ ~is ~nonvoid and ~$\mathfrak{c}\mbox{-}$\foreignlanguage{english}{\emph{dense
in itself}}, ~thus so is ~$\mbox{CMin}\big(B(\varGamma_{_{n}},\,\delta_{_{n}})\big)$
~since ~$B(\varGamma_{_{n}},\,\delta_{_{n}})\subset U\setminus K$
~is open and ~$\varGamma_{_{n}}\in\mbox{CMin}$$(M)$~; ~also
by 5),\[
\mbox{CMin}\big(B(\varGamma_{_{n}},\delta_{_{n}})\big)=\mbox{Per}\big(B(\varGamma_{_{n}},\delta_{_{n}})\big)\,\sqcup\,\mbox{Am}\big(B(\varGamma_{_{n}},\delta_{_{n}})\big)\]
hence in virtue of ~3),~ the \emph{Cantor-Dirichlet} \emph{Principle}
implies \smallskip{}

\textbf{\hspace*{3mm}}7)\textbf{\hspace*{6mm}}$\big(\,\,$$X\in\mbox{CMin}\big(\varGamma_{_{n}},\,\delta_{_{n}}\big)\,\,\,\,\,\,\,\,\mbox{and}\,\,\,\,\,\,\,\,$$\epsilon>0\,\,\big)\,\,\,\,\,\implies\,\,\,\,\,$$\#\,\mbox{Per}\big(B(X,\,\epsilon)\big)=\mathfrak{c}$\smallskip{}

In particular, ~~$P_{_{n}}:=\mbox{Per}\big(B(\varGamma_{_{n}},\,\delta_{_{n}})\big)$
~is an $d_{_{H}}-$open, ~$\mathfrak{c}\mbox{-}$\emph{dense in
itself} ~subset of ~$\mbox{Per}(M)$~ $d_{_{H}}-$accumulating
in ~$\varGamma_{_{n}}$~$\big($by 7 and Lemma 4$\,\big)$.~ Let
~$P:=\underset{n\in\mathbb{N}}{\bigcup}P_{_{n}}$. ~Then since ~$P$
~$d_{_{H}}-$accumulates in ~$\varGamma_{_{n}}$ ~and ~$\varGamma_{_{n}}\overset{d_{_{H}}}{\longrightarrow}Q\in\mathfrak{S}(\mbox{bd\,}K)$,
~it follows that ~\smallskip{}

$\bullet$~~~$P\subset\mbox{Per}(M\setminus K)$ ~\emph{is a}
~$\mathfrak{c}\mbox{-}$\emph{dense in itself,} \emph{~$d_{_{H}}-$open
subset of} ~$\mbox{Per}(M)$, ~\emph{$d_{_{H}}-$accumulating in}
~$\mathfrak{S}(\mbox{bd\,}K)$.

\smallskip{}

Moreover,\smallskip{}

$\bullet$~~~$K$~\emph{~is} ~\emph{bi-stable} \emph{in relation
to} ~$P^{*}=\underset{\gamma\in P}{\bigcup}\gamma\,\,\,\,$:\smallskip{}

$P^{*}$ ~is a union of periodic orbits and hence invariant. ~Given
any $V\in\mathcal{N}_{_{K}}$ ~let ~$\lambda>0$~~be such that
~~$B(K,\,\lambda)\subset V.$~~Since ~$Q\subset K,$ ~$d_{_{H}}(\varGamma_{_{n}},\, Q)\longrightarrow0$
~and ~$\delta_{_{n}}<d_{_{H}}(\varGamma_{_{n}},\, Q)$, ~there
is a ~$n_{_{0}}\geq1$ ~such that:\[
n>n_{_{0}}\implies d_{_{H}}(\varGamma_{_{n}},\, Q)<\lambda/2\implies\]
\[
\Longrightarrow\varGamma_{_{n}}\subset B(Q,\,\lambda/2)\subset B(K,\,\lambda/2)\mbox{\,\,\ and \,\,\ensuremath{\delta_{_{n}}}<\ensuremath{\lambda}/2\ensuremath{\implies}}\]
\[
\implies B(\varGamma_{_{n}},\,\delta_{_{n}})\subset B(K,\,\lambda)\implies\]
\[
\Longrightarrow P_{_{n}}:=\mbox{Per}\big(B(\varGamma_{_{n}},\,\delta_{_{n}})\big)\subset\mbox{Per}\big(B(K,\,\lambda)\big)\subset\mbox{Per}(V)\]

Since ~\emph{K }is closed, by 2) there is a ~$0<\delta<\lambda/2$
~such that\begin{equation}
B(K,\,\delta)\,\bigcap\,\left(\underset{1\leq n\leq n_{_{0}}}{\bigcup}B[\,\varGamma_{_{n}},\,\delta_{_{n}}]\right)=\varnothing\end{equation}
Therefore,\[
x\in B(K,\,\delta)\,\cap\, P^{*}\implies\big(x\in P_{_{n}}^{*}\mbox{\,\ for some \,}n>n_{_{0}}\big)\implies\mathcal{O}(x)\subset V\]
because ~$\mathcal{O}(x)\in P_{_{n}}$ ~and ~$P_{_{n}}\subset\mbox{Per}(V).$
~The \emph{bi-stability} of ~\emph{K} ~in relation to ~\emph{$P^{*}$~}is
proved.\smallskip{}

$\bullet$~~~\emph{for any sequence} ~$\gamma_{_{n}}\in P$, ~~$\mbox{dist}(\gamma_{_{n}},K)\longrightarrow0\implies\mbox{period}(\gamma_{_{n}})\longrightarrow+\infty\,\,\,$:\smallskip{}

By ~2), ~each ~$B[\,\varGamma_{_{n}},\,\delta_{_{n}}]$~is a compact
disjoint from ~\emph{K, ~}hence given any ~$n_{_{0}}\geq1$,~
there is a ~$\delta>0$~ satisfying the identity (8) above, therefore
by 6),\[
\big(\gamma\in P\mbox{\,\,\ and\,\,\,}\gamma\,\cap\, B(K,\,\delta)\neq\emptyset\big)\implies\big(\gamma\in P_{_{n}}\mbox{ \,\ for some\,\ }n>n_{_{0}}\big)\implies\]
\[
\implies\,\,\mbox{period}(\gamma)>n>n_{_{0}}\]

2nd case:~~~There is a sequence~~$\varLambda_{_{n}}\in\mbox{Per}(U\setminus K)$~~
$d_{_{H}}-$accumulating in ~\foreignlanguage{english}{$\mathfrak{S}(\mbox{bd\,}K)$\@.}\smallskip{}

Since by hypothesis condition \textbf{10.2} is not true, an argument
completely similar to that used in the 1st\emph{ }case proves that\smallskip{}

$\bullet$~~\emph{~there is a} $\mathfrak{c}-$\emph{dense in itself}
\emph{set} ~$A\subset\mbox{Am}(M\setminus K),$ ~\emph{$d_{_{H}}$-open
in} ~$\mbox{Am}(M),$ ~\emph{$d_{_{H}}$-accumulating in} ~$\mathfrak{S}(\mbox{bd\,}K)$
~\emph{and such that} ~$K$ \emph{~is} \emph{bi-stable} \emph{in
relation to} ~$A^{*}.$ \medskip{}

3rd case:~~~There is a sequence~~$\varLambda_{_{n}}\in\mbox{Eq}(U\setminus K)$
~~$d_{_{H}}$-accumulating in ~\foreignlanguage{english}{$\mbox{Eq}(\mbox{bd\,}K)$}.\smallskip{}

Since by hypothesis ~\textbf{10.1} ~is not true, there are necessarily
sequences ~$z_{_{n}}\in M$~and ~$\epsilon_{_{n}}>0$ ~such that:\[
\{z_{_{n}}\}\in\mbox{Eq}(U\setminus K)\,\,\,\,\,\,\,\,\,\,\,\,\mbox{and}\,\,\,\,\,\,\,\,\,\,\,\mbox{dist}(z_{_{n}},\, K)\longrightarrow0\,\,\,\,\,\,\,\,\,\,\,\mbox{and}\,\,\,\,\,\,\,\,\,\,\,\,\mbox{\#\ensuremath{\,}Eq}\big(B(z_{_{n}},\,\epsilon_{_{n}})\big)<\mathfrak{c}\]
By the $d_{_{H}}-$closeness of ~$\mbox{Eq}(M)$ ~in conjunction
with Lemma 11,\emph{ }we may suppose that ~$\{z_{_{n}}\}\overset{d_{_{H}}}{\longrightarrow}\{z\}$
~for some ~$\{z\}\in\mbox{Eq}(K)$. ~Clearly ~$z\in\mbox{bd\,}K$,
~thus ~$\{z\}\in\mbox{Eq}(\mbox{bd\,}K)$. ~Now\emph{ }~\emph{$\mbox{CMin}(U\setminus K)$
~is~ $\mathfrak{c}-$dense in itself,~} $\{z_{_{n}}\}\in\mbox{CMin}(U\setminus K)$~
for all\emph{ ~$n\geq1,$ }~$\mbox{CMin}(U\setminus K)\,=$$\mbox{Eq}(U\setminus K)\,\,\sqcup\,\,$$\mbox{Per}(U\setminus K)\,\,\sqcup\,\,$$\mbox{Am}(U\setminus K)$
~and moreover $\mbox{\#\ensuremath{\,}Eq}\big(B(z_{_{n}},\epsilon_{_{n}})\big)<\mathfrak{c}$,
~therefore we infer using\emph{ Cantor-Dirichlet Principle }that
there is a subsequence ~$\big(z_{_{n_{i}}}\big)$ ~such that:\[
\mbox{\#\ensuremath{\,}Am}\left(B\big(z_{_{n_{i}}},\epsilon\big)\right)=\mathfrak{c}\,\,\,\,\forall\, i\geq1,\,\,\,\epsilon>0\,\,\,\,\mbox{\ or\,\,\,\,\,}\mbox{\#\ensuremath{\,}Per}\left(B\big(z_{_{n_{i}}},\epsilon\big)\right)=\mathfrak{c}\,\,\,\,\,\mbox{ \ensuremath{\forall}\,}i\geq1\,,\,\,\,\epsilon>0\,\]

Thus by Lemma 4,\[
\big\{ z_{_{n_{i}}}\big\}\in\mbox{cl\ensuremath{_{_{H}}}Am}(U\setminus K)\,\,\,\,\,\,\,\forall\, i\geq1\mbox{\,\,\,\,\,\,\,\,\ or\,\,\,\,\,\,\,}\big\{ z_{_{n_{i}}}\big\}\in\mbox{cl\ensuremath{_{_{H}}}Per}(U\setminus K)\,\,\,\,\,\,\,\forall\, i\geq1\]
and since ~$\big\{ z_{_{n_{i}}}\big\}\overset{d_{_{H}}}{\longrightarrow}\big\{ z\big\}\in\mbox{Eq}(\mbox{bd\,}K)\subset\mathfrak{S}(\mbox{bd\,}K)$,
~it follows that~ \emph{there is a sequence in }~$\mbox{Am}(U\setminus K)$
~$d_{_{H}}$-\emph{accumulating in} ~$\mathfrak{S}(\mbox{bd\,}K)$
\emph{~or~~there is a sequence in} ~$\mbox{Per}(U\setminus K)$
~$d_{_{H}}$\emph{-accumulating in} ~$\mathfrak{S}(\mbox{bd\,}K)$.
~Thus, the 3rd case implies the ~1st~ or the ~2nd$.$~On the
other hand, as we have seen, the 1st case clearly implies the ~2nd~
and vice versa, hence the ~1st~ and ~2nd~ cases always occur.~
Therefore if conditions \textbf{~10.1}, \textbf{~10.2} ~and~ \textbf{10.3~}
are false then \textbf{~10.4~} is true. The proof of ~Theorem 1~
is complete.\hfill{}\emph{Q.E.D.}

\medskip{}

\noun{7. INDEPENDENT REALIZATIONS.~ EXAMPLES.}\medskip{}

~~~With Theorem 1 established, the question of whether all the
28 cases it describes are \emph{realizable} naturally arises. ~Furthermore
we may doubt whether all these cases are \emph{mutually} \emph{independent}.
~Let ~$(M,\theta)$ ~be a ~$C^{\,0}$ ~flow on a compact, connected
metric space ~and ~$\varSigma$ ~one of the twenty eight conditions
of Theorem 1. ~We say that ~$(M,\theta)$ ~is an \emph{independent
realization }of ~$\varSigma$ ~ if there ~a compact, invariant,
proper subset~ $K\subset M$ ~such that condition ~$\varSigma$
~is satisfied for this choice of ~$M$, ~$K$ ~and ~$\theta$
~but none of the remaining conditions of Theorem 1 holds, for the
same ~$M$, ~$K$ ~and ~$\theta$.%
\footnote{Sometimes we actually identify a specific compact invariant proper
subset ~$K$, ~saying that ~$(M,\theta)$, ~$K$ ~constitutes
an independent realization of the condition in question.%
}

We can give the following answer to the questions raised above: ~all
conditions, except the twelve involving ~$X\mbox{-}trees$, ~$X\mbox{-}\alpha shells$
~and ~$X\mbox{-}\omega shells$ ~(that is, all but the {}``exceptional''
conditions\textbf{ ~4}, \textbf{~5}, \textbf{~6}, \textbf{~7.4~}
to ~\textbf{7.6}, \textbf{~8.4~} to ~\textbf{8.6} ~and \textbf{~9.4~}
to ~\textbf{9.6}), admit independent realizations\emph{ }by ~$C^{\,\infty}$
~flows on compact, connected manifolds $\big($in fact on ~$\mathbb{S}^{n}$
~for some ~$1\leq n\leq4\,)$. ~Each of the exceptional conditions
admits a ~$\big(C^{\,0}\big)$ \emph{~}independent realization on
which\emph{ ~$M$ }~is\emph{ }a compact, connected invariant subset
of a ~$C^{\,\infty}$ ~flow ~$\phi$ ~on ~$\mathbb{R}^{n}$,
~$\theta$ ~is ~$\phi|_{\mathbb{R}\times M}$, ~$K\subsetneq M$
~is an equilibrium orbit of ~$\phi$ ~and ~$3\leq n\leq6$. ~A
moment of reflexion shows that this implies that in the last sentence
we may substitute ~$\mathbb{R}^{n}$ ~by any (2nd countable, Hausdorff)
~$C^{\,\infty}$ ~manifold ~$\mathcal{M}$ ~of dimension ~$m\geq n$.
~We call such an independent realization\emph{ }of one of the twelve
exceptional conditions of Theorem 1 a \emph{subsmooth independent
realization} and call the manifold ~$\mathcal{M}\supset M$ ~carrying
the ~$C^{\,\infty}$ ~flow ~$\phi$ ~the \emph{ambient manifold
}of this realization\emph{.} ~Whether each of these twelve exceptional
conditions admits an independent realization by a ~$C^{\, r}$ ~$(r\geq1)$
~flow on a compact, connected ~$C^{\,\infty}$ ~manifold requires
further exam. ~In \cite{teix} we give an independent realization\emph{
}for each one of the 28 conditions of Theorem 1, in the line of what
was stated above. ~Here, for the sake of brevity, we will confine
our attention to the twelve exceptional conditions,%
\footnote{\emph{i.e. }to conditions\textbf{ ~4}, \textbf{~5}, \textbf{~6},
\textbf{~7.4~} to ~\textbf{7.6}, \textbf{~8.4~} to ~\textbf{8.6},
\textbf{~9.4~} to ~\textbf{9.6}.%
} since these present greater difficulties than the others and display
less known and more interesting dynamical phenomena. ~Actually the
main difficulty lies in obtaining subsmooth independent realization\emph{
}of conditions ~\textbf{4}, ~\textbf{5} ~and ~\textbf{6} ~with
~$K$ ~an equilibrium orbit, a periodic orbit and a compact aperiodic
minimal. ~Once this done, we have the essential {}``dynamical pieces''
for the construction of subsmooth independent realizations of the
remaining nine exceptional conditions ~\textbf{7.4}~- \textbf{7.6},
~\textbf{8.4} - \textbf{8.6 }and ~\textbf{9.4} - \textbf{9.6}, ~which
are then easily achieved through well known vector field constructions
on \emph{trivial} normal tubular neighbourhoods.

\begin{figure}[H]
\noindent \begin{centering}
\includegraphics[scale=0.9]{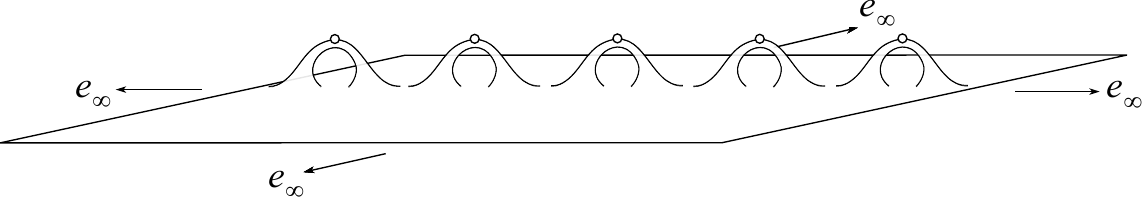}
\par\end{centering}

\caption{{\small The infinite genus smooth surface ~$S\subset\mathbb{R}^{3}.$}}
\end{figure}

We shall first produce examples of subsmooth independent realizations
for condition ~\textbf{6} ~with ~$K$, ~respectively, an \emph{equilibrium
orbit}, a \emph{periodic orbit }and a \emph{compact aperiodic minimal
}(a ~$\mathbb{T}^{2}$ ~with an irrational\emph{ }linear flow).%
\footnote{\emph{Remark. }~Vector fields on submanifolds ~$M\subset\mathbb{R}^{n}$
~will be always represented in the usual abridged form ~$X:M\longrightarrow\mathbb{R}^{n}$
~\emph{i.e.} instead of considering ~$\upsilon:M\longrightarrow TM\subset T\mathbb{R}^{n}=\mathbb{R}^{n}\times\mathbb{R}^{n}$,
~we work with ~$X=\pi_{_{2}}\circ\upsilon$, ~where ~$\pi_{_{2}}$
~is the projection onto the 2nd factor.~ Where doubts may arise
concerning what flow is aimed at, we indicate it inside brackets or
as a subscript \emph{e.g. }we write ~$\mbox{Min}(X,\theta)$ ~and~
$\omega_{_{\theta}}(z)$ ~instead of ~$\mbox{Min}(X)$ ~and ~$\omega(z)$.
~When giving examples of flows generated by ~$C^{\, r}$ ~$(r\geq1)$
~vector fields on manifolds ~$M$, often the vector field in question
is indicated as a subscript \emph{e.g. }we write ~$\omega_{_{X}}(z)$
~for the ~$\omega$-limit set of ~$z$ ~on the flow ~$X^{t}$
~generated by ~$X\in\mathfrak{X}^{\, r}(M)$.\emph{ ~}We use these
notations freely (with the subscript indicating either the flow or
the generating vector field) since no risk of ambiguity arises.%
}\medskip{}

Example 1\emph{.} ~\emph{Subsmooth independent realization of condition}
~\textbf{6 }~with ~$M$ ~an orbit closure of a ~$C^{\,\infty}$
~flow ~$\zeta^{t}$ ~on ~$\mathbb{S}^{3}\subset\mathbb{R}^{4}$,
~$\theta=\zeta^{t}|_{\mathbb{R}\times M}$ ~and ~$K=\big\{(0,0,0,1)\big\}\subsetneq M$
~an \emph{equilibrium orbit}.

Our point of departure is a beautiful example, due to Beniere and
Meigniez \cite{beni}, of a ~$C^{\,\infty}$ ~complete vector field
~$\upsilon$ ~generating a \emph{flow without minimal sets }on a
non-compact, orientable surface ~$\mathfrak{M}$ ~of infinite genus.%
\footnote{The first example of a ~$C^{\,\infty}$ ~flow on a manifold, without
minimal sets was, to our knowledge, constructed on a surface by Takashi
Inaba in 1995 (see \cite{inab}). ~This achievement is, in our opinion,
one of the highlights of exceptional dynamics in the second half of
the last century. %
} ~The set of \emph{end points} ~$E(\mathfrak{M})$ ~is homeomorphic
to ~$\varDelta:=$$\{0\}\,\cup\,\{1/n:\, n\in\mathbb{N}\}\subset\mathbb{R}$
~and all end points are \emph{flat},%
\footnote{an end ~$e\in E(\mathfrak{M})$ ~is \emph{flat} if it has a neighbourhood
homeomorphic to ~$\mathbb{R}^{2}$ ~in the \emph{end-points compactification}
~$\mathfrak{M}^{\propto}=\mathfrak{M}\,\sqcup\, E(\mathfrak{M})$
~of ~$\mathfrak{M}$. ~Richards \cite{rich} calls such an end
point \emph{{}``planar}''. ~Beniere and Meigniez \cite{beni} designate
by ~$M$ ~our surface ~$\mathfrak{M}$.%
} except the non-isolated one. ~We shall first construct a smoothly
$(C^{\,\infty})$ ~embedded surface ~$S\subset\mathbb{R}^{3}$ ~that
is ~$C^{\,\infty}$ ~diffeomorphic to ~$\mathfrak{M}$. ~Operate
the following ~$C^{\,\infty}$ ~surgery within the ambient manifold
~$\mathbb{R}^{3}$: ~to the plane ~$\mathbb{R}^{2}\times\{0\}\subset\mathbb{R}^{3}$
~smoothly add denumerably many handles%
\footnote{the resulting surface is known as the \emph{infinite Loch-Ness monster.}%
} as shown in fig. 16. ~From each handle remove one point ~$e_{_{n}},$
~$n\in\mathbb{N}.$ ~We obtain a non-closed, smooth 2-submanifold
~$S\subset\mathbb{R}^{3}$. ~As ~$\mathfrak{M}$, ~$S$ ~is an
orientable $C^{\,\infty}$ surface of infinite genus with all ends
isolated and flat except one, ~$e_{_{\infty}}$,~ which is both
non-isolated and non-flat. ~This implies that its end points set\emph{
~$E(S)$} ~is also homeomorphic to ~$\varDelta$ ~(see above),
hence there is a homeomorphism ~$\xi:E(\mathfrak{M})\longrightarrow E(S)$
~sending the unique non-flat\emph{ }end of ~$\mathfrak{M}$ ~to
the unique non-flat\emph{ }end of ~$S$. ~By \emph{Kerekjarto Theorem
}(see \emph{e.g.} \cite{rich}, p.262 and \cite{beni}, p.26) the
surfaces ~$\mathfrak{M}$ ~and ~$S$ ~are homeomorphic,%
\footnote{we need not to care with non-orientable ends since there are none:
both ~$\mathfrak{M}$ ~and ~$S$ ~are orientable.%
} hence, as it is well known, ~$C^{\,\infty}$ ~diffeomorphic. ~Let
~$f:\mathfrak{M}\longrightarrow S$ ~be a ~$C^{\,\infty}$ ~diffeomorphism
defining a ~$C^{\,\infty}$ ~embedding ~$\mathfrak{M}\hookrightarrow\mathbb{R}^{3}\supset S$
~and inducing the \emph{$C^{\,\infty}$ ~complete }tangent vector
field\emph{ }~$X:=f_{_{*}}\upsilon$ ~on ~$S\subset\mathbb{R}^{3}$.
~As occurs with ~$\upsilon$, ~the flow ~$X^{t}$ ~\emph{has
no minimal sets} ~($f$ ~realizes a ~$C^{\,\infty}$\emph{ ~}flow
conjugation).

\textbf{Definition.} ~Let ~$\theta$ ~be a ~$C^{\,0}$ ~flow
on a metric space~ $M$. ~A point ~$x\in M$ ~is called a \emph{limit
point }of ~$(M,\theta)$ ~if ~$x$ ~belongs to the ~$\alpha$-limit
set or to the $\omega$-limit set of some point of ~$M$. ~In this
case the orbit ~$\mathcal{O}(x)$ ~is called a \emph{limit orbit}
of the flow. ~We denote the set of \emph{limit points }of the flow
~$(M,\theta)$ ~by ~$\Upsilon_{\theta}$ ~and if the flow is given
by a vector field ~$\upsilon$, ~by ~$\Upsilon_{\upsilon}$.

From the inductively constructed tangentially orientable foliated
atlas of ~$\mathfrak{M}$ ~corresponding to the vector field ~$\upsilon$
~(given in \cite{beni}), it is easily seen that:

- ~each \emph{limit point} ~$x\in\mathfrak{M}$ ~has nonvoid ~$\alpha$-limit
and ~$\omega$-limit sets and both the positive and negative orbit
of ~$x$ ~accumulate in the unique non-isolated end of ~$\mathfrak{M}$
~and in no other end of this surface.

Since a homeomorphism between non-compact surfaces uniquely extends
to a homeomorphism between their respective end-points compactifications,\emph{
}it follows that in the flow ~$X^{t}$~ generated by ~$X\in\mathfrak{X}^{\,\infty}(S)$,
~both the ~$\alpha_{_{X}}$-limit ~and $\omega_{_{X}}$-limit sets
of each \emph{limit point }of ~$(S,X^{t})$ ~are\emph{ closed, unbounded
}%
\footnote{both these sets are closed subsets of ~$S$ ~that do not accumulate
in the isolated ends ~$e_{_{n}},$ ~$n\in\mathbb{N}$ ~of this
surface and the closure of ~$S$ ~in ~$\mathbb{R}^{3}$ ~equals
~$S\,\sqcup\,\{e_{_{n}}:\, n\in\mathbb{N}\}.$ ~Their\emph{ unboundeness
}also follows from the fact that a closed, bounded subset of ~$\mathbb{R}^{3}$
~is compact and thus if it is a nonvoid, invariant subset of the
flow ~$X^{t}$, ~then it must contain a minimal set\emph{ }of it\emph{.
~}But ~$X^{t}$ ~has no\emph{ }minimal sets.%
} subsets of ~$\mathbb{R}^{3}$. ~Let ~$U$ ~be a normal tubular
neighbourhood of ~$S$ ~in ~$\mathbb{R}^{3}$ ~(indeed a trivial
1-dimensional vector bundle over ~$S\,)$. ~Extend ~$X:=f_{_{*}}\upsilon\in\mathfrak{X}^{\,\infty}(S)$
~to a non-singular\emph{ }vector field ~$X_{_{0}}\in\mathfrak{X}^{\,\infty}(U)$
~defining,\[
\begin{array}{lll}
X_{_{0}}:U & \longrightarrow & \mathbb{R}^{3}\\
\,\,\,\,\,\,\,\,\,\,\,\,\,\,\, z & \longmapsto & X\circ\pi(z)\end{array}\]
where ~$\pi:U\longrightarrow S$~ is the canonical ~$C^{\,\infty}$
~submersion (orthogonal projection of ~$U$ ~over ~$S$).%
\footnote{\emph{i.e .}~$U=\underset{x\in S}{\bigsqcup}\mathfrak{D}_{_{\lambda(x)}}(x)$,
~where ~$\lambda\in C^{\,\infty}\big(S,\mathbb{R}^{+}\big)$, ~$\mathfrak{D}_{_{\lambda(x)}}(x)$
~is the 1-dimensional affine open disk with centre ~$x\in S$ ~and
radius ~$\lambda(x)$ ~orthogonal to ~$S$ ~at ~$x$ ~and ~$\pi:U\longrightarrow S$
~is the ~$C^{\,\infty}$ ~submersion defined by:\[
\pi\big(\mathfrak{D}_{_{\lambda(x)}}(x)\big)=\{x\}\,\,\,\,\,\,\,\mbox{for all \,}x\in S\]
}~ Let ~$p:=(0,0,0,1)$, ~$O:=(0,0,0,0)$, ~$\varphi:\mathbb{R}^{3}\longrightarrow\mathbb{S}^{3}\setminus\{p\}$
~the inverse stereographic projection ~$\big(\mathbb{R}^{3}$ ~identified
with ~$\mathbb{R}^{3}\times\{-1\}\,\big)$. ~$\varphi$ ~induces
the ~$C^{\,\infty}$ vector field ~$\varphi_{_{*}}X_{_{0}}$ ~on
the open subset ~$\varphi(U)$ ~of ~$\mathbb{S}^{3}$.

\begin{figure}[H]
\noindent \begin{centering}
\includegraphics[scale=0.3]{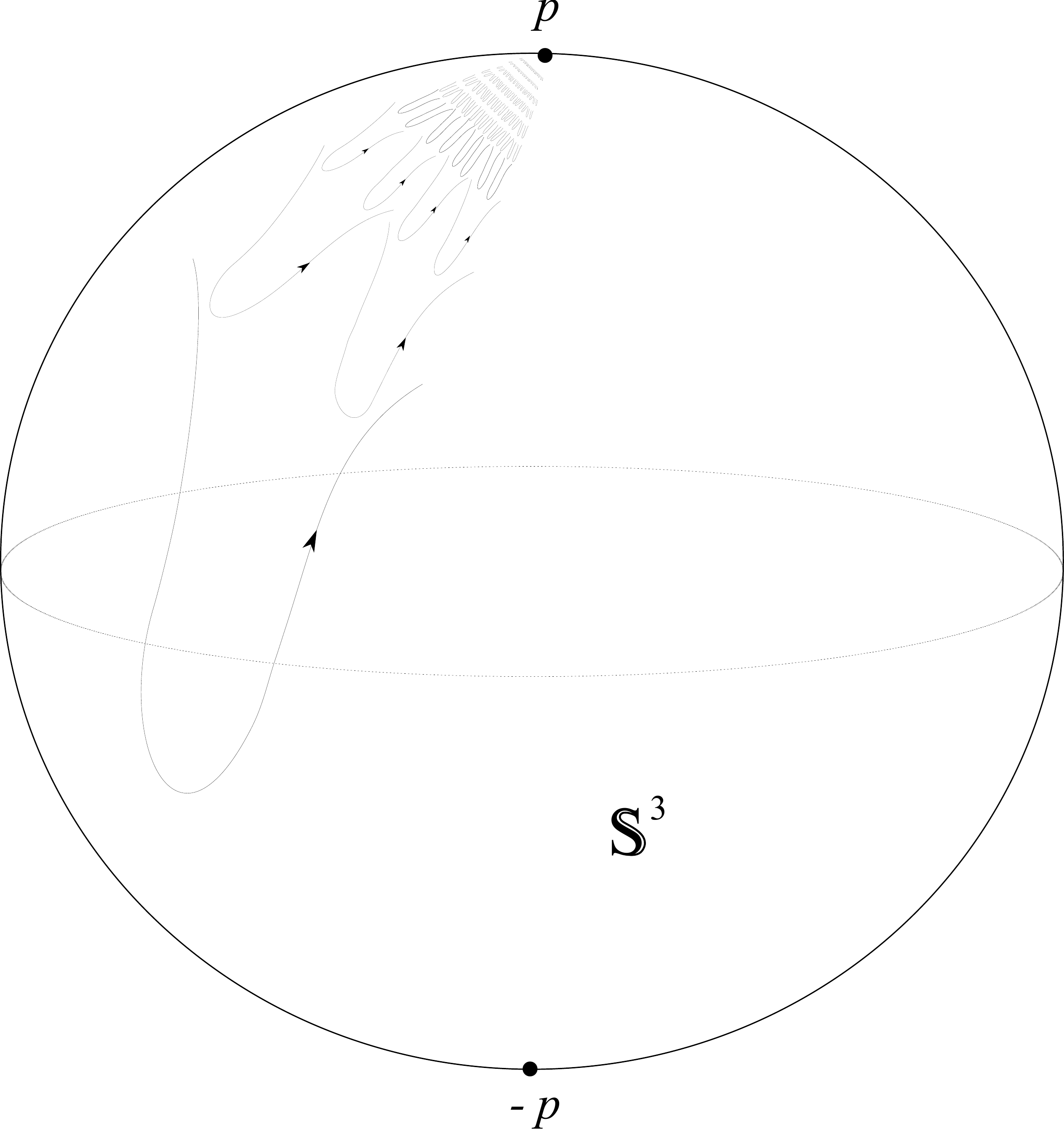}
\par\end{centering}

\caption{$K\mbox{-}tree$ ~in the flow ~$\big(\mathbb{S}^{3},\,\zeta^{t}\,\big)$.}
\end{figure}

By \emph{Kaplan Smoothing Theorem} (Kaplan \cite{kapl}, p.157), there
is a scalar function ~$\lambda\in C^{\,\infty}\big(\mathbb{S}^{3},[\,0,1\,]\big)$
~with ~$\lambda^{-1}(0)=\mathbb{S}^{3}\setminus\varphi(U)$ ~and
such that \[
\begin{array}{llll}
\zeta:\mathbb{S}^{3}\subset\mathbb{R}^{4} & \longrightarrow & \mathbb{R}^{4}\\
\,\,\,\,\,\,\,\,\,\,\,\,\, s & \longmapsto & O & \mbox{\,\,\,\ on\,\,\,\,}\mathbb{S}^{3}\setminus\varphi(U)\\
\,\,\,\,\,\,\,\,\,\,\,\,\, s & \longmapsto & \lambda\varphi_{_{*}}X_{_{0}}(s) & \mbox{\,\,\,\ on\,\,\,\,}\varphi(U)\end{array}\]
defines a ~$C^{\,\infty}$~ vector field on ~$\mathbb{S}^{3}$.
~The smoothly embedded surface $\varphi(S)$ ~is invariant under
the flow ~$\zeta^{t}$ ~$\big(\,\varphi|_{S}$ ~realizes a ~$C^{\,\infty}$
~flow conjugation between the global flow ~$\big(S,X^{t}\big)$
~and ~$\big(\varphi(S),(\varphi_{_{*}}X)^{t}\big)$; ~moreover
~$\mbox{im\,}\lambda|_{\varphi(S)}\subset]\,0,1\,],$ ~thus ~$\lambda\varphi_{_{*}}X=\zeta|_{\varphi(S)}$
~is necessarily a \emph{complete }vector field, topologically equivalent
to ~$X\in\mathfrak{X}^{\,\infty}(S)$ ~via the ~smooth diffeomorphism
~$\varphi\,\big)$.

~Let ~$q\in\Upsilon_{_{X}}$ ~and ~$z:=\varphi(q)$. ~Recall
that ~$\alpha_{_{X}}(q)$, ~$\omega_{_{X}}(q)$ ~and ~$\mbox{cl\,}\mathcal{O}_{_{X}}(q)$
~are unbounded, closed subsets of\emph{ }~$\mathbb{R}^{3}$, ~therefore
for all such ~$q$ ~and ~$z$,\begin{equation}
\begin{array}{c}
\alpha_{_{\zeta}}(z)=\varphi\big(\alpha_{_{X}}(q)\big)\,\sqcup\,\{p\}\mbox{\,\,\,\,\,\,\,\,\ and\,\,\,\,\,\,\,\,\,}\omega_{_{\zeta}}(z)=\varphi\big(\omega_{_{X}}(q)\big)\,\sqcup\,\{p\}\\
\,\\
\mbox{cl\,}\mathcal{O}_{_{\zeta}}(z)=\varphi\big(\mbox{cl\,}\mathcal{O}_{_{X}}(q)\big)\,\sqcup\,\{p\}\end{array}\end{equation}
Let ~$M:=\mbox{cl\,}\mathcal{O}_{_{\zeta}}(z)$, ~~$\theta:=\zeta^{t}|_{\mathbb{R}\times M}$
~and ~$K:=\{p\}.$ ~$\theta$ ~is a ~$C^{\,0}$ ~flow on the
compact, connected metric space ~$M\subset\mathbb{S}^{3}$ ~$\big($with
the euclidean metric of ~$\mathbb{R}^{4}\supset\mathbb{S}^{3}\,\big)$
~and ~$K$ ~is a compact, invariant proper subset of ~$M$. ~Now
with respect to the (sub)flow ~$(M,\theta)$, ~it is clear from
(9) that every ~$y\in M\setminus K=M\setminus\{p\}$ ~belongs to
~$A_{_{\theta}}^{-}(K)\,\cap\, A_{_{\theta}}^{+}(K)$ ~as ~$M=\varphi\big(\mbox{cl\,}\mathcal{O}_{_{X}}(q)\big)\,\sqcup\,\{p\}$
~and ~$\mbox{cl\,}\mathcal{O}_{_{X}}(q)=\mathcal{O}_{_{X}}(q)\,\cup\,\alpha_{_{X}}(q)\,\cup\,\omega_{_{X}}(q)\subset\Upsilon_{_{X}}$.
~Obviously ~$K$ ~is \emph{isolated from minimals} in ~$(M,\theta)$
~and no ~$x\in M\setminus K$ ~has its ~$\alpha_{_{\theta}}$-limit
set or its ~$\omega_{_{\theta}}$-limit set contained in ~$K$ ~\emph{$\big($i.e.
}equal to\emph{ $\{p\}\big)$. ~}This immediately implies that, for
this choice of ~$M$, ~$\theta$ ~and ~$K$, ~none of conditions
of Theorem 1, with the exception of condition \emph{~}\textbf{6},
~can hold.\textbf{~ }Therefore by Theorem 1 the above ~$M$, ~$\theta$~~and~
$K$ ~necessarily provide a \emph{subsmooth independent realization
}of condition ~\textbf{6}, ~with\emph{ ambient manifold} ~$\mathcal{M}=\mathbb{S}^{3}$.~
This can be easily verified directly: the existence of a ~$K\mbox{-}tree$,
~with ~$K=\{p\}$ ~and ~$\gamma_{_{0}}=\mathcal{O}_{_{\theta}}(z)$
~is now a straightforward inductive consequence of Lemma 3, since
every ~$y\in M\setminus K$ ~belongs to ~$A_{_{\theta}}^{-}(K)\,\cap\, A_{_{\theta}}^{+}(K)$
~and ~$K$ ~is consequently \emph{non-stagnant }in ~$(M,\theta)$. 

The example above provides the essential {}``dynamical piece'' for
the construction of a subsmooth independent realization\emph{ }of
condition ~\textbf{7.6}$.$ ~But in order to do the same for\emph{
}conditions\emph{ }\textbf{~8.6~ }and ~\textbf{9.6 ~}we need to
construct a subsmooth independent realization of condition ~\textbf{6}
~with ~$K$ ~a \emph{periodic orbit }and one with ~$K$ ~a \emph{compact
aperiodic minimal}, respectively\emph{.}

\medskip{}

Example 2. \emph{subsmooth independent realization }of condition ~\textbf{6
}~with ~$M$ ~an orbit closure of a ~$C^{\,\infty}$ ~flow ~$\nu^{t}$
~on the \emph{ambient} \emph{manifold} ~$\mathcal{M}=\mathbb{S}^{3}\times\mathbb{S}^{1}$,
~$\theta=\nu^{t}|_{\mathbb{R}\times M}$ ~and ~$K=\{p\}\times\mathbb{S}^{1}\subsetneq M$
~a\emph{ periodic orbit}.

Endow the compact manifold ~$\mathcal{M}=\mathbb{S}^{3}\times\mathbb{S}^{1}$
~with the ~$C^{\,\infty}$ ~vector field ~$\nu=\zeta+\frac{\partial}{\partial\varrho}$
~where ~$\varrho$ ~is the ~$\mathbb{S}^{1}$ ~coordinate. ~Call
the circles ~$\{z\}\times\mathbb{S}^{1}$, ~$z\in\mathbb{S}^{3}$
~the \emph{parallels }of ~$\mathcal{M}.$~ The parallel ~$\gamma=\{p\}\times\mathbb{S}^{1}$
~is a \emph{periodic orbit} of ~$\nu$ ~(with period ~$2\pi)$,
~since ~$p$ ~is a singularity of ~$\zeta\in\mathfrak{X}^{\,\infty}\big(\mathbb{S}^{3}\big)$.
~Observe that the flow ~$\big(\mathcal{M},\nu^{t}\big)$ ~projects
to the flow ~$\big(\mathbb{S}^{3},\zeta^{t}\big)$ ~via the projection
onto the 1st factor ~$\pi_{_{1}},$ ~more precisely, for every ~$t\in\mathbb{R}$
~and $u\in\mathbb{S}^{3}$, ~$s\in\mathbb{S}^{1},$\begin{equation}
\pi_{_{1}}\circ\nu^{t}(u,s)=\zeta^{t}\big(\pi_{_{1}}(u,s)\big)=\zeta^{t}(u)\end{equation}

\begin{figure}[H]
\noindent \begin{centering}
\includegraphics[scale=0.7]{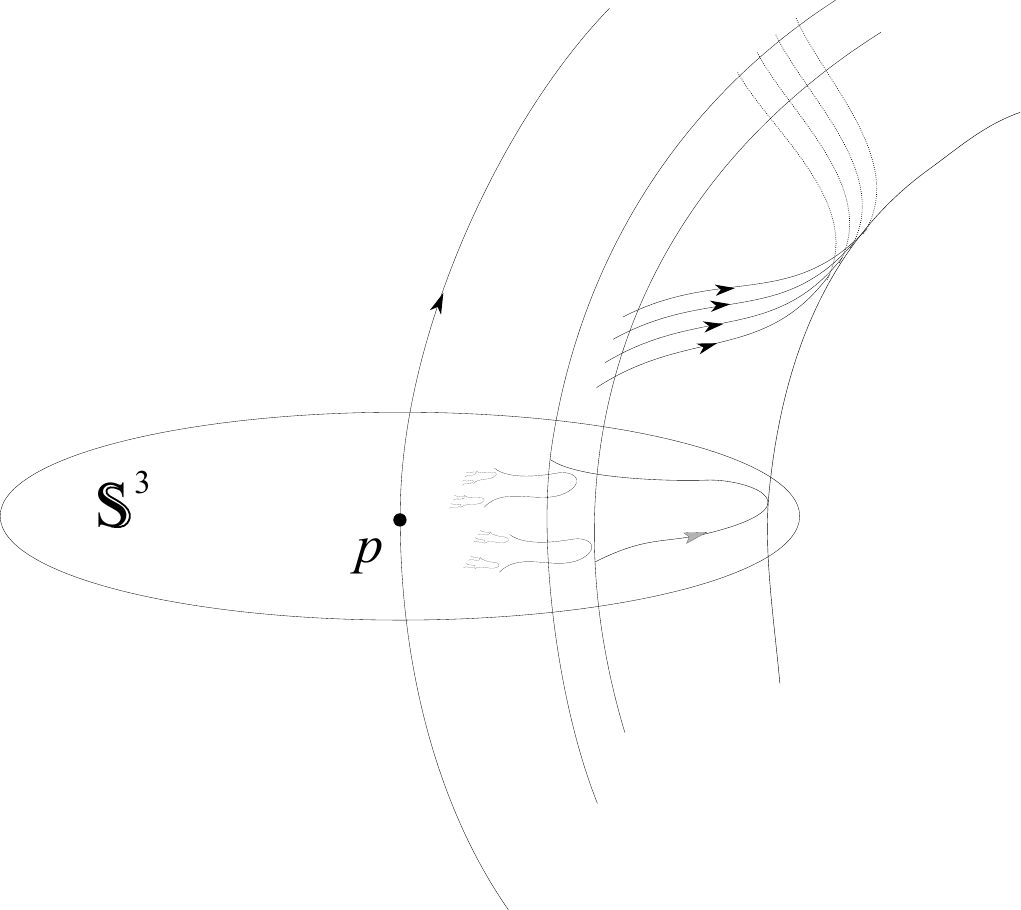}
\par\end{centering}

\caption{{\small Example 2, showing the projection ({}``shadow'') of the
flow $\big(\mathbb{S}^{3}\times\mathbb{S}^{1},v^{t}\big)$ to the
flow $\big(\mathbb{S}^{3},\zeta^{t}\big)$$.$}}
\end{figure}

and similarly ~$\big(\mathcal{M},\nu^{t}\big)$ ~projects to the
flow ~$\big(\mathbb{S}^{1},\left(\partial/\partial\varrho\right){}^{t}\big)$
~via the projection onto the 2nd factor ~$\pi_{_{2}},$ ~and thus
identifying ~$\mathbb{S}^{1}$ ~with ~$\mathbb{R}/2\pi\mathbb{Z}=\mathbb{R}\,(\mbox{mod\,2\ensuremath{\pi}})$,\begin{equation}
\nu^{t}(u,s)=\big(\zeta^{t}(u),\, s+t\,(\mbox{mod\,2\ensuremath{\pi}})\big)\end{equation}

We will deduce all the facts needed about the dynamics of ~$\nu^{t}$
~from these equalities.~ Observe that for every ~$u\in\mathbb{S}^{3}$,
~$s\in\mathbb{S}^{1},$

\hspace{3mm}a)\hspace{3mm}$\pi_{_{1}}\big(\mathcal{O}_{_{\nu}}(u,s)\big)=\mathcal{O}_{_{\zeta}}(u)$

\hspace{3mm}b)\hspace{3mm}$\pi_{_{1}}\big(\omega_{_{\nu}}(u,s)\big)=\omega_{_{\zeta}}(u)$

\hspace{3mm}c)\hspace{3mm}$\pi_{_{1}}\big(\alpha_{_{\nu}}(u,s)\big)=\alpha_{_{\zeta}}(u)$

\hspace{3mm}d)\hspace{3mm}$\pi_{_{1}}\big(\mbox{cl\,}\mathcal{O}_{_{\nu}}(u,s)\big)=\mbox{cl\,}\mathcal{O}_{_{\zeta}}(u)$

a) is an immediate consequence of (10); ~b) and c) are time symmetric
facts. ~To prove b) observe that the inclusion ~$\subset$~ follows
immediately from (10) and the reciprocal ~$\supset$ ~follows from
the compactness of ~$\mathbb{S}^{1}$: ~if ~$q\in\omega_{_{\zeta}}(u)$,
~then for any ~$s\in\mathbb{S}^{1}=\mathbb{R}\,(\mbox{mod\,}2\pi)$,
~$\omega_{_{\nu}}(u,s)$ ~necessarily intercepts the \emph{parallel
~$\{q\}\times\mathbb{S}^{1}$}. ~To see this, let ~$t_{_{n}}\longrightarrow+\infty$
~be such that ~$\zeta^{t_{_{n}}}(u)\longrightarrow q\in\omega_{_{\zeta}}(u)$.
~The sequence ~$s+t_{_{n}}(\mbox{mod\,}2\pi)$ ~has a subsequence
~$\big(s+t_{_{n_{i}}}(\mbox{mod\,2\ensuremath{\pi}})\big)_{i\in\mathbb{N}}$
~converging to some ~$s_{_{0}}$ ~belonging to the compact ~$\mathbb{S}^{1}=\mathbb{R}\,(\mbox{mod\,2\ensuremath{\pi}})$
~and thus by (11),\[
\underset{i\rightarrow+\infty}{\mbox{lim}}\nu^{t_{_{n_{i}}}}(u,s)=\underset{i\rightarrow+\infty}{\mbox{lim}}\big(\zeta^{t_{_{n_{i}}}}(z),\, s+t_{_{n_{i}}}\,(\mbox{mod\,2\ensuremath{\pi}})\big)=(q,s_{_{0}})\in\mathbb{S}^{3}\times\mathbb{S}^{1}\]
 therefore ~$(q,s_{_{0}})\in\omega_{_{\nu}}(u,s)$ ~since ~$t_{_{n_{i}}}\longrightarrow+\infty$.
~ Finally, d) is an immediate consequence of a), b) and c).

Now let ~$z\in\varphi(\Upsilon_{_{X}})\in\mathbb{S}^{3}$ ~(see
example 1) and ~$N:=\mbox{cl\,}\mathcal{O}_{_{\nu}}(z,b)\subset\big(\mbox{cl\,}_{_{\zeta}}(z)\big)\times\mathbb{S}^{1},$
~$\phi:=\nu^{t}|_{\mathbb{R}\times N}$ ~and ~$J:=\gamma=\{p\}\times\mathbb{S}^{1}$
(a periodic orbit, see above), where ~$b$~ is an arbitrary point
of ~$\mathbb{S}^{1}$. ~Observe that from the four identities above
it follows that\[
(u,s)\in\mbox{cl\,}\mathcal{O}_{_{\nu}}(z,b)\setminus\gamma\implies u\in\mbox{cl\,}\mathcal{O}_{_{\zeta}}(z)\setminus\{p\}\implies\]
\[
\implies u\in A_{_{\zeta}}^{-}\big(\{p\}\big)\,\cap\, A_{_{\zeta}}^{+}\big(\{p\}\big)\implies(u,s)\in A_{_{\nu}}^{-}(J)\,\cap\, A_{_{\nu}}^{+}(J)\implies\]
\[
\implies(u,s)\in A_{_{\phi}}^{-}(J)\,\cap\, A_{_{\phi}}^{+}(J)\]

Therefore ~$J$ ~is \emph{isolated from minimals} in ~$(N,\phi)$
~and no ~$x\in N\setminus J$ ~has its ~$\alpha_{_{\phi}}$-limit
set or its ~$\omega_{_{\phi}}$-limit set contained in ~$J$ ~\emph{$\big($i.e.
}equal to\emph{ }the periodic orbit\emph{ ~$\gamma=\{p\}\times\mathbb{S}^{1}$$\big)$.
~}This immediately implies that, for this choice of ~$M:=N$, ~$\theta:=\phi$
~and ~$K:=J$, ~none of conditions of Theorem 1, with the exception
of condition \emph{~}\textbf{6}, ~can hold.\textbf{~ }Therefore
by Theorem 1 the above ~$N$, ~$\phi$~~and~ $J$ ~necessarily
provide a subsmooth independent realization\emph{ }of condition ~\textbf{6
~}with\textbf{ ~$J$ ~}a \emph{periodic orbit, }the\emph{ }ambient
manifold being ~$\mathcal{M}=\mathbb{S}^{3}\times\mathbb{S}^{1}.$

\medskip{}

Example 3. \emph{subsmooth independent realization }of condition ~\textbf{6
}~with ~$M$ ~an orbit closure of a ~$C^{\,\infty}$ ~flow ~$\upsilon^{t}$
~on the \emph{ambient} \emph{manifold} ~$\mathcal{N}=\mathbb{S}^{3}\times\mathbb{S}^{1}\times\mathbb{S}^{1}$,
~$\theta=\upsilon^{t}|_{\mathbb{R}\times M}$ ~and ~$K=\{p\}\times\mathbb{S}^{1}\times\mathbb{S}^{1}\subsetneq M$
~a\emph{ compact aperiodic minimal}.

Endow the compact manifold ~$\mathcal{N}=\mathbb{S}^{3}\times\mathbb{S}^{1}\times\mathbb{S}^{1}$
~with the ~$C^{\,\infty}$ ~vector field ~$\upsilon=\zeta+\frac{\partial}{\partial\rho}+\sqrt{2}\frac{\partial}{\partial\sigma}=\nu+\sqrt{2}\frac{\partial}{\partial\sigma}$
~where ~$\rho$ ~and ~$\sigma$ ~are, respectively, the coordinates
of the left and of the right ~$\mathbb{S}^{1}$. ~Since ~$\{p\}\times\mathbb{S}^{1}$
~is a \emph{periodic orbit} of ~$\nu\in\mathfrak{X}^{\,\infty}\big(\mathbb{S}^{3}\times\mathbb{S}^{1}\big)$~
with constant speed 1, ~the 2-torus ~$\{p\}\times\mathbb{S}^{1}\times\mathbb{S}^{1}\subset\mathcal{N}$
~is invariant under ~$\upsilon^{t}$ ~and carries an \emph{irrational}
linear flow with slope ~$\sqrt{2}$. ~Again as in example 2, the
flow ~$\big(\mathbb{S}^{3}\times\mathbb{S}^{1}\times\mathbb{S}^{1},\upsilon^{t}\big)$
~projects to ~$\big(\mathbb{S}^{3}\times\mathbb{S}^{1},\nu^{t}\big)$
~via the projection ~$\pi_{_{1,2}}$ ~onto the two first factors,
and to ~$\big(\mathbb{S}^{1},(\sqrt{2}\partial/\partial\sigma)^{t}\big)$
~via the projection ~$\pi_{_{3}}$ ~onto the 3rd factor, hence
for every ~$(u,s,w)\in\mathcal{N}$ ~and ~$t\in\mathbb{R},$ \begin{equation}
\pi_{_{1,2}}\circ\upsilon^{t}(u,s,w)=\nu^{t}\big(\pi_{_{1,2}}(u,s,w)\big)=\nu^{t}(u,s)\end{equation}
\begin{equation}
\upsilon^{t}(u,s,w)=\big(\nu^{t}(u,s),\,\sqrt{2}t+w\,(\mbox{mod\,2\ensuremath{\pi}})\big)\end{equation}

By an argument entirely analogous to that employed in example 2, it
follows that taking ~$z\in\varphi(\Upsilon_{_{X}})\in\mathbb{S}^{3}$
~(see example 1), and letting ~$M:=\mbox{cl\,}\mathcal{O}_{_{\upsilon}}(z,b,c)\subset\mbox{cl\,}\mathcal{O}_{_{\nu}}(z,b)\times\mathbb{S}^{1}$,
~$\theta:=\upsilon^{t}|_{\mathbb{R}\times M}$ ~and ~$K:=\{p\}\times\mathbb{S}^{1}\times\mathbb{S}^{1}$
~, where ~$b,\, c$ ~are arbitrary points of the first (left) and
second (right) ~$\mathbb{S}^{1}$, respectively, we get a subsmooth
independent realization\emph{ }of condition~ \textbf{6},\textbf{
~}with ~$K$ ~a \emph{compact aperiodic minimal }(see above), the
ambient manifold\emph{ }being ~$\mathcal{N}=\mathbb{S}^{3}\times\mathbb{S}^{1}\times\mathbb{S}^{1}$.\medskip{}

~~~We will now briefly indicate how to obtain a subsmooth independent
realization of condition \textbf{~5} ~with ~$K$ ~an \emph{equilibrium
orbit}. ~For ~$K$ ~a \emph{periodic orbit }or a \emph{compact
aperiodic minimal }we only have to proceed as in examples 2 and 3.
~The analogue subsmooth independent realizations of condition \textbf{~4~
}are obtained time-reversing the previously mentioned realizations
of condition \textbf{~5}.

A subsmooth independent realization of condition \textbf{~5~ }with
\textbf{~$K$ }~an \emph{equilibrium orbit} and ambient manifold
~$\mathbb{S}^{3}$ ~is achieved through a simple (and obvious) modification
of Beniere and Meigniez's construction: in their paper (\cite{beni},
p.23, bottom), the cut-and-paste operation that defines ~$M_{_{1}}$
~is performed \emph{only} for each ~$p\in\mathbb{Z}_{_{0}}^{+}$
~(instead of for all ~$p\in\mathbb{Z}\,\big)$. ~Following otherwise
their construction, we finally obtain a ~$C^{\,\infty}$ orientable
surface ~$\mathfrak{N}$ ~of infinite genus, again with its end
points set homeomorphic to ~$\varDelta:=\{0\}\,\cup\,\{1/n:\, n\in\mathbb{N}\}\subset\mathbb{R}$
~and all ends \emph{flat }except the non-isolated one. ~By \emph{Kerekjarto}
\emph{Theorem} (see \emph{e.g.} \cite{rich}, p.262 and \cite{beni},
p.26), this surface is thus ~$C^{\,\infty}$ ~diffeomorphic to the
surface ~$\mathfrak{M}$~ of example 1 (\emph{i.e.} to the original
surface carrying a flow without minimal set constructed in \cite{beni})
and hence to ~$S\subset\mathbb{R}^{3}$. ~$\mathfrak{N}$ ~carries
a ~$C^{\,\infty}$ ~vector field ~$\upsilon$ ~whose flow is\emph{
no longer without minimal sets, }as there are points ~$x\in\mathfrak{N}$
~with ~$\alpha(x)=\emptyset=\omega(x)$. ~However for each \emph{limit
point }~$x\in\mathfrak{N}$ ~of the flow ~$\upsilon^{t}$, ~it
is easily seen that ~$\alpha(x)=\emptyset$, ~$\omega(x)\neq\emptyset$
~and both the positive and negative orbits of ~$x$ ~accumulate
in the unique non-isolated end of ~$\mathfrak{N}$ ~and in no other
end of this surface. ~As in example 1, we have again a ~$C^{\,\infty}$
~diffeomorphism ~$f:\mathfrak{N}\longrightarrow S\subset\mathbb{R}^{3}$~
defining a ~$C^{\,\infty}$ ~embedding ~$\mathfrak{N}\hookrightarrow\mathbb{R}^{3}$~
and inducing a complete tangent vector field ~$X:=f_{_{*}}\upsilon$
~on ~$S\subset\mathbb{R}^{3}$. ~Then for each \emph{limit point
~}$x$ ~of ~$\big(S,X^{t}\big)$, ~$\omega_{_{X}}(x)$ ~is a
~\emph{closed, unbounded} subset of ~$\mathbb{R}^{3}$ ~and ~$\underset{t\rightarrow-\infty}{\mbox{lim}}\|X^{t}(x)\|=+\infty$~
$\big(\,\|\centerdot\|$ ~being the euclidean norm on ~$\mathbb{R}^{3}\,\big)$,
~\emph{i.e.} the point ~$x$ ~\emph{{}``escapes''} to infinite
on ~$\mathbb{R}^{3}$ ~when ~$t\longrightarrow-\infty$. ~Then
proceeding exactly as in example 1, ~for each ~$q\in\varUpsilon_{_{X}}$,
~$z:=\varphi(q)$,\[
\begin{array}{c}
\alpha_{_{\zeta}}(z)=\{p\}\mbox{\,\,\,\,\,\,\,\,\ and\,\,\,\,\,\,\,\,\,}\omega_{_{\zeta}}(z)=\varphi\big(\omega_{_{X}}(q)\big)\,\sqcup\,\{p\}\\
\,\\
\mbox{cl\,}\mathcal{O}_{_{\zeta}}(z)=\varphi\big(\mbox{cl\,}\mathcal{O}_{_{X}}(q)\big)\,\sqcup\,\{p\}\end{array}\]

Letting ~$M:=\mbox{cl\,}\mathcal{O}_{_{\zeta}}(z)\subset\mathbb{S}^{3}$,
~~$\theta:=\zeta^{t}|_{\mathbb{R}\times M}$ ~and ~$K:=\{p\}$,
~we then have\[
y\in M\setminus K\implies y\in B_{_{\theta}}^{-}(K)\,\cap\, A_{_{\theta}}^{+}(K)\]
It is now immediate to verify that for these ~$M$,~ $\theta$,
~$K$, ~none of the 28 conditions of Theorem 1\emph{, }with the
exception of condition ~\textbf{5},\emph{ ~}can hold and therefore
~$(M,\theta)$, ~$K$ ~provide a subsmooth independent realization
of condition ~\textbf{5} (see example 1). 

We cannot enter here the details of the construction of subsmooth
independent realization for the nine exceptional condition ~\textbf{7.4
~}to ~\textbf{9.6}. ~In \cite{teix} such realizations are given
for conditions ~\textbf{7.4~} to ~\textbf{7.6} ~with ambient manifold
~$\mathbb{R}^{4}$, ~for ~\textbf{8.4} ~to ~\textbf{8.6} ~with
ambient manifold ~$\mathbb{R}^{5}$ ~and for ~\textbf{9.4} ~to
~\textbf{9.6} ~with ambient manifold ~$\mathbb{R}^{6}$, ~the
procedure being the same in all cases. ~The crucial fact is that
~$\mathbb{S}^{3}$, ~$\mathbb{S}^{3}\times\mathbb{S}^{1}$, ~$\mathbb{S}^{3}\times\mathbb{S}^{1}\times\mathbb{S}^{1}$
~embed, respectively, as codimension one, closed ~$C^{\,\infty}$
~submanifolds of ~$\mathbb{R}^{4}$, ~$\mathbb{R}^{5}$ ~and~
$\mathbb{R}^{6}$ ~$\big($for the last two the \emph{open book decompositions
}of ~$\mathbb{R}^{5}$ ~and ~$\mathbb{R}^{6}$ ~are used), and
consequently have \emph{trivial} normal tubular neighbourhoods in
the corresponding ambient manifolds. \bigskip{}

\noun{8. PROOF OF LEMMA 7. ~TWO TOPOLOGICAL LEMMAS.}

\bigskip{}

\textbf{Lemma 7.} \emph{~Let ~$M$ ~be a locally compact metric
space} \emph{with a ~$C^{\,0}$} ~\emph{flow.}~\noun{ }\emph{If
~$A\subset M$~ is open and ~}$\mbox{CMin}$\emph{$(D)$ ~is ~$d_{_{H}}-$dense
in itself ~then ~}$\mbox{CMin}(D)$\emph{~ is ~$\mathfrak{c}-$dense
in itself. ~~If ~$\mathfrak{D}$ ~is a ~$d_{_{H}}-$open and
dense in itself subset of ~}$\mbox{CMin}(M)$\emph{ ~then ~$\mathfrak{D}$
~is ~$\mathfrak{c}-$dense in itself.}\medskip{}

~~~Before entering the proof, a few technical definitions will
be needed.~ Recall that $\mbox{F}:=\{0,1\}$.~ For each ~$n\in\mathbb{N},\,\,\,\,\,\,\mbox{F}^{\, n}=\{0,1\}^{n}=$~the
set of finite sequences of 0's and 1's with length $n.$ ~Again,
since no risk of ambiguity arises, commas and brackets are omitted
in the representation of the elements of ~$\mathcal{F},$~ thus
we write, for example, ~$01$~ instead of ~$(0,1)$~ and ~$\mbox{F}^{1}$
~is naturally identified with ~$\mbox{F}$.

~ Let\[
\mbox{F}^{\,\leq n}:=\underset{1\leq m\leq n}{\bigsqcup}\mbox{F}^{\, m}\,\,\,\,\,\,\,\,\,\,\,\,\,\,\,\,\,\,\,\,\,\,\,\,\,\,\,\,\,\mathcal{F}:=\underset{n\in\mathbb{N}}{\bigsqcup}\,\mbox{F}^{\, n}\,\,\,\,\,\,\,\,\,\,\,\,\,\,\,\,\,\,\,\,\,\,\,\,\,\,\,\,\,\,\mathcal{F}_{_{\emptyset}}:=\mathcal{F}\,\sqcup\,\{\emptyset\}\]
$\big(\,\mathcal{F}_{_{\emptyset}}$~is the set of finite sequences
of 0's and 1's including the empty sequence ~$\emptyset\,\big)$.~
If ~$a,b\in\mathcal{F}_{_{\emptyset}},$ ~$ab$ ~represents, as
usual, the element of ~$\mathcal{F}_{_{\emptyset}}$ ~obtained by
adjoining ~$b$~ to the right end of ~$a$~ $\big($naturally,
~$a\emptyset=a=\emptyset a$ ~for every ~$a\in\mathcal{F}_{_{\emptyset}}\,\big)$.
~For any ~$n\in\mathbb{N}\,\,\,\,\mbox{and}\,\,\, a\in\mbox{F}^{\, n},\,\,\,\,|a|:=n$
~(the \emph{length} of ~$a).$ ~We now define the operators ~$*,+,\bullet,-$
~on ~$\mathcal{F}$:

$\,\,\,\,\,\,\,\,\,\,\,\,\,\,\,\,0^{*}:=1,\,\,\,\,\,\,\,\,\,1^{*}:=0,\,\,\,\,\,\,\,\,\,(bc)^{*}:=bc^{*}$

$\,\,\,\,\,\,\,\,\,\,\,\,\,\,\,$$(bc)^{+}:=bcc^{*},$ ~~~~$(bc)^{\bullet}:=bcc$
~~~~for any ~$b\in\mathcal{F}_{_{\emptyset}},\,\,\,$$c\in\mbox{F,}$

$\,\,\,\,\,\,\,\,\,\,\,\,\,\,\,\,0^{-}:=0,\,\,\,\,\,\,\,\,\,1^{-}:=0$
~~and ~~$(bc)^{-}:=b$ ~~~~for every ~$b\in\mathcal{F},\,\, c\in\mbox{F.}$

Obvious facts about ~$\mathcal{F}$ ~and its operators such as ~
$\mbox{F}^{\, n+1}=\{ac:\,\, a\in\mbox{F}^{\, n},\,\, c\in\mbox{F}\}=(\mbox{F}^{\, n})^{+}\,\cup\,(\mbox{F}^{\, n})^{\bullet}$
~will be implicitly used without mention.\medskip{}

\emph{Proof of Lemma 7. ~}The first part of the lemma clearly implies
the second: suppose ~$\mathfrak{D}$ ~is a ~$d_{_{H}}-$open and
dense in itself subset of ~$\mbox{CMin}(M).$ ~By lemma 5 (section
5), for each ~$X\in\mathfrak{D}$, ~there is a ~$\delta_{_{X}}>0$
~such that ~$\mbox{CMin}\big(B(\, X,\,\delta_{_{X}})\big)\subset\mathfrak{D}$
~(as ~$\mathfrak{D}$ ~is ~$d_{_{H}}-$open ~in ~$\mbox{CMin}(M)$,
~there is an ~$\epsilon_{_{X}}>0$ ~such that ~$B_{_{H}}(\, X,\,\epsilon_{_{X}})\,\cap\,\mbox{CMin}(M)\subset\mathfrak{D}$.
~By lemma 5, there is a ~$\delta_{_{X}}>0$ ~such that ~$\mbox{CMin}\big(B(\, X,\,\delta_{_{X}})\big)\subset\mbox{Ci}\big(B(\, X,\,\delta_{_{X}})\big)\subset B_{_{H}}(\, X,\,\epsilon_{_{X}})$,
~hence ~$\mbox{CMin}\big(B(\, X,\,\delta_{_{X}})\big)\subset\mathfrak{D}\,\big).$
~By lemma 8, ~$\mbox{CMin}\big(B(\, X,\,\delta_{_{X}})\big)$ ~is
a ~$d_{_{H}}-$open subset of ~$\mbox{CMin}(M)$ ~and hence of
~$\mathfrak{D}.$ ~$\mathfrak{D}$ ~is ~$d_{_{H}}-$dense in itself
(by hypothesis), and so is ~$\mbox{CMin}\big(B(\, X,\,\delta_{_{X}})\big)$,
~since it is $d_{_{H}}-$open in ~$\mathfrak{D}.$ ~By the first
part of lemma 7, ~$\mbox{CMin}\big(B(\, X,\,\delta_{_{X}})\big)$
~is actually ~$\mathfrak{c}-$\emph{dense in itself}, and obviously
so is ~$\mathfrak{D}=\underset{X\in\mathfrak{D}}{\bigcup}\mbox{CMin}\big(B(\, X,\,\delta_{_{X}})\big)$.

\emph{Proof of the 1st part of lemma 7.~ }Since ~\emph{$A$ ~}is
open and minimal sets are compact, given any ~$\varLambda_{_{0}}\in\mbox{CMin}(A)$,
~$B(\varLambda_{_{0}},\,\epsilon_{_{0}})\subset A$ ~for sufficiently
small ~$\epsilon_{_{0}}>0,$ ~hence by Lemma 6 (section 5) it is
sufficient to prove that for any ~$\varLambda_{_{0}}\in\mbox{CMin}(A)$
~and ~$\epsilon_{_{0}}>0,$\[
\#\,\mbox{CMin}\left(B(\,\varLambda_{_{0}},\,\epsilon_{_{0}})\right)=\mathfrak{c}\]

~~~The demonstration is based on a generalisation of the idea lying
behind the construction of Cantor's ternary set: for each ~$\varLambda\in\mbox{CMin}(A)$~~and
~$\epsilon>0$ ~such that ~$B[\,\varLambda,\,\epsilon\,\,]$ ~is
a compact subset of ~$A$, ~the existence of a continuum of $d_{_{H}}-$Cauchy
sequences of compact minimal sets ~$X\in\mbox{CMin}\big(B(\,\varLambda,\,\epsilon\,)\big)$
~$d_{_{H}}-$ converging to a continuum of \emph{mutually disjoint,
}nonvoid, compact, connected invariant sets ~$Y\subset B(\,\varLambda,\,\epsilon\,)$~
is proved. ~The result then follows since each such ~\emph{Y} ~contains
at least one compact minimal set of the flow.

~~~Let ~$\varLambda_{_{0}}\in\mbox{CMin}(A)$ ~and ~$\epsilon_{_{0}}>0$
~such that ~$B[\,\varLambda_{_{0}},\epsilon_{_{0}}\,]\subset A$
~is compact. ~Since ~$\mbox{CMin}$\emph{$(A)$} ~is $d_{_{H}}$-dense
in itself, there is a ~$\varLambda_{_{1}}\in\mbox{CMin}(A)$ ~such
that

\[
\varLambda_{_{1}}\in B_{_{H}}(\varLambda_{_{0}},\epsilon_{_{0}})\setminus\{\varLambda_{_{0}}\}\]

This implies ~$\varLambda_{_{1}}\subset B(\varLambda_{_{0}},\epsilon_{_{0}})$
~and, ~$\big|\varLambda_{_{1}}\big|_{\varLambda_{_{0}}}<\epsilon_{_{0}}$
~since ~$\varLambda_{_{1}}$ ~is closed. ~Also ~$\varLambda_{_{0}}\,\cap\,\varLambda_{_{1}}=\emptyset$
~since distinct minimal sets are disjoint.~~Let

\[
\epsilon_{_{1}}:=\mbox{min}\left\{ \mbox{dist}(\varLambda_{_{0}},\varLambda_{_{1}})\,\,,\,\,\epsilon_{_{0}}-\big|\varLambda_{_{1}}\big|_{\varLambda_{_{0}}}\right\} \Big/3\]

Clearly, $\epsilon_{_{1}}>0,$~~~$B[\,\varLambda_{_{0}},\epsilon_{_{1}}\,]\,\cap\, B[\,\varLambda_{_{1}},\epsilon_{_{1}}\,]=\emptyset$~~
and by the triangle inequality for the metric~ $d$~ of ~$M,$~
it follows%
\footnote{by the same argument as that given in footnote {[}50{]} ahead, changing
$n$ , $b$ to $0$ and $n+1$ , $b^{+}$ to $1$.%
} that ~~$B[\,\varLambda_{_{1}},\epsilon_{_{1}}\,]\subset B(\,\varLambda_{_{0}},\epsilon_{_{0}}\,)$.
~Also ~$0<\epsilon_{_{1}}<\epsilon_{_{0}}/3$ ~hence ~$B[\,\varLambda_{_{0}},\epsilon_{_{1}}\,]\subset B(\,\varLambda_{_{0}},\epsilon_{_{0}}\,)$.

~~~Consider now the following proposition ~$\Theta(n)$ ~in the
variable ~$n\in\mathbf{\mathbb{N}}:$\vspace{2mm}

~~~~$\varLambda_{_{a}}\in\mbox{CMin}(A)$ ~~and ~~$\epsilon_{_{m}}>0$
~are defined for all ~$\,\, a\in\mbox{F}^{\,\leq n}$ ~and ~$\,\,1\leq m\leq n$,
~respectively, and satisfy:\smallskip{}

$\bullet\,\,\,\,$for all ~$1\leq m\leq n$,\smallskip{}

~~~~~a)~~~~~$\epsilon_{_{m}}<\epsilon_{_{m-1}}\big/3$ \smallskip{}

$\bullet\,\,\,\,$for all ~$a\in\mbox{F}^{\,\leq n}$,\smallskip{}

~~~~~b)~~~~~$B\big[\,\varLambda_{_{a}},\epsilon_{_{|a|}}\,\big]\,\,\cap\,\, B\big[\,\varLambda_{_{a^{*}}},\epsilon_{_{|a|}}\,\big]=\emptyset$\smallskip{}

~~~~~c)~~~~~$B\big[\,\varLambda_{_{a}},\epsilon_{_{|a|}}\,\big]\,\,\subset\,\, B\big(\varLambda_{_{a^{-}}},\epsilon_{_{|a|-1}}\big)$

\smallskip{}

Note that the truth of $\Theta(1)$ is already established. We prove
that $\Theta(n)\implies\Theta(n+1)$. ~Assume $\Theta(n)$ is true.
~Observe that it is enough%
\footnote{since all the other conditions are guaranteed to be true by the hypothesis
that ~$\Theta(n)$ ~is true.%
} to define ~$\epsilon_{_{n+1}}>0$~ and $\,\varLambda_{_{a}}\in\mbox{CMin}(A)\mbox{\,\ for all }$$a\in\mbox{F}^{\, n+1}$~
and show that ~$\epsilon_{_{n+1}}<\epsilon_{_{n}}\big/3$~ and that
\[
B\left[\,\varLambda_{_{a}},\epsilon_{_{n+1}}\,\right]\,\,\cap\,\, B\left[\,\varLambda_{_{a^{*}}},\epsilon_{_{n+1}}\,\right]=\emptyset\,\,\,\,\,\,\,\,\mbox{and\,\,\,\,\,\,\,\,}B\left[\,\varLambda_{_{a}},\epsilon_{_{n+1}}\,\right]\,\,\subset\,\, B\big(\varLambda_{_{a^{-}}},\epsilon_{_{n}}\,\big)\]

for all ~$a\in\mbox{F}^{\, n+1}$. ~~Define ~$\varLambda_{_{b^{\bullet}}}:=\varLambda_{_{b}}$
~for every ~$b\in\mbox{F}^{\, n}$. ~~Again since ~$\mbox{CMin}(A)$
~is $d_{_{H}}$-dense in itself, ~$\varLambda_{_{b}}\in\mbox{CMin}(A)$~
and ~$\epsilon_{_{n}}>0,$ ~we may select ~$\varLambda_{_{b^{+}}}\in\mbox{CMin}(A)$
~satisfying \begin{equation}
\varLambda_{_{b^{+}}}\in B_{_{H}}(\varLambda_{_{b}},\epsilon_{_{n}})\setminus\{\varLambda_{_{b}}\}\end{equation}

With each ~$a\in\mbox{F}^{\, n+1}$ ~we have thus associated a compact
minimal set contained in ~$A$.~ Now let ~$\epsilon_{_{n+1}}:=\mbox{min}(\delta_{_{n+1}},\lambda_{_{n+1}})$
~where

\[
\delta_{_{n+1}}:=\mbox{min}\left\{ \mbox{dist\ensuremath{\left(\varLambda_{_{b}},\varLambda_{_{b^{+}}}\right)}}:\, b\in\mbox{F}^{\, n}\right\} \Big/3\,\]
\[
\lambda_{_{n+1}}:=\mbox{min}\left\{ \epsilon_{_{n}}-\big|\varLambda_{_{b^{+}}}\big|_{\varLambda_{_{b}}}:\, b\in\mbox{F}^{\, n}\right\} \Big/3\]

By (14), ~$\varLambda_{_{b}}$~and ~$\varLambda_{_{b^{+}}}$ ~are
distinct (compact) minimal sets hence disjoint, thus ~$\mbox{dist\ensuremath{\big(\varLambda_{_{b}}},\ensuremath{\varLambda_{_{b^{+}}}\big)}}>0\,;$
~~$\varLambda_{_{b^{+}}}\subset B\left(\varLambda_{_{b}},\epsilon_{_{n}}\right)$
~and ~$\varLambda_{_{b^{+}}}$~is closed hence ~$\big|\varLambda_{_{b^{+}}}\big|_{\varLambda_{_{b}}}<\epsilon_{_{n}}$,
~therefore ~$\epsilon_{_{n+1}}>0.$ ~Also, ~$\epsilon_{_{n+1}}<\epsilon_{_{n}}/3$
~since ~$\epsilon_{_{n+1}}\leq\lambda_{_{n+1}}<\epsilon_{_{n}}/3$.
~Observe that for any ~$a\in\mbox{F}^{\, n+1},\,\,\,\,\,\,\,\{a,a^{*}\}=\{b^{\bullet},b^{+}\}$~
where ~$b:=a^{-}\in\mbox{F}^{\, n},$~ thus ~$\{\varLambda_{_{a}},\varLambda_{_{a^{*}}}\}=\{\varLambda_{_{b}},\varLambda_{_{b^{+}}}\}$
~since ~$\varLambda_{_{b^{\bullet}}}=\varLambda_{_{b}}$; ~therefore
\[
B\left[\,\varLambda_{_{a}},\epsilon_{_{n+1}}\,\right]\,\,\cap\,\, B\left[\,\varLambda_{_{a^{*}}},\epsilon_{_{n+1}}\,\right]=\emptyset\]
since ~$\epsilon_{_{n+1}}\leq\delta_{_{n+1}}.$ ~To complete the
induction it remains to show that \begin{equation}
B\left[\,\varLambda_{_{a}},\epsilon_{_{n+1}}\,\right]\,\,\subset\,\, B\big(\varLambda_{_{a^{-}}},\epsilon_{_{n}}\big)\,\,\,\,\mbox{for all \,}\, a\in\mbox{F}^{\, n+1}\end{equation}
Let ~$b:=a^{-}$; ~if ~$a=b^{\bullet}$~then ~ $\varLambda_{_{a}}=\varLambda_{_{b^{\bullet}}}=\varLambda_{_{b}}$
~thus the inclusion is trivial since ~$\epsilon_{_{n+1}}<\epsilon_{_{n}}\big/3<\epsilon_{_{n}};$~
if ~$a=b^{+}$ ~then%
\footnote{by the triangle inequality for the the metric ~$d$ ~of ~$M$,
~~$\varLambda_{_{b^{+}}}\subset B\left[\varLambda_{_{b}},\,\big|\varLambda_{_{b^{+}}}\big|_{\varLambda_{_{b}}}\right]$
~hence ~$B\Big[\varLambda_{_{b^{+}}},\epsilon_{_{n+1}}\Big]\subset B\left[\varLambda_{_{b}},\,\big|\varLambda_{_{b^{+}}}\big|_{\varLambda_{_{b}}}+\epsilon_{_{n+1}}\right]$
~and ~since ~$\big|\varLambda_{_{b^{+}}}\big|_{\varLambda_{_{b}}}<\epsilon_{_{n}},$
~it follows that ~~$\big|\varLambda_{_{b^{+}}}\big|_{\varLambda_{_{b}}}+\epsilon_{_{n+1}}\leq\big|\varLambda_{_{b^{+}}}\big|_{\varLambda_{_{b}}}+\left(\epsilon_{_{n}}-\big|\varLambda_{_{b^{+}}}\big|_{\varLambda_{_{b}}}\right)\Big/3<\epsilon_{_{n}}$,
~thus ~$B\Big[\varLambda_{_{b^{+}}},\epsilon_{_{n+1}}\Big]\subset B\left[\varLambda_{_{b}},\,\big|\varLambda_{_{b^{+}}}\big|_{\varLambda_{_{b}}}+\epsilon_{_{n+1}}\right]\subset B\left(\varLambda_{_{b}},\epsilon_{_{n}}\right).$%
} ~$B\big[\,\varLambda_{_{b^{+}}},\epsilon_{_{n+1}}\big]\subset B\big(\varLambda_{_{b}},\epsilon_{_{n}}\big)$.~
The proof that ~$\Theta(n)\implies$$\Theta(n+1)$ ~is complete.
~Therefore ~$\,\varLambda_{_{a}}\in\mbox{CMin}(A)$ ~and ~$\epsilon_{_{n}}>0$
~are defined for all ~$a\in$$\mathcal{F}$ ~and ~$n\in\mathbf{\mathbb{N},}$
~respectively and are such that condition a) of $\Theta(n)$ is true
for all ~$m\in\mathbb{N}$\textbf{ ~}and conditions b) and c) are
true for all ~$a\in\mathcal{F}.$ ~Using induction again, it is
now simple to deduce from c) that\[
B\big[\varLambda_{_{a}},\epsilon_{_{|a|}}\,\big]\,\subset\,\, B\big(\varLambda_{_{0}},\epsilon_{_{0}}\big)\,\,\,\,\mbox{for all \,}\, a\in\mathcal{F}\,\]
\[
B\big[\,\varLambda_{_{ab}},\epsilon_{_{|ab|}}\big]\,\,\subset\,\, B\big(\varLambda_{_{a}},\epsilon_{_{|a|}}\big)\,\,\,\,\mbox{for all \,}\, a,b\in\mathcal{F}\,\]

~~~To conclude the proof we associate with each infinite sequence
of ~0's ~and ~1's, ~$w\in\mbox{F}^{\,\mathbb{N}}$, ~a compact
minimal set ~$\varGamma_{_{w}}\in\mbox{CMin}\left(B(\,\varLambda{}_{_{0}},\epsilon_{_{0}}\,)\right)$
~and show that for all ~~$w,v\in\mbox{F}^{\,\mathbb{N}}$$,\mbox{ \,\,\,\,\,}w\neq v$
~implies ~$\varGamma_{_{w}}\neq\varGamma_{_{v}}$, ~therefore proving
the existence of a \emph{continuum} of compact minimal sets contained
in ~$B(\,\varLambda{}_{_{0}},\epsilon_{_{0}}\,),$~ since ~$\#\,\mbox{F}^{\,\mathbb{N}}=\mathfrak{c}$.
~For each ~$w\in\mbox{F}^{\,\mathbb{N}}$ ~let ~$w_{_{n}}\in\mbox{F}^{\, n}$
~denote the sequence of the first \emph{n} digits of ~\emph{w. ~}We
claim that ~$\big(\varLambda_{_{w_{_{n}}}}\big)$ ~is a $d_{_{H}}-$Cauchy
sequence in the compact metric space ~$\mathfrak{S}\big(B[\,\varLambda_{_{0}},\,\epsilon_{_{0}}]\big)$.
~Observe that together, the triangle inequality for the ~$d_{_{H}}$~
metric of ~$\mbox{C}(M)$, ~~$\varLambda_{_{bc}}\in B_{_{H}}\big(\varLambda_{_{b}},\epsilon_{_{|b|}}\big)$
~for all ~$b\in\mathcal{F},\, c\in\mbox{F}$ ~and ~$\epsilon_{_{n+1}}<\epsilon_{_{n}}/3$
~~imply that for all ~$n,k\in\mathbb{N}$,\[
\varLambda_{_{w_{_{n+k}}}}\in B_{_{H}}\left(\varLambda_{_{w_{_{n}}}},\epsilon_{_{n}}\cdot\left(\sum_{j=0}^{k-1}\,\frac{1}{3^{^{j}}}\right)\right)\]
and since ~$\epsilon_{_{n}}<\epsilon_{_{0}}/3^{n}$ ~it follows
that for any ~$p,q>n$\[
\varLambda_{_{w_{_{p}}}},\,\varLambda_{_{w_{_{q}}}}\in B_{_{H}}\left(\varLambda_{_{w_{_{n}}}},\frac{\epsilon_{_{0}}}{2\cdot3^{n-1}}\right)\]
hence \[
d_{_{H}}\Big(\varLambda_{_{w_{_{p}}}},\,\varLambda_{_{w_{_{q}}}}\Big)<2\cdot\frac{\epsilon_{_{0}}}{2\cdot3^{n-1}}=\frac{\epsilon_{_{0}}}{3^{n-1}}\]

therefore the sequence ~$\big(\varLambda_{_{w_{_{n}}}}\big)$~ is
clearly $d_{_{H}}-$Cauchy and thus $d_{_{H}}$-convergent. ~Since
for all ~$n\in\mathbb{N}$\textbf{,\[
\varLambda_{_{w_{_{n}}}}\subset B\big[\,\varLambda_{_{w_{_{1}}}},\epsilon_{_{1}}\big]\subset B\big(\varLambda_{_{0}},\epsilon_{_{0}}\big)\]
}it follows that\[
\mbox{lim}\,\varLambda_{_{w_{_{n}}}}=:\varLambda_{_{w}}\in\mathfrak{S}\left(B(\varLambda_{_{0}},\epsilon_{_{0}})\right)\]

On the other hand, if ~$v,w\in\mbox{F}^{\,\mathbb{N}}$~and ~$v\neq w$~
then there is a~ $m\in\mathbb{N}$~ such that ~$v_{_{m}}=w_{_{m}}^{*}$,
~thus\[
\begin{array}{c}
\Big(\varLambda_{_{w{}_{_{n}}}}\in B\big(\varLambda_{_{w_{_{m}}}},\epsilon_{_{m}}\big)\,\,\,\,\,\mbox{for any}\,\,\,\,\,\, n>m\Big)\implies\varLambda_{_{w}}\subset B\big[\,\varLambda_{_{w_{_{m}}}},\epsilon_{_{m}}\,\big]\\
\,\\
\Big(\varLambda_{_{v{}_{_{n}}}}\in B\big(\varLambda_{_{v_{_{m}}}},\epsilon_{_{m}}\big)\,\,\,\,\,\mbox{for any}\,\,\,\,\,\, n>m\Big)\implies\varLambda_{_{v}}\subset B\big[\,\varLambda_{_{v_{_{m}}}},\epsilon_{_{m}}\,\big]=B\big[\,\varLambda_{_{w_{_{m}}^{*}}},\epsilon_{_{m}}\,\big]\end{array}\]

therefore ~$\varLambda_{_{w}}\,\,\cap\,\,\varLambda_{_{v}}=\emptyset$
~since ~$B\big[\,\varLambda_{_{w_{_{m}}}},\epsilon_{_{m}}\,\big]\,\,\cap\,\, B\big[\,\varLambda_{_{w_{_{m}}^{*}}},\epsilon_{_{m}}\,\big]=\emptyset.$
~This establishes the existence of a \emph{continuum} of mutually
disjoint, nonvoid, compact, connected, invariant sets ~$\varLambda_{_{w}}\in\mathfrak{S}\left(B(\varLambda_{_{0}},\epsilon_{_{0}})\right),$~
$w\in\mbox{F}^{\,\mathbb{N}}$ ~and each such ~$\varLambda_{_{w}}$
contains necessarily at least one compact minimal set of the flow
~$\varGamma_{_{w}}\in\mbox{CMin}\left(B(\,\varLambda{}_{_{0}},\epsilon_{_{0}}\,)\right)\subset\mbox{CMin}(A)$.%
\footnote{\emph{Remark: }Observe that even assuming the \emph{Continuum Hypothesis},
a standard Baire category argument is not enough, without further
assumptions, to prove Lemma 7. ~Suppose that ~$\varLambda_{_{0}}\in\mbox{CMin}(A)$
~and ~$\epsilon_{_{0}}>0$ ~are such that ~$B[\,\varLambda_{_{0}},\epsilon_{_{0}}]$
~is a compact subset of ~$A$. ~Since ~$\mbox{CMin}(A)$~ is
$d_{_{H}}-$dense in itself, ~$\mathfrak{I}:=\mbox{cl}_{_{H}}\,\mbox{CMin}\left(B(\,\varLambda{}_{_{0}},\epsilon_{_{0}}\,)\right)$\emph{~}
is easily seen to be a $d_{_{H}}-$dense in itself and compact subset
of ~$\mathfrak{S}\left(B[\,\varLambda_{_{0}},\epsilon_{_{0}}]\right).$
~By Baire Theorem, if ~$\mathfrak{I}$ ~is nonvoid then ~$\mathfrak{I}$
cannot be countable, therefore assuming the \emph{Continuum Hypothesis},
~$\#\,\mathfrak{I}\geq\mathfrak{c}$, ~i.e. ~$\mathfrak{I}$ ~contains
a continuum of \emph{distinct }compact, connected, invariant subsets
of ~$B[\,\varLambda_{_{0}},\epsilon_{_{0}}].$~ However \emph{a
priori }nothing guarantees that these sets are \emph{mutually disjoint},
thus the possibility of associating with each ~$X\in\mathfrak{I}$
~a \emph{distinct} (and hence disjoint) compact minimal set contained
in ~$X$ ~is compromised. %
}\hfill{}$\blacksquare$\medskip{}

\textbf{8.2}\textbf{\emph{ ~}}Recall that every locally compact,
connected mertic space is separable.

\textbf{Proposition. ~}\emph{If ~$M$ ~is a locally compact, separable
metric space then so is ~}$\big[\,\mbox{C}(M),\, d_{_{H}}\big]$\emph{,
~where ~}$\mbox{C}(M)$\emph{ ~is the set of all nonvoid compact
subsets of ~$M$.}

\emph{Proof.} ~Claim 1:~$\big[\,\mbox{C}(M),\, d_{_{H}}\big]$ ~\emph{is
locally compact.} ~Since ~$M$ ~is locally compact, given any ~$K\in\mbox{C}(M)$
~there is an ~$\epsilon>0$~ such that ~$B[\, K,\,\epsilon\,]$
~is compact. ~By \emph{Blaschske Theorem }(see section 2) ~$\mbox{C}\big(B[\, K,\,\epsilon\,]\big)$
~is ~$d_{_{H}}-$compact. ~Since~ $X\in B_{_{H}}(K,\,\epsilon)\,\Longrightarrow\, X\subset B(K,\,\epsilon)\,\Longrightarrow\, X\in\mbox{C}\big(B[\, K,\,\epsilon\,]\big)$,
~$\mbox{C}\big(B[\, K,\,\epsilon\,]\big)$ ~is a ~$d_{_{H}}-$compact
neighbourhood of ~$K$ ~in ~$\big[\,\mbox{C}(M),\, d_{_{H}}\big]$,
~thus ~this last metric space is locally compact.

Claim 2: ~$\big[\mbox{\,\ C}(M),\, d_{_{H}}\big]$ ~\emph{is separable}.
~Since ~$M$ ~is locally compact and separable, it admits an equivalent
\emph{boundedly} \emph{compact} \emph{metric }~$d'$ ~(every closed
and bounded subset is compact, see \emph{e.g. }\cite{lima}, p.278).
~Fixing ~$O\in M$, ~$M=\underset{n\geq1}{\bigcup}B_{_{d'}}[\, O,\, n\,]$
~defines an exhaustion of ~$M$ ~by compact subsets. ~Again by
\emph{Blaschke Theorem, ~$\mbox{C}\big(B_{_{d'}}[\, O,\, n\,]\big)$
~}is ~$d'_{_{H}}-$compact for every ~$n\geq1$ ~and hence ~$d'_{_{H}}-$separable,
thus ~$\big[\mbox{\,\ C}(M),\, d'_{_{H}}\big]$ ~is separable since
~$\mbox{C}(M)=\underset{n\geq1}{\bigcup}\mbox{C}\big(B_{_{d'}}[\, O,\, n\,]\big)$.
~Now ~$d_{_{H}}$ ~and ~$d'_{_{H}}$ ~are equivalent metrics
on ~$\mbox{C}(M)$ ~since ~$d$ ~and ~$d'$ ~are equivalent
on ~$M$, ~hence ~$\big[\mbox{\,\ C}(M),\, d{}_{_{H}}\big]$ ~is
separable.\hfill{}$\blacksquare$

\medskip{}

\textbf{8.3}\textbf{\emph{ ~}}\emph{Proof of}\textbf{\emph{ ~$CH\Longleftrightarrow\mathfrak{c}DH$}}\emph{:}

$\big(CH\Longrightarrow\mathfrak{c}DH\,\big)$\emph{:} ~recall that
a separable metric space ~$L$ ~has at most a \emph{continuum} of
points. ~Suppose ~$L$ ~is \emph{uncountable}. ~Let ~$\mathfrak{I}$
~be the set of points of ~$L$ ~for which there is an $\epsilon>0$
~such that ~$B(x,\,\epsilon)$ ~is \emph{countable}. ~For each
~$x\in\mathfrak{I}$ ~define\[
\epsilon_{_{x}}:=\mbox{sup}\big\{\epsilon>0:\,\, B(\, x,\,\epsilon)\mbox{ \,\ is countable\,\ensuremath{\big\}}}\]
Now observe that in the proof of Theorem 3, from the beginning until
the end of claim 1 we have only used the fact that ~$\big[\mbox{\,\ CMin}(M),\, d_{_{H}}\big]$
~is an uncountable separable metric, hence all facts proved until
there are valid for arbitrary uncountable separable metric spaces.
~In particular ~$\mathfrak{I}$ ~is a countable, open subset of
~$L$ ~and ~$\mathfrak{D}:=L\setminus\mathfrak{I}$ ~is a set
such that every neighbourhood ~$U_{_{z}}$ ~of each ~$z\in\mathfrak{D}$
~contains uncountably many points of ~$\mathfrak{D}.$ ~As the
cardinal of ~$L\supset\mathfrak{D}$ ~is at most ~$\mathfrak{c}=2^{\aleph_{_{0}}}$,
~the \emph{Continuum Hypothesis }actually implies ~$\#\,\big(U_{_{z}}\,\cap\,\mathfrak{D}\big)=\mathfrak{c}$.
~Hence ~$\mathfrak{D}$ ~is ~$\mathfrak{c}-$\emph{dense in itself.}

$\big(\,\neg\, CH\Longrightarrow\neg\,\mathfrak{c}DH\,\big)$\emph{:
~}assume there is a cardinal ~$\aleph_{_{0}}<\beta<\mathfrak{c}.$
~Given the bijection between the cardinal ~$\mathfrak{c}\supset\beta$
~and ~$\mathbb{R}$, ~there is a set ~$S\subset\mathbb{R}$ ~with
~$\aleph_{_{0}}<\#\, S=\beta<\mathfrak{\mathfrak{c}}$. ~With the
euclidean metric inherited from ~$\mathbb{R},$~ $S$ ~is an uncountable
separable metric space. ~Removing from ~$S$ ~an arbitrary \emph{countable}
set ~$\mathfrak{I}$ ~we obtain again a set ~$\mathfrak{D}=S\setminus\mathfrak{I}$
~of cardinal ~$\#\, S=\#\,(S\setminus\mathfrak{I})=\beta$. ~Hence
~$\aleph_{_{0}}<\#\,\mathfrak{D}=\beta<\mathfrak{c}$ ~and therefore,
since ~$\mathfrak{D}$ ~is nonvoid, it cannot be ~$\mathfrak{c}-$\emph{dense
in itself}.\hfill{}$\blacksquare$

\_\_\_\_\_\_\_\_\_\_\_\_\_\_\_\_\_\_\_\_\_\_\_\_\_\_\_\_\_\_\_\_\_\_\_\_\_\_\_\_\_\_\_\_\_\_\_\_\_\_\_\_\_\_\_\_\_\_\_\_\_\_\_\medskip{}

\emph{Acknowledgement: }the author wishes to thank Jorge Rocha (FCUP
- Portugal) and Vítor Araújo (UFB - Brasil), for their helpful comments
and suggestions, leading to an improvement of the presentation. \medskip{}

{\small \medskip{}
}{\small \par}

\emph{E-mail address: pedro.teixeira@fc.up.pt}

Centro de Matemática da Universidade do Porto \\
Rua do Campo Alegre, 687 \\
4169-007 Porto \\
Portugal
\end{document}